\documentclass[11pt]{amsart}

\calclayout
\makeatletter
\def\l@subsection{\@tocline{2}{0pt}{2.5pc}{5pc}{}}
\def\l@subsubsection{\@tocline{2}{0pt}{5pc}{7.5pc}{}}
\makeatother

\usepackage[T1]{fontenc}
\usepackage[utf8]{inputenc}

\usepackage{graphicx}

\usepackage[english]{babel}
\usepackage[
backend=biber,
style=alphabetic,
sorting=nyt,
giveninits=true,
maxbibnames=99,
isbn=false,
doi=true,
url=true,
eprint=true
]{biblatex}

\DeclareNameAlias{default}{given-family}

\DeclareFieldFormat[article,incollection,inproceedings]
{title}{\mkbibemph{#1}}
\DeclareFieldFormat{journaltitle}{#1}
\DeclareFieldFormat[article]{volume}{\mkbibbold{#1}}
\DeclareFieldFormat[article]{pages}{#1}

\DeclareDelimFormat{nametitledelim}{\addcomma\space}

\renewbibmacro*{in:}{%
	\ifentrytype{article}
	{}
	{\printtext{\bibstring{in}\intitlepunct}}}

\DeclareFieldFormat{labelalphawidth}{#1}

\AtEveryBibitem{%
	\iffieldundef{doi}
	{}
	{\clearfield{url}}}

\usepackage{amsthm}
\usepackage{amssymb}
\usepackage{dirtytalk}
\usepackage[all]{xy}
\usepackage[]{enumitem}
\usepackage[dvipsnames,table]{xcolor}
\DeclareMathAlphabet{\pazocal}{OMS}{zplm}{m}{n}
\usepackage{calrsfs}
\usepackage{faktor}
\usepackage{tikz-cd}
\usepackage[
colorlinks=true,
linkcolor=blue,
citecolor=purple,
urlcolor=blue
]{hyperref}
\usepackage[nameinlink]{cleveref}
\usepackage{array,booktabs,tabularx,pdflscape}
\usepackage{caption}
\usepackage{geometry}

\newtheorem{thm}{Theorem}
\newtheorem*{thm*}{Theorem}
\Crefname{thm}{Theorem}{Theorems}
\newtheorem{thmmain}{Theorem}

\crefname{thmmain}{theorem}{theorems}
\Crefname{thmmain}{Theorem}{Theorems}
\newtheorem{prop}[thm]{Proposition}
\newtheorem*{prop*}{Proposition}
\newtheorem{lem}[thm]{Lemma}
\newtheorem*{lem*}{Lemma}
\newtheorem{cor}[thm]{Corollary}
{Main Corollary}
\newtheorem{dfn}[thm]{Definition}
\newtheorem*{dfn*}{Definition}

\newtheorem{rem}[thm]{Remark}
\newtheorem*{rem*}{Remark}

\newtheorem{ex}
{Example}

\newcommand{\Z}{\mathbb{Z}}
\newcommand{\Q}{\mathbb{Q}}
\newcommand{\R}{\mathbb{R}}
\newcommand{\C}{\mathbb{C}}

\renewcommand{\i}{\mathrm{i}}
\newcommand{\pr}{\operatorname{pr}}

\renewcommand{\S}{\mathbb{S}}
\newcommand{\RP}{\mathbb{RP}}

\renewcommand{\P}{\mathbb{P}}
\renewcommand{\O}{\pazocal{O}}
\newcommand{\Pic}{\operatorname{Pic}}
\newcommand{\Bl}{\operatorname{Bl}}
\newcommand{\NN}{\pazocal{N}}

\makeatletter
\newcommand{\leftfaktor}[2]{%
	\mathchoice
	{\leftfaktor@impl\displaystyle{#1}{#2}}%
	{\leftfaktor@impl\textstyle{#1}{#2}}%
	{\leftfaktor@impl\scriptstyle{#1}{#2}}%
	{\leftfaktor@impl\scriptscriptstyle{#1}{#2}}%
}
\newcommand{\leftfaktor@impl}[3]{%
	\settoheight{\dimen@}{$#1#2$}%
	\raisebox{-0.4\dimen@}{$#1#2$}
	\mkern-4mu\diagdown\mkern-4mu
	\settoheight{\dimen@}{$#1#3$}%
	\raisebox{0.4\dimen@}{$#1#3$}
}
\makeatother

\begin{document}
	
\title[]{Topology of some smooth real Fano threefolds arising from blow-ups along rational curves.}
\author{Guillaume Kineider}
\address{Aix-Marseille University, CNRS, I2M, Marseille, France}
\email{guillaume.kineider@univ-amu.fr}

\thanks{}

\begin{abstract}
	We describe the topological blow-up of a fibre of a Seifert fibration. We then apply this result to determine the diffeomorphism types of the real loci of smooth real Fano threefolds obtained by blowing up either $\P^3_\R$ or a smooth real quadric threefold along a smooth real rational curve. As an application, we also determine the topology of the real loci of the two real forms of the quintic del Pezzo threefold $V_5$, completing the topological classification of smooth real rational del Pezzo threefolds. 
\end{abstract}

\maketitle
\tableofcontents

\section{Introduction.}

Let $X$ be a smooth real projective variety. Its real locus $X(\R)$, endowed with the Euclidean topology, is a smooth closed real manifold whenever it is nonempty. Certain topological and geometric properties of this manifold reflect subtle algebro-geometric features of the underlying variety. For instance, Comessatti proved in 1914 \cite{comessattiSullaConnessioneSuperficie1914} that, in dimension two, an orientable hyperbolic real locus forms an obstruction for rationality, even geometrically. In a series of papers \cite{kollarRealAlgebraicThreefolds1998, kollarRealAlgebraicThreefolds1999, kollarRealAlgebraicThreefolds1999a, kollarRealAlgebraicThreefolds2016}, Koll\'ar aimed, among other things, at extending the result of Comessatti to threefolds. Using the Minimal Model Program, he showed that the topology of the real locus of a smooth rationally connected threefold $X$ can be studied, provided it is orientable or satisfies certain weaker assumptions, through the topology of a suitable birational model \cite[Theorem~1.2]{kollarRealAlgebraicThreefolds1998}. Such a minimal model must be a Mori fibre space and is therefore either a conic bundle, a del Pezzo fibration, or a Fano threefold. The topology of the first two classes was treated in \cite{kollarRealAlgebraicThreefolds1999a} and \cite{kollarRealAlgebraicThreefolds2016}, leaving the finitely many families of Fano threefolds as the remaining building blocks. Some families have already been studied, for instance real cubic threefolds and intersection of two quadrics by Krasnov (\cite{krasnovRigidIsotopyClassification2006, krasnovIntersectionsTwoReal2018}), and orientable prime Fano threefolds of genus $12$ by Koll\'ar and Schreyer \cite{KS04}. For most of the other families, the question remains open.

The purpose of this article is first to study this question for every Fano blow-ups of $\P_\R^3$ and quadric threefolds along smooth rational curves, except for quadric threefolds of signature $(3,2)$ blown-up along quartics. This latter class will be studied in a forthcoming paper. The following results summary the main results in this direction.

\begin{thmmain}[Reduced version of \Cref{Th:BlowUpProjSpaceAlongRatCurves}]
	Let $c$ be a smooth $\R$-rational curve in $\P^3$ of degree $d$. Then,
	\begin{enumerate}
		\item[$(d = 1):$] $\Bl_c \P^3(\R)\approx \RP^2\times \S^1$;
		\item[$(d = 2):$] $\Bl_c \P^3(\R)\approx \RP^3\#\, \S^2\widetilde{\times}\S^1$;
		\item[$(d = 3):$] $\Bl_c \P^3(\R)\approx \RP^2\times \S^1$;
		\item[$(d = 4):$] $\Bl_c \P^3(\R)\approx \RP^3\#\, \S^2\widetilde{\times}\S^1$ or $M\big(\RP^2; 0 ; \overline{0};(2,1),(4,1)\big)$.
	\end{enumerate}
\end{thmmain}

The manifold $\S^2\widetilde{\times}\S^1$ is described in \Cref{ex:CircleBundles} and the notation $M\big(\dots\big)$ represents the diffeomorphism type of the Seifert fibred bundle described by the given invariants. The details are set in \Cref{Sec:SeifertFibrations}.

\begin{thmmain}[Reduced version of \Cref{Th:BlowUpQuadricThreefold32AlongRatCurves} and \Cref{Th:BlowUpQuadricThreefold41AlongRatCurves}]
	Let $Q$ be a smooth real quadric hypersurface of $\P^4$ and let $c$ be a smooth real rational curve in $Q$ with degree $d$. If $Q$ has signature $(3,2)$, the real locus of $\Bl_c Q$ has the following description:
	\begin{enumerate}
		\item[$(d = 1):$] $\Bl_cQ(\R)\approx \mathbb{RP}^2\times \S^1$;
		\item[$(d = 2):$] $\Bl_c Q(\R)\approx \S^2\widetilde{\times}\S^1 \#\, \S^2\widetilde{\times}\S^1$ or $M\big(\mathrm{Kl}^2; 0; \bar{1}\big)$ or $\S^2\widetilde{\times}\S^1$,
		\item[$(d = 3):$] $\Bl_cQ(\R)\approx \mathbb{RP}^2\times \S^1$ or $M\big(\overline{\mathbb D^2};0;0;(2,1),(3,1)\big)$ where $\overline{\mathbb D^2}$ denotes the closed $2$-disc whose boundary is a reflector circle.
	\end{enumerate} 
	If $Q$ has signature $(4,1)$, the real locus of $\Bl_c Q$ has the following description:
	\begin{enumerate}
		\item[$(d = 2):$] $\Bl_{c} Q (\R) \approx \S^2\widetilde{\times}\S^1$;
		\item[$(d = 4):$] $\Bl_{c} Q (\R) \approx \S^2\widetilde{\times}\S^1$.
	\end{enumerate}
\end{thmmain}

The interest of these computations goes beyond the individual examples. Koll\'ar's work suggests that birational geometry may be a powerful tool to study the topology of real algebraic varieties, especially for Fano varieties since their birational geometry is highly constrained. Thus, gaining a better comprehension of the elementary birational operations appearing in the Minimal Model Program or the Sarkisov Program should be highly rewarding. An example of this is that the degree $3$ case of \Cref{Th:BlowUpQuadricThreefold32AlongRatCurves} allows us to determine the topology of the trivial real form of the quintic del Pezzo threefold $V_5$. Its other real form is naturally a homogeneous space whose topology is clear: this is \Cref{Th:TopologyOfV5} we reproduce below. This result completes the computation of the real loci of smooth real rational del Pezzo threefolds.

\begin{thmmain}[Reduced version of \Cref{Th:TopologyOfV5}]
	Let $X_5$ and $Y_5$ be the real forms of $V_5$. Then
	\[X_5(\R)\approx\RP^3 \qquad \text{whereas} \qquad Y_5(\R) \approx \faktor{\operatorname{SO}_3(\R)}{O}.\]
	Here $O\simeq\mathfrak S_4$ is the rotational octahedral group.
\end{thmmain}
 
The diffeomorphism type of the blow-up along a subvariety $V\subset X$ is specified by the embedding $V(\R)\hookrightarrow X(\R)$ up to smooth isotopy. Then, the topological operation induced by the blow-up, which shall be called \emph{topological blow-up} is some kind of surgery we recall in \Cref{Sec:ConstructionTopBlowUp}. The challenge is therefore to understand the global picture after this surgery. Our strategy is to treat the specific case where $V(\R)$ appears as a fibre  of some fibration, since heuristically the local picture should mix smoothly with the global one. The most general statement for this principle is \Cref{prop:TopBlowUpAlongFibreInLocTrivialBundle}: blowing up a fibre in a locally trivial bundle produces the pullback bundle over the blow-up of the base.

Our specific case of interest is curves in threefolds. Thus our framework relies on circle bundle structures which admit a generalisation, namely Seifert fibre bundles which play in addition a significant and natural role in the classification of $3$-manifolds. We thus extend the above generic result to Seifert fibrations, this is \Cref{Th:TopBlowUpSeifertFibrationAlongFibre}. This is the main topological tool used throughout the paper.

Determining the diffeomorphism type of the blow-up along a curve $c\subset X$ requires understanding the smooth isotopy type of $c(\R)$ inside $X(\R)$. The classification of such links is a problem in its own right. For example, the classification of real algebraic links in $\P^3$ has has been carried out in \cite{bjorklundRealAlgebraicKnots2011} for rational curves of degree at most $5$ and later extended to curves of genus $g\leq 1$ and degree at most $6$ in \cite{mikhalkinRealAlgebraicKnots2016}. The case of smooth rational curves in $\P^3$ is quite simple, as they lie on a hyperboloid and are therefore analogues of torus knots in $\S^3$. We shall call them \emph{hyperboloidal} knots. It is known that hyperboloidal knots, with one exception, may be realised as regular fibres of Seifert bundle structures on $\RP^3$ so we may apply \Cref{Th:TopBlowUpSeifertFibrationAlongFibre}. It should be emphasised that our method does not restrict to Fano blow-ups as it computes the diffeomorphism type, or at least a small list of possible ones, for the blow-up of $\P^3$ along some smooth curves of any degree, namely those whose real loci are isotopic to a hyperboloidal knot. For the quadric threefolds, it is a case-by-case study, using a similar line of idea. This article is a first step towards a systematic study of the topology of real loci of rational Fano threefolds. We intend to extend this investigation to further families in forthcoming work.

There is also a geometric motivation from Thurston's geometrization. Koll\'ar proved that, apart from finitely many possible exceptions arising from Fano threefolds, closed orientable hyperbolic $3$-manifolds do not occur as real loci of smooth rationally connected threefolds. These remaining exceptions were actually ruled out by Viterbo, using different methods from symplectic field theory (see \cite{KhBourbaki}). This also showed that the orientability assumption may be dropped. To go further in that direction, the conclusion of the series of papers from Koll\'ar yields the following. 
\begin{thm}[Theorem~6.6 in \cite{kollarTopologyRealComplex}]
	Let $X$ be a smooth, real, projective threefold. Assume that $X$ is uniruled and $X(\R)$ is orientable. Then every component of $X(\R)$ is among the following:
	\begin{enumerate}
		\item Seifert fibred;
		\item connected sum of lens spaces
		\item torus bundle over $\S^1$, or double covered by a torus bundle over $\S^1$;
		\item finitely many other possible exceptions, or
		\item obtained from the above by repeatedly taking connected sum with $\RP^3$ and $\S^2\times \S^1$.
	\end{enumerate}
\end{thm}
Again, the fourth item corresponds to Fano varieties. The third possibility has been ruled out in \cite{mangolteUniruledSixmanifoldsContain2012}. Thus, one expects the real loci of rational Fano threefolds to carry relatively restricted three-dimensional geometries. The descriptions given in our main results gives readily the JSJ decomposition and the geometries (see \Cref{Sec:SeifertFibrations}). The tables below summarise \Cref{Th:BlowUpProjSpaceAlongRatCurves,Th:BlowUpQuadricThreefold32AlongRatCurves,Th:BlowUpQuadricThreefold41AlongRatCurves}.

\smallskip

\textit{Acknowledgements.}
I would like to thank Frédéric Mangolte and Andrea Fanelli for their supervision and guidance throughout this project. I am also grateful to Adrien Dubouloz and Jules Chenal for their interest in my work and their valuable insights. Finally, I thank the topologists in Marseille---Benjamin Audoux, Luisa Paoluzzi, Léo Bénard and Michel Boileau---for their helpful discussions and advice.

\smallskip
\textit{Use of generative AI.}
GPT-6 Astra (non-Pro version) was used to assist with the redaction of this article.  
\newpage

\newgeometry{margin=1.7cm}
\begin{landscape}	
	\begin{table}[p]
		\centering
		\small
		\setlength{\tabcolsep}{4pt}
		\setlength{\parskip}{0pt}
		\setlength{\doublerulesep}{2pt}
		\renewcommand{\arraystretch}{1.25}
		
		\caption{Diffeomorphism type and geometries of the real loci. When, for one degree, there are several classes of curves, we use normal forms introduced in the text. Otherwise, we use the generic notation $c_d$ for any smooth $\R$-rational curve of degree $d$ in $Z$. The first rows concern $Z(\R)$ and the two last rows of each block concern $(\Bl_C Z)(\R)$. In this table $M_0:=M\big(\overline{\mathbb D^2};0;0;(2,1),(3,1)\big)$.}
		\label{tab:topology-blowups}
		
		\begin{tabularx}{\linewidth}{
				|p{2.7cm}|
				*{5}{>{\centering\arraybackslash}X|}
			}
			\hline
			Threefold $Z$
			& \multicolumn{5}{c|}{$\P^3$}\\
			\hline
			
			Diffeomorphism type
			& \multicolumn{5}{c|}{$\RP^3$}\\[0.5em]
			\hline
			
			Thurston geometries
			& \multicolumn{5}{c|}{$\S^3$}\\[0.5em]
			\hline
			
			Degree of $C$
			& $1$ & $2$ & $3$
			& \multicolumn{2}{c|}{$4$}\\
			\hline
			
			Curve $C$
			& $c_1$ & $c_2$ & $c_3$
			& $\Gamma_{4,1}$ & $\Gamma_{4,3}$\\
			\hline
			
			Diffeomorphism type
			& $\RP^2\times\S^1$
			& $\RP^3\#\,\S^2\widetilde{\times}\S^1$
			& $\RP^2\times\S^1$
			& $\RP^3\#\,\S^2\widetilde{\times}\S^1$
			& $M\big(\RP^2;0;\overline{0};(2,1),(4,1)\big)$\\[0.8em]
			\hline
			
			Thurston geometries
			& $\S^2\times\R$
			& $\S^3,\;\S^2\times\R$
			& $\S^2\times\R$
			& $\S^3,\;\S^2\times\R$
			& $\mathbb{H}^2\times\R$\\[0.8em]
			\hline\hline
		\end{tabularx}\par\nointerlineskip
		\begin{tabularx}{\linewidth}{
				|p{2.7cm}|
				*{8}{>{\centering\arraybackslash}X|}
			}
			Threefold $Z$
			& \multicolumn{8}{c|}{$\pazocal Q_{3,2}$}\\
			\hline
			
			Diffeomorphism type
			& \multicolumn{8}{c|}{$\S^2\widetilde{\times}\S^1$}\\[0.5em]
			\hline
			
			Thurston geometries
			& \multicolumn{8}{c|}{$\S^2\times\R$}\\[0.5em]
			\hline
			
			Degree of $C$
			& $1$
			& \multicolumn{2}{c|}{$2$}
			& \multicolumn{4}{c|}{$3$}
			& $4$\\
			\hline
			
			Curve $C$
			& $c_1$
			& $C_{2,+}$ & $C_{2,-}$
			& $C_{3,-2}$ & $C_{3,0}$
			& $C_{3,2}$ & $C_{3,s}$
			& $c_4$\\
			\hline
			
			Diffeomorphism type
			& $\RP^2\times\S^1$
			& $\S^2\widetilde{\times}\S^1 \#\,\S^2\widetilde{\times}\S^1$
			& $M\big(\mathrm{Kl}^2;0;\bar{1}\big)$
			& $\RP^2\times\S^1$
			& $\RP^2\times\S^1$
			& $M_0$
			& $\RP^2\times\S^1$
			& Unknown\\[0.8em]
			\hline
			
			Thurston geometries
			& $\S^2\times\R$
			& $\S^2\times\R,\;\S^2\times\R$
			& $\mathbb{E}^3$
			& $\S^2\times\R$
			& $\S^2\times\R$
			& $\mathbb{H}^2\times\R$
			& $\S^2\times\R$
			& Unknown\\[0.8em]
			\hline\hline
		\end{tabularx}\par\nointerlineskip
		\begin{tabularx}{\linewidth}{
				|p{2.7cm}|
				*{10}{>{\centering\arraybackslash}X|}
			}
			Threefold $Z$
			& \multicolumn{2}{c|}{$\pazocal Q_{4,1}$}
			& \multicolumn{5}{c|}{$X_5$}
			& \multicolumn{3}{c|}{$Y_5$}\\
			\hline
			
			Diffeomorphism type
			& \multicolumn{2}{c|}{$\S^3$}
			& \multicolumn{5}{c|}{$\RP^3$}
			& \multicolumn{3}{c|}{
				$M\big(\S^2;0;-1;(2,1);(3,1);(4,1)\big)$
			}\\[0.5em]
			\hline
			
			Thurston geometries
			& \multicolumn{2}{c|}{$\S^3$}
			& \multicolumn{5}{c|}{$\S^3$}
			& \multicolumn{3}{c|}{$\S^3$}\\[0.5em]
			\hline
			
			Degree of $C$
			& $2$ & $4$
			& \multicolumn{3}{c|}{$1$} & $2$ & $3$
			& $1$ & $2$ & $3$\\
			\hline
			
			Curve $C$
			& $c_2$ & $c_4$
			& $\ell_{\mathrm{int}}$ & $\ell_{\mathrm{ext}}$ & $\ell_{\mathrm{sp}}$
			& $c_2$ & $c_3$
			& $c_1$ & $c_2$ & $c_3$\\
			\hline
			
			Diffeomorphism type
			& $\S^2\widetilde{\times}\S^1$
			& $\S^2\widetilde{\times}\S^1$
			& $\RP^2\times\S^1$
			& $\RP^2\times\S^1$
			& $\RP^2\times\S^1$
			& Unknown
			& Unknown
			& $M_0$
			& Unknown
			& Unknown\\[0.8em]
			\hline
			
			Thurston geometries
			& $\S^2\times\R$
			& $\S^2\times\R$
			& $\S^2\times\R$
			& $\S^2\times\R$
			& $\S^2\times\R^1$
			& Unknown
			& Unknown
			& $\mathbb{H}^2\times\R$
			& Unknown
			& Unknown\\[0.8em]
			\hline
		\end{tabularx}
	\end{table}
\end{landscape}

\restoregeometry
\newpage 

\section{Seifert fibrations and their invariants.}\label{Sec:SeifertFibrations}

Throughout, we use \emph{Seifert fibration} in the sense of \cite[\textsection~3]{scottGeometries3Manifolds1983}, allowing non-orientable total spaces. Since conventions vary, we record those used below and refer to Scott's survey for a comprehensive introduction to the general theory. For us, such a fibration is a decomposition into circles, each having a saturated neighbourhood isomorphic to a standard fibred solid torus or to a fibred solid Klein bottle. The latter is the mapping torus of a reflection of a disc, equipped with its suspension foliation. The quotient map
\[ p\colon M\longrightarrow B \]
has a two-dimensional orbifold as its base. For the closed connected manifolds considered here, its singularities are cone points and possibly reflector circles. Reflector circles form the boundary of the underlying surface $|B|$, although $M$ itself has no boundary.

This broader definition is useful here: most of our manifolds are non-orientable, arising from topological blow-ups along knots with trivial normal bundle, and reflector circles occur naturally. Our use of the Seifert fibrations is twofold. First, they provide explicit descriptions of the prime manifolds encountered below. These are particularly efficient to determine their model geometries; see below. They also represent our main tool to describe globally the effect of blowing up along knots; see \Cref{Th:TopBlowUpSeifertFibrationAlongFibre}. Let us also mention that Epstein's theorem \cite{epsteinPeriodicFlowsThreeManifolds1972} asserts that every smooth foliation by circles of a compact $3$-manifold is a Seifert fibration in this sense.

We describe fibrations using 
\[ \big(B;w_1;e;(\alpha_1,\beta_1),\ldots,(\alpha_k,\beta_k)\big). \]
\begin{enumerate}
	\item $B$ is the base orbifold, with underlying surface $|B|$. Cone orders may be omitted from the notation for $B$ when they are already recorded by the $\alpha_i$; reflector circles are always indicated.
	\item The class $w_1\in H^1(|B|;\Z_2)$ records reversal of the orientation of the fibre under transport around loops. It is the Stiefel--Whitney class of the fibration, not of its total space. In the absence of reflectors,
	\[ M\text{ is orientable} \quad\Longleftrightarrow\quad w_1=w_1(T|B|). \]
	A reflector circle always makes $M$ non-orientable.	
	\item Each pair $(\alpha_i,\beta_i)$ describes an isolated exceptional fibre, with $\alpha_i\geq 2$ and $\gcd(\alpha_i,\beta_i)=1$. Our oriented local model is the quotient of $\mathbb D^2\times\S^1$ by the cyclic action generated by
	\[ (z,u)\longmapsto (e^{2\pi i/\alpha_i}z,e^{-2\pi i\beta_i/\alpha_i}u). \]
	Thus Scott's orbit invariants $(\alpha_i,q_i)$ satisfy $\beta_iq_i\equiv1\pmod{\alpha_i}$. For oriented $M$, local orientations are induced by that of $M$, and we normalise $0<\beta_i<\alpha_i$. For non-orientable $M$, local orientations may be reversed independently, and we normalise $0<\beta_i\leq\alpha_i/2$.
	\item Once these local data are normalised, the Euler number $e$ records the remaining global gluing obstruction. If $B$ has no reflectors,
	\[ e\in	\begin{cases}
		\Z, & M\text{ orientable},\\
		\Z_2, & M\text{ non-orientable}.
	\end{cases} \]
	In the non-orientable case, we take $e=0$ if some $\alpha_i=2$: the two possible values then describe isomorphic fibrations. If $B$ has reflector circles, we also set $e=0$, since there is no independent global obstruction parameter.
\end{enumerate}
Thus, our conventions follow Scott, except that we write $e=-b$ for his global obstruction invariant; $e_{\mathbb Q}$ below denotes his rational Euler number. When $M$ is orientable, we choose an orientation to compute the invariants, although the notation for its total space forgets this orientation. The notation
\[ M\big(B;w_1;e;(\alpha_1,\beta_1),\ldots,(\alpha_k,\beta_k)\big) \]
denotes the diffeomorphism type of the total space, without a chosen orientation. We use it when no simpler description is available, and omit the empty list for genuine circle bundles. Different fibrations can have diffeomorphic total spaces. The classification and the extension across reflector circles are discussed in \cite[\textsection~3, pp.~434--437]{scottGeometries3Manifolds1983}.

\begin{rem}
	For a genuine circle bundle over a closed connected surface $B$, there is an Euler class
	\[\mathcal E(p)\in H^2(B;\Z_{w_1}), \qquad H^2(B;\Z_{w_1})\simeq
	\begin{cases}
		\Z, & w_1=w_1(TB),\\
		\Z_2, & w_1\neq w_1(TB).
	\end{cases}\]
	Here $\Z_{w_1}$ is the local coefficient system with monodromy $(-1)^{w_1(\gamma)}$ along a loop $\gamma$. The displayed isomorphisms follow from Poincar\'e duality with local coefficients. Under these identifications, $\mathcal E(p)$ corresponds to the Euler number $e$; in the integer-valued case, its sign is fixed by the chosen orientation of $M$. In particular, $e=0$ exactly when the bundle admits a section. Like the Stiefel--Whitney class, the Euler class is natural under pullback.
	
	To make this naturality explicit for Seifert fibrations, notably in the proof of \Cref{Th:TopBlowUpSeifertFibrationAlongFibre}, we use a cohomological description of their invariants.	For bases with only cone-point singularities, the integral orbifold Euler class
	\[ \mathcal E(p)\in H^2_{\mathrm{orb}}(B;\Z_{w_1}) \]
	records the global twisting together with the exceptional-fibre data. Here orbifold cohomology means the cohomology of the orbifold's classifying space. For an orbifold map $f\colon B'\to B$, naturality takes the form
	\[ \mathcal E(f^*p)=f^*\mathcal E(p) \in H^2_{\mathrm{orb}}(B';\Z_{f^*w_1}). \]
	Our survey \cite{kineiderClassificationSeifertFibrations2026} gives the construction of this class, its functorial properties, and its relation to Seifert invariants.
\end{rem}

We now recall the relation with Thurston geometries. A closed connected Seifert manifold always admits a (unique) geometry model which is determined by
\[\chi_{\mathrm{orb}}(B) = \chi(|B|)-\sum_{i=1}^k \left(1-\frac{1}{\alpha_i}\right),\qquad e_{\mathbb Q}(p)= \begin{cases}
	e-\displaystyle\sum_{i=1}^k\frac{\beta_i}{\alpha_i}, & M\text{ oriented},\\
	0, & M\text{ non-orientable},
\end{cases}\]
according to the following table:
\[\begin{array}{c|ccc}
	& \chi_{\mathrm{orb}}(B)>0 & \chi_{\mathrm{orb}}(B)=0 & \chi_{\mathrm{orb}}(B)<0\\
	\hline
	e_{\Q}=0	& \S^2\times\R & \mathbb E^3 & \mathbb H^2\times\R\\
	e_{\Q}\neq 0 & \S^3 & \mathrm{Nil} & \widetilde{\mathrm{SL}_2(\R)}
\end{array}\]
Conversely, the closed manifolds admitting one of these six geometries are precisely the Seifert manifolds; see \cite[Theorem~5.3]{scottGeometries3Manifolds1983}. Notice that $e_{\mathbb Q}=0$ for a non-orientable total space does not imply $e=0$.

\begin{ex}\label{ex:CircleBundles}
	Let $a$ denote the nonzero element of $H^1(\RP^2;\Z_2)$.
	\begin{itemize}
		\item The quotient
		\[\faktor{(\S^2\times\S^1)}{\bigl((x,z)\sim(-x,\overline z)\bigr)} \longrightarrow \RP^2, \qquad [x,z]\longmapsto[x], \]
		is a circle bundle with $w_1=a$ and a section $[x]\mapsto[x,1]$, hence $e=0$. Its total space is $\RP^3\#\,\RP^3$: a fundamental domain is $\S^2\times[0,\pi]$, with antipodal identification on each end, and cutting at the middle gives two copies of $\RP^3$ minus an open ball. Thus
		\[ \RP^3\#\,\RP^3\approx M\big(\RP^2;a;0\big).\]
		Note that it is the only non-trivial connected sum admitting a geometry, modelled on $\S^2\times \R$.
		\item Throughout the article, we write
		\[\S^2\widetilde{\times}\S^1 :=	\faktor{(\S^2\times[0,1])}{\bigl((x,1)\sim(-x,0)\bigr)}.\]
		This mapping torus is non-orientable and prime, but not irreducible. The map $[x,t]\mapsto[x]$ makes it a genuine circle bundle over $\RP^2$. Its fibre orientation is preserved, so $w_1=0$. Its Euler class is the nonzero element of $H^2(\RP^2;\Z)\simeq\Z_2$, since the zero class would give $\RP^2\times\S^1$, whereas
		\[ \pi_1(\S^2\widetilde{\times}\S^1)\simeq\Z, \qquad \pi_1(\RP^2\times\S^1)\simeq\Z_2\times\Z. \]
		Consequently,
		\[ \S^2\widetilde{\times}\S^1\approx M\big(\RP^2;0;1\big).\]
		It also carries a geometry modelled on $\S^2\times\R$.
	\end{itemize}
\end{ex}
    
\section{Topological blow-up.}

Heuristically, blowing up consists in separating undirected normal directions along a subvariety. When the blown-up subvariety is a real, smooth and has real points, it induces a surgery on the real locus we recall in the section below. Our presentation follows \cite{weiWuClassesTopological2011}; see also \cite{mikhalkinBlowupEquivalenceSmooth1997} for a comprehensive exposition. Throughout the rest of the paper, the unqualified term \say{blow-up} will refer to the \emph{algebraic} blow-up ; when we mean the construction below we will write \emph{topological} blow-up.
 
\subsection{Construction.}\label{Sec:ConstructionTopBlowUp}

Let $M$ be a smooth closed manifold and $N$ a closed submanifold of codimension $r>1$. Fix a tubular embedding $\tau:\pazocal N_{N/M}\hookrightarrow M$ identifying $\pazocal N_{N/M}$ with an open tubular neighbourhood of $N$. Consider its real projectivisation $p:\P\NN_{N/M}\to N$, together with the tautological bundle $\pazocal{L}$ over it defined by :
\[\pazocal{L}=\left\{\big((x,\ell),(x,v)\big)\in p^*\NN_{N/M}\subset \P\NN_{N/M}\times_N \NN_{N/M}\;| \; v\in\ell\right\}.\]
\[\xymatrix{\pazocal{L}\ar@{^{(}->}[r] \ar[rdd] & p^*\NN_{N/M} \ar[dd] \ar[rd]& & \\
	& & \NN_{N/M} \ar[dd] \ar@{^{(}->}[rdd]^\tau& \\
	& \P\NN_{N/M} \ar[dr]^p & & \\
	& & N \ar@{^{(}->}[r] & M\\}\]    

\begin{dfn}
	The \emph{topological blow-up of $M$ along $N$} is the smooth manifold
	\[\Bl^{\mathrm{top}}_N M = (M- N)\cup_{\phi} \pazocal L,\]
	with
	\begin{equation}\label{eq:GluingTopBlowUp}
		\begin{array}{c | c c c}
			\phi :& \pazocal{L}- \P\NN_{N/M}&\to &\tau(\NN_{N/M}- N)\subset M- N\\
			& ((x,\ell),(x,v))&\mapsto & \tau(x,v)
		\end{array}.
	\end{equation}
\end{dfn}

The map $p$ extends to the blow-down map $p:\Bl^{\mathrm{top}}_N M\to M$, equal to the identity on $M\setminus N$, and given on $\pazocal L$ by
\[((x,\ell),(x,v))\longmapsto \tau(x,v).\]
The subset
\[E:=p^{-1}(N)=\P\pazocal N_{N/M}.\]
is a smooth hypersurface in $\Bl^{\mathrm{top}}_N M$, called the \emph{exceptional divisor}, whose normal bundle is
\[\pazocal N_{E/\Bl^{\mathrm{top}}_N M}\simeq \pazocal L.\]
The diffeomorphism type of $\Bl_N^{\mathrm{top}}M$ is invariant under smooth isotopy $N\hookrightarrow M$. For more details, see \cite{mikhalkinBlowupEquivalenceSmooth1997}. A well-known first example is the following. 

\begin{prop}\label{prop:TopBlowUpInAPoint}
	Let $M$ be a smooth connected $n$-manifold and $x$ a point in $M$. Then,
	\[\Bl^{\mathrm{top}}_x M\approx M\#\, \RP^n.\] 
	If moreover $M$ is oriented and $n$ is odd, then $\Bl_x^{\mathrm{top}}M$ is naturally oriented, and the above diffeomorphism may be chosen orientation-preserving with $\RP^n$ is endowed with the orientation induced by the chosen orientation of $T_xM$.
\end{prop}
\begin{proof}
	By definition, the topological blow-up of $M$ at $x$ is obtained by	replacing a small open ball around $x$ by the tautological line bundle over
	\[ \P\NN_{x/M}\simeq\P(T_xM)\simeq \RP^{n-1}.\]
	Thus it suffices to identify this replacement piece with $\RP^n$ minus an open $n$-ball. But recall that $\RP^n$ is obtained from $\RP^{n-1}$ by attaching an $n$-cell by the standard degree-$2$ map $\S^{n-1}\to \RP^{n-1}$. Equivalently, a tubular neighbourhood of the standard embedded $\RP^{n-1}\subset \RP^n$ has normal bundle equal to its tautological line bundle, finishing the proof of the first statement.
	
	Assume now that $M$ is oriented and that $n$ is odd. Put $V=T_xM$. The orientation of $M$ induces an orientation of $V$. We identify $\RP^n\simeq \P(\R\oplus V)$, where $\R\oplus V$ is oriented by putting the positive basis of $\R$ before an oriented basis of $V$. The unit sphere $\S(\R\oplus V)\simeq \S^n$ is oriented as the boundary of the unit ball. Since $\RP^n=\S^n/\{\pm 1\}$, this orientation descends to $\RP^n$ because $n$ is odd and the antipodal map on $\S^n$ is orientation-preserving. With this convention, the identification of the tautological line bundle over $\P(T_xM)$ with $\RP^n\setminus \mathring{\mathbb{D}}^n$ is compatible with the boundary orientations.
\end{proof}

\subsection{Topological blow-up of some fibrations along fibres.}

\subsubsection{Along a fibre of any locally trivial bundle.}

In a locally trivial bundle $p\colon M\to B$, any fibre $F$ has a trivial normal bundle $\NN_{p(F)/B}\times F\simeq T_{p(F)} B \times F$. This product structure is actually compatible with the bundle given by $p$ away from the centre $F$. Here is a precise statement.

\begin{prop}\label{prop:TopBlowUpAlongFibreInLocTrivialBundle} 
	Let $F$ be a fibre of a locally trivial bundle $p\colon M\to B$ with $\dim M\geq 3$. Let $b=p(F)$ and denote by $\pi_F$ and $\pi_b$ the topological blow-ups along $F\subset M$ and $b\in B$. Then there is a natural diffeomorphism $\Bl_F^{\mathrm{top}}M \approx \pi_b^*M$ over $\Bl_b^{\mathrm{top}}B$. We denote the corresponding bundle by $\widetilde{p}\colon\Bl_F^{\mathrm{top}}M\to \Bl_b^{\mathrm{top}}B$.
\end{prop}
\begin{proof}
	Using the notations as in \Cref{Sec:ConstructionTopBlowUp}, we have a canonical identification
	\[\NN_{F/M}\simeq (p_{|F})^*\NN_{b/B}\simeq F\times \NN_{b/B}.\]
	and, choosing appropriate tubular embedding $\tau_F$ and $\tau_b$ of $F$ and $b$, an isomorphism of bundles 
	\[\xymatrix{
	F\times \NN_{b/B}\ar[r]^{\tau_F} \ar[d]^{\pr_2} & \tau_F(\NN_{F/M})\subset M \ar[d]_{\pi_|}\\
	\NN_{b/B} \ar[r]^{\tau_b} & \tau_b(\NN_{b/B})\subset B} \]
	We then get
	\[\P\NN_{F/M}\simeq F\times \P\NN_{b/B} \qquad \textnormal{and} \qquad \pazocal L_F\simeq F\times \pazocal L_b. \]
	Hence, near the exceptional divisor, $\Bl_F^{\mathrm{top}} M\simeq F\times \Bl_b^{\mathrm{top}} B$ and the map $\widetilde{p}$ is locally just the second projection. Away from the exceptional divisors it is the original map $p$. These two definitions agree on the overlap, so they glue to a smooth bundle. This also makes commutativity with the blow-down maps immediate.
	
	The commutative square induces a canonical map
	\[\Phi\colon \Bl_F^{\mathrm{top}}M \longrightarrow \pi_b^*M,\qquad y\longmapsto (\widetilde{p}(y),\pi_F(y)).\]
	It is clear that $\Phi$ is a diffeomorphism and therefore $\Bl_F^{\mathrm{top}}M$ is naturally identified with the total space of $\pi_b^*\pi$.
\end{proof}

\subsubsection{Along a fibre of a Seifert fibration.}

The following is one of our main result and the best topological tool we have to compute the topological blow-up along knots in $3$-manifolds. 

\begin{thm}\label{Th:TopBlowUpSeifertFibrationAlongFibre}
	Let $p\colon M\to B$ be a Seifert fibration over a closed connected orbifold with only cone points, and let $F$ be a fibre. Then $\Bl^{\mathrm{top}}_F M$ admits a Seifert fibration $\widetilde p$, compatible with $p$ through the blow-down map $\pi_F:\Bl^{\mathrm{top}}_F M\to M$.
	
	Moreover, suppose that the Seifert invariants of $p$ are
	\[\big(B; w_1; e; (\alpha_1,\beta_1), \dots, (\alpha_k,\beta_k)\big).\]
	If $F$ is exceptional, we may assume that its local invariants are $(\alpha_k,\beta_k)$. If $F$ is regular, we formally assign to it the invariants $(\alpha_k,\beta_k)=(1,0)$, adding this pair to the	list above.
	
	Let $\bar e\in\Z_2$ denote the reduction of $e$ modulo $2$ when $e$ is integer-valued, and $e$ itself when it is already $\Z_2$-valued. Then the Seifert invariants of $\widetilde p$ are as follows.
	\begin{enumerate}
		\item If $\alpha_k$ is odd, they are
		\[\big(B^\#;\, (w_1,0);\, \bar e + \beta_k\mod 2;\, (\alpha_1,\beta_1), \dots, (\alpha_{k-1},\beta_{k-1})\big),\]
		where $B^\#$ is obtained from $B$ by replacing a small orbifold disc around $p(F)$ by a M\"obius band. In particular, its underlying surface is $|B^\#|\approx |B|\#\,\RP^2$.		
		\item If $\alpha_k$ is even, they are
		\[\big(B^{\mathrm{mir}};\, w_1 ;\, 0;\, (\alpha_1,\beta_1), \dots, (\alpha_{k-1},\beta_{k-1})\big),\]
		where $B^{\mathrm{mir}}$ is obtained from $B$ by removing a small open orbifold disc around $p(F)$ and declaring the resulting boundary circle to be a reflector boundary.
	\end{enumerate}
	In particular, when $F$ is regular, $\widetilde p$ is just the pullback of $p$ under the topological blow-down
	\[\pi_{p(F)}\colon \Bl^{\mathrm{top}}_{p(F)}|B|\longrightarrow |B|.\]
\end{thm}

\begin{rem}
	Let us make some remarks on the above theorem.
	\begin{enumerate}
		\item In the odd-multiplicity case, let $q\colon |B^\#|\to |B|$ denote the topological blow-down and let $C\subset |B^\#|$ be the image of the exceptional divisor. Then
		\[ H^1(|B^\#|;\Z_2)\simeq q^*H^1(|B|;\Z_2) \oplus \Z_2\langle\eta_C\rangle,\]
		where $\eta_C$ is the Poincar\'e dual of $[C]\in H_1(|B^\#|;\Z_2)$. The notation $(w_1,0)$ refers to this decomposition: the fibre-orientation class of $\widetilde p$ is $q^*w_1$ and vanishes on the new crosscap generator.
		
		In the even-multiplicity case, the underlying surface of $B^{\mathrm{mir}}$ is $|B|$ with an open disc removed. Thus
		\[ H^1(|B^{\mathrm{mir}}|;\Z_2) \simeq H^1(|B|;\Z_2),\]
		and the fibre-orientation class is unchanged. However, because of the reflector circle, the total space is always non-orientable.
		\item Recall that
		\[ H^2(|B|;\Z_{w_1}) \simeq	\begin{cases}
			\Z,&w_1=w_1(T|B|),\\
			\Z_2,&w_1\neq w_1(T|B|).
		\end{cases}\]
		In the odd-multiplicity case, the addition of a crosscap gives $w_1(\widetilde p) \neq w_1(T|B^\#|)$, and hence 
		\[H^2(|B^\#|;\Z_{w_1(\widetilde p)}) \simeq \Z_2.\]
		In the even-multiplicity case, the arising of a reflector circle makes consideration about the Euler class irrelevant.
		\item In the above theorem, we assumed that the starting base space has no reflector circle. It is just for notational simplification. The result also holds in this case, provided the blown-up fibre is a regular or an isolated exceptional one. The only change concerns the Euler number. Indeed, in these settings, it is set initially to $0$ and remains $0$ whatsoever.
		\item We did not seek a description of the blow-up along an exceptional fibre with solid Klein bottle local model, the need for one did not show up.
	\end{enumerate}
\end{rem}

\begin{proof}
	We first show that the Seifert fibration lifts through the topological blow-up. Choose a standard fibred neighbourhood $N_{\alpha,\beta}$ of $F$, written as
	\[ N_{\alpha,\beta} = \frac{\mathbb D^2_z\times \S^1_u}{\langle\gamma\rangle}, \qquad \gamma\cdot (z,u) = \left(e^{2\pi i/\alpha}z, u-\frac{\beta}{\alpha}\right),\]
	where $\S^1_u=\R/\Z$. Since $\alpha$ and $\beta$ are coprime, the action of $\langle\gamma\rangle$ is free, and the quotient map
	\[ q\colon\mathbb D^2_z\times \S^1_u \longrightarrow N_{\alpha,\beta}\]
	is a finite covering of degree $\alpha$. Locally on $N_{\alpha,\beta}$, the Seifert fibration is given by the action of a circle group $\S^1_{\mathrm{fib}}$: the projection of the translation on the second factor of the covering space: $t\cdot(z,u)=(z,u+t)$. This action commutes with $\gamma$, so it descends to $N_{\alpha,\beta}$. In particular, the central fibre is $F=q(\widehat F)$, where $\widehat F=\{0\}\times \S^1_u$.
	
	Now blow up $\widehat F$ in the covering space via \Cref{prop:TopBlowUpAlongFibreInLocTrivialBundle}, yielding the circle bundle isomorphism:
	\[ \Bl^{\mathrm{top}}_{\widehat F} (\mathbb D^2_z\times \S^1_u) \simeq \Bl^{\mathrm{top}}_0\mathbb D^2_z\times\S^1_u.\]
	With the immediate identification
	\[ \Bl^{\mathrm{top}}_0\mathbb D^2 = \left\{(z,\ell)\in\mathbb D^2\times\RP^1 \mid z\in\ell \right\},\]
	the circle bundle on $ \Bl^{\mathrm{top}}_0\mathbb D^2_z\times\S^1_u$ is induced by the action of $\S^1_{\mathrm{fib}}$:
	\[t \cdot (z,\ell,u) = (z, \ell, u+t).\]
	The exceptional divisor in this covering space is $\RP^1\times \S^1_u$. The action of $\gamma$ lifts naturally to a free action on this blow-up:
	\[\gamma \cdot (z,\ell,u) =	\left(e^{2\pi i/\alpha}z, e^{2\pi i/\alpha}\ell,u-\frac{\beta}{\alpha}\right).\]
	Moreover, because $q$ is a finite covering, the normal bundle of $F$ in $N_{\alpha,\beta}$ is the quotient of the normal bundle of $\widehat F$ in $\mathbb D^2_z\times \S^1_u$ by $\langle\gamma\rangle$ yielding the natural identification
	\[\Bl^{\mathrm{top}}_F N_{\alpha,\beta} \simeq\frac{\Bl^{\mathrm{top}}_0\mathbb D^2_z\times \S^1_u}{\langle\gamma\rangle},\]
	compatible with the respective blow-down maps. 	
	This action of $\S^1_{\mathrm{fib}}$ still commutes with the one of $\langle\gamma\rangle$, and hence descends to an action on $\Bl^{\mathrm{top}}_F N_{\alpha,\beta}$ for which the blow-down map is equivariant. Its orbit space is
	\[\widetilde{D}_\alpha = \frac{\Bl^{\mathrm{top}}_0\mathbb D^2_z}{\langle z\mapsto e^{2\pi i/\alpha}z\rangle}.\]	
	Outside the exceptional divisor, the blow-down is a diffeomorphism, and the lifted action agrees with the original Seifert fibration under this identification. On the exceptional divisor the stabiliser of a point represented by $(0,\ell,u)$ consists of those translations corresponding to powers $\gamma^j$ such that $e^{2\pi i j/\alpha}\ell=\ell.$ Thus the stabiliser is trivial if $\alpha$ is odd, and has order $2$ if $\alpha$ is even. In the latter case, the nontrivial element acts by reflection in the direction transverse to the exceptional divisor. Consequently this local action defines a Seifert fibration, allowing reflector fibres in the even case and glueing with $p$ away from $F$, yielding the desired fibration $\widetilde p$.
	
	We now determine the base orbifold and the Seifert invariants of $\widetilde p$. First, suppose that $\alpha$ is even. The half-turn $z\mapsto-z$, corresponding to the element $\gamma^{\alpha/2}$, fixes every projective direction $\ell\in\RP^1$, but reverses the normal direction to the exceptional circle in $\Bl^{\mathrm{top}}_0\mathbb D^2$. No other nontrivial rotation fixes a projective direction. It follows that $\widetilde D_\alpha$ is an annulus whose inner boundary is a reflector circle, while its outer boundary is ordinary and that $\widetilde B$ is obtained from $B$ by removing a small open orbifold disc around $p(F)$ and declaring
	the resulting boundary circle to be a reflector boundary. The Stiefel--Whitney class is just the restriction of the starting one, and the global Euler number is set to $0$ due to the presence of a reflector circle. The other local orbital invariants remain unchanged.
	
	Now, suppose that $\alpha$ is odd. In this case, the induced action of $\Z_\alpha$ on the exceptional circle $\RP^1$ is free. Moreover, the quotient $\widetilde D_\alpha$ is a M\"obius band. Indeed, $\Bl^{\mathrm{top}}_0\mathbb D^2$ is itself a M\"obius band, and its exceptional circle is a deformation retract. The quotient map restricts to a covering of degree $\alpha$ between the exceptional circles. Since $\alpha$ is odd, the normal line bundle of the quotient circle remains nonorientable. Consequently, the base $\widetilde B$ is obtained from $B$ by replacing the orbifold disc around $p(F)$ by a M\"obius band. The lifted circle action preserves the
	orientation of its fibres over the M\"obius band, since it is induced by translations on $\S^1_u$, hence the claimed Stiefel--Whitney class.
	
	It remains to compute the global Euler number, showing it equals $\bar e+ \beta_k\pmod 2$, which admits the following geometric interpretation. The local invariant $\beta_k$ measures the twisting of the Seifert fibration around the exceptional fibre, relative to the chosen boundary section. Blowing up this fibre removes its exceptional nature, but does not eliminate the corresponding twisting, which now contributes to the global Euler invariant. Meanwhile, the replacement of the orbifold disc by a M\"obius band makes this invariant binary and thus independent of the boundary section. To establish this formula, we shall compare the boundary section associated with the standard local invariant $(\alpha,\beta)$ with a section of the circle bundle over the new M\"obius band.
	
	Choose angular coordinates
	\[z = re^{2\pi i\theta}, \qquad \theta\in\R/\Z,\]
	on the punctured transverse disc. Denote by $h$ the homology class of a regular Seifert fibre,
	represented in the covering space by one turn around $\S^1_u$. On the boundary of the original Seifert fibred solid torus, consider the section curve $\sigma$ represented by
	\[t \longmapsto \left(re^{2\pi i t/\alpha}, -\frac{\beta t}{\alpha}\right), \qquad 0\leq t\leq 1.\]	
	On the exceptional divisor of the blown-up covering space, a projective direction is represented by
	\[\ell_\theta=\R e^{2\pi i\theta}, \qquad \ell_{\theta+1/2}=\ell_\theta.\]
	We consider the following section over the core of the M\"obius band
	\[t \longmapsto \left[0, \ell_{t/(2\alpha)}, -\frac{\beta t}{2}-\frac{\beta t}{2\alpha}\right], \qquad 0\leq  t\leq 1,\]
	Its endpoints differ by the action of $\gamma^{\frac{\alpha+1}{2}}$, using the identification of opposite projective directions. Denote its homology class by $c$. As the lifted circle fibration is principal over this band, we may extend the section over the whole band.
	
	We now compare $\sigma$ and $c$ in the homology of the blown-up neighbourhood. Traversing $c$ twice
	corresponds in the covering coordinates to the increments
	\[ \Delta\theta=\frac{1}{\alpha}, \qquad \Delta u	= -\frac{\beta}{\alpha}-\beta\]
	while the curve $\sigma$ has increments
	\[\Delta\theta = \frac{1}{\alpha}, \qquad \Delta u = -\frac{\beta}{\alpha}.\]
	Consequently, we obtain the relation
	\[ 2c = \sigma - \beta h\]
	in the first homology group of the blown-up	neighbourhood.
	
	Since the boundary of a M\"obius band represents twice the class of its core, the chosen section over the new M\"obius band restricts to the boundary as $2c$. Thus the relation above shows that its	boundary section differs from the original standard section $\sigma$ by $\beta$ turns around the	regular fibre.
	
	This correction is independent of the choice of sections, when taken modulo $2$. Indeed, two sections over the M\"obius band differ by a map to $\S^1$, and the restriction homomorphism
	\[ H^1(\mathrm{M\ddot ob};\Z) \longrightarrow H^1(\partial\mathrm{M\ddot ob};\Z)\]
	is multiplication by $2$. Hence changing the section changes its boundary winding by an even integer. This finishes the proof.
\end{proof}

We may use \Cref{Th:TopBlowUpSeifertFibrationAlongFibre} to compute the following basic blow-up, which will arise a lot throughout.

\begin{prop}\label{prop:TopBlowUpAlongUnknot}
	Let $K\subset M$ be a trivial knot in a $3$-manifold. Then,
	\[\Bl_K^{\mathrm{top}} M\approx M\#\,\S^2\widetilde{\times} \S^1.\]
\end{prop}

\begin{proof}
	The trivial knot $K$ is contained in a $3$-ball. Blowing-up being a local operation, we focus on the transformation of the $3$-ball. Gluing another $3$-ball along its boundary yields $\S^3$ on which we may choose a Hopf fibration realising $K$ as a fibre. The Seifert invariants of this bundle are $\big(\S^2;0;-1\big)$ (see for instance \cite{geigesSeifertFibrationsLens2017} or \cite{kineiderClassificationSeifertFibrations2026}). Applying \Cref{Th:TopBlowUpSeifertFibrationAlongFibre}, the blow-up along $K$ produces $M\big(\RP^2;0;\overline{1}\big)\approx \S^2\widetilde{\times} \S^1$. Thus, the topological blow-up of the $3$-ball along $K$ yields $\S^2\widetilde{\times}\S^1$ minus a $3$-ball. Gluing this along the boundary $2$-sphere of the $3$-ball containing $K$ produces the claimed connected sum.
\end{proof}

\begin{rem}
	We may give an alternate proof of \Cref{prop:TopBlowUpAlongUnknot} using birational geometry. We take the concrete model $K=c_{2,r}(\R)$ with $c_{2,r}=\{x_0^2+x_1^2-x_3^2=x_2=x_4=0\}\subset \P^3$. The result is then a combination of the \Cref{prop:BlowUpPointQuadricThreefolds} below and \Cref{prop:TopBlowUpInAPoint}.
\end{rem}
The following description is well-known.

\begin{prop}\label{prop:BlowUpPointQuadricThreefolds}
	Let $\pazocal Q_{r,s}$ be the smooth quadric hypersurface
	\[\pazocal Q_{r,s}:=\{x_0^2+\dots+x_{r-1}^2-x_r^2-\dots-x_{r+s-1}^2=0\}\subset \P^{r+s-1}.\]
	The blow-up of $\P^3_\R$ along the smooth $\R$-rational conic $c_{2,r} = \{x_0^2+x_1^2-x_3^2=x_2=x_4=0\}$ is isomorphic to $\pazocal{Q}_{3,2}$ blown-up in a real point.\\
	The blow-up of $\P^3_\R$ along the smooth conic $c_{2,\emptyset}=\{x_0^2+x_1^2+x_2^2=x_3=x_4=0\}$ with no real point is isomorphic to $\pazocal{Q}_{4,1}$ blown-up in a real point.
\end{prop}

\section{Topological blow-ups along torus and hyperboloidal knots.}

In this section, we introduce several classes of knots and determine the diffeomorphism types of the topological blow-ups along them. Our terminology is adapted to the algebro-geometric viewpoint adopted throughout the paper.

A \emph{knot} in a smooth $3$-manifold $M$ is a smoothly embedded circle. Knots are regarded as unoriented. Two knots are \emph{equivalent} if there is a diffeomorphism of the ambient manifold taking one onto the other. We do not require this map to preserve orientation. In particular, a knot and its mirror image are equivalent. We focus on knots, although links arise naturally when passing to covering spaces.

\subsection{Torus knots in $\S^3$.}\label{Sec:TorusKnotS3}

A \emph{torus knot} is a knot in $\S^3$ lying on the boundary of an unknotted solid torus. Fix an unknotted torus $\mathbb T\subset\S^3$. Its complement has two components, whose closures are solid tori. A \emph{meridian--longitude pair} for $\mathbb T$ is a pair $(\mu,\lambda)$ of oriented simple closed curves on $\mathbb T$, intersecting transversely once, such that $\mu$ bounds a disc in one of these solid tori and $\lambda$ bounds a disc in the other. Their homology classes give a basis
\[ H_1(\mathbb T;\Z)=\Z\mu\oplus\Z\lambda.\]

\begin{dfn}
	For coprime integers $r,s$, an \emph{$(r,s)$-torus knot} is a knot on $\mathbb T$ which can be oriented so as to represent the class
	\[r\mu+s\lambda\in H_1(\mathbb T;\Z).\]
\end{dfn}

Every essential simple closed curve on $\mathbb T$ represents a primitive homology class and hence is an $(r,s)$-torus knot for some coprime integers $r,s$. An inessential simple closed curve bounds a disc on $\mathbb T$ and is therefore trivial.

An $(r,s)$-torus knot is trivial if and only if $|r|=1$ or $|s|=1$. Since our equivalence relation allows orientation-reversing maps, changing the signs of the parameters or interchanging them does not
change the equivalence class. Thus every nontrivial torus knot is equivalent to one with coprime parameters $2\leq s\leq r$. The following classification is classical.

\begin{lem}\label{lem:ClassesTorusKnots}
	For coprime integers $2\leq s\leq r$, the equivalence classes of $(r,s)$-torus knots are pairwise distinct.
\end{lem}

Every trivial torus knot is equivalent to the $(1,1)$-torus knot. Consequently, positive parameters suffice for the following statement.

\begin{prop}\label{prop:TopBlowUpTorusKnotS3}
	Let $K$ be an $(r,s)$-torus knot in $\S^3$, where $r,s$ are coprime positive integers. If $r=1$ or $s=1$, then
	\[\Bl^{\mathrm{top}}_K\S^3 \approx \S^2\widetilde{\times}\S^1.\]
	Otherwise,
	\[\Bl^{\mathrm{top}}_K\S^3 \approx M\bigl(\RP^2;0;0;(r,\beta_r),(s,\beta_s)\bigr),\]
	where $\beta_r$ and $\beta_s$ are uniquely determined by
	\[\begin{aligned}
		&0<\beta_r\leq r/2,
		&&s\beta_r\equiv \pm1\pmod r,\\
		&0<\beta_s\leq s/2,
		&&r\beta_s\equiv \pm1\pmod s.
	\end{aligned}\]
\end{prop}

\begin{proof}
	If $r=1$ or $s=1$, the knot is trivial and the assertion follows from \Cref{prop:TopBlowUpAlongUnknot}.
	
	Assume henceforth that $r,s\geq2$. The sphere $\S^3$ admits a Seifert fibration for which $K$ is a regular fibre and whose two exceptional fibres have multiplicities $r$ and $s$; see \cite{kineiderClassificationSeifertFibrations2026} or \cite{geigesSeifertFibrationsLens2017}. Choose integers $u,v$ satisfying
	\[0<u<r,\qquad 0<v<s,\qquad su\equiv1\pmod r,\qquad rv\equiv1\pmod s. \]
	Since $r$ and $s$ are coprime, these congruences give $su+rv\equiv1\pmod{rs}$. The inequalities $1<su+rv<2rs$ therefore imply
	\[ su+rv=rs+1.\]
	With our conventions, the Seifert invariants of this fibration are
	\[ \left(\S^2;0;\frac{1-su-rv}{rs};(r,u),(s,v)\right) = \bigl(\S^2;0;-1;(r,u),(s,v)\bigr).\]
	Applying \Cref{Th:TopBlowUpSeifertFibrationAlongFibre} to the regular fibre $K$ yields
	\[\Bl^{\mathrm{top}}_K\S^3 \approx M\bigl(\RP^2;0;1;(r,u),(s,v)\bigr).\]
	It remains to normalise the local invariants. For this nonorientable Seifert manifold, replacing a pair $(\alpha,\beta)$ by $(\alpha,\alpha-\beta)$ changes the obstruction by $1$ modulo $2$. Exactly one of the inequalities $u>r/2$ and $v>s/2$ holds. Indeed, if neither held, we would have $su+rv\leq rs$.	If both held, then
	\[2=s(2u-r)+r(2v-s)\geq r+s>2,\]
	which is impossible. Thus, setting
	\[\beta_r=\min(u,r-u),\qquad \beta_s=\min(v,s-v),\]
	normalises exactly one local pair by replacing its second entry by its complement. The obstruction consequently becomes $0$. The resulting coefficients satisfy the stated bounds and congruences, which determine them uniquely.
\end{proof}

\begin{rem}[Orbifold Euler classes of the blow-ups, torus case]\label{rem:OrbifoldEulerClassesTorus}
	The above Seifert invariants may appear quite \say{messy}. As the computation in \cite{kineiderClassificationSeifertFibrations2026} shows, the above generalised Hopf fibration on $\S^3$ admits a much cleaner description in terms of orbifold Euler class. The Seifert invariants formalism is responsible for the artificial technicity. Here we reinterpret the result and its proof using the orbifold Euler class formalism.	All the fibrations considered below have $w_1=0$, so their coefficient systems are constant. 
	
	For a trivial knot in $\S^3$, one may use the genuine Hopf fibration: its Euler class $\mathcal E\in H^2(\S^2;\Z)\simeq \Z$ is a generator and pulls back to the nonzero class in $H^2(\RP^2;\Z)\simeq\Z_2$, giving $\S^2\widetilde{\times}\S^1$.
	
	Write $\RP^2(m,n)$ for the projective plane with cone points of orders $m,n$. Let
	\[ f_{m,n}\colon\RP^2(m,n)\longrightarrow\S^2(m,n) \]
	be the topological blow-down at a regular point. If $\varpi$ is a Seifert fibration over $\S^2(m,n)$ and $F$ is a regular fibre, then \Cref{Th:TopBlowUpSeifertFibrationAlongFibre} identifies the fibration on $\Bl^{\mathrm{top}}_F M$ with $f_{m,n}^*\varpi$. Denoting $\mathcal{E}(\varpi)\in H^2_{\mathrm{orb}}(\S^2(m,n);\Z)$ the Euler class of the initial fibration, the naturality under pullback yields:
	\[\mathcal E(f_{m,n}^*\varpi)=f_{m,n}^*\mathcal E(\varpi).\]
	
	In the presentation
	\[H^2_{\mathrm{orb}}(\S^2(m,n);\Z) =\left\langle c,q_m,q_n \;\middle|\; mq_m=c,\; nq_n=c\right\rangle, \]
	the class $c$ is the surface class. The pullback $f_{m,n}^*$ is the quotient obtained by adding the relation $2c=0$:
	\[ H^2_{\mathrm{orb}}(\RP^2(m,n);\Z) \simeq \frac{H^2_{\mathrm{orb}}(\S^2(m,n);\Z)}{\langle2c\rangle}.\]
	Indeed, the surface class pulls back to the nonzero element of $H^2(\RP^2;\Z)\simeq\Z_2$, while the local generators and their relations are preserved. We denote pullback classes by a hat.
	
	Now, consider that $K$ is an $(r,s)$-torus knot in $\S^3$, with coprime integers $r,s\geq2$. By the calculation in \cite{kineiderClassificationSeifertFibrations2026},
	\[ H^2_{\mathrm{orb}}(\S^2(r,s);\Z)=\Z u, \qquad c=rsu,\quad q_r=su,\quad q_s=ru,\]
	and $\mathcal E(\pi_{r,s})=u$. Thus the induced Seifert fibration $\widetilde\pi_{r,s}\colon \Bl^{\mathrm{top}}_K\S^3 \longrightarrow \RP^2(r,s)$ has Euler class
	\[\mathcal E(\widetilde\pi_{r,s})=\widehat u,\]
	where $\widehat u$ is just a generator of $H^2_{\mathrm{orb}}(\RP^2(r,s);\Z)\simeq \Z_{2rs}$; in particular it has order $2rs$.
\end{rem}

\subsection{Hyperboloidal knots in $\RP^3$.}\label{Sec:HyperboloidalKnots}

We follow the description of hyperboloidal knots and links in \cite{mikhalkinRealAlgebraicKnots2016}, with a slight change of terminology: we distinguish knots contained in a smooth real quadric from knots that are merely equivalent to such knots. Throughout this subsection, we identify $\RP^3=\P^3(\R)$.

\begin{dfn}
	A \emph{hyperboloid} in $\RP^3$ is an embedded torus realising a genus-$1$ Heegaard splitting of $\RP^3$.
\end{dfn}

It is classical that all such tori are smoothly ambient isotopic. The real locus of any smooth real quadric surface in $\P^3$ of signature $(2,2)$ is a hyperboloid.

\begin{dfn}
	A knot in $\RP^3$ is \emph{hyperboloidal} if it is contained in the real locus of a smooth real quadric surface of signature $(2,2)$. It is \emph{topologically hyperboloidal} if it is equivalent to a hyperboloidal knot.
\end{dfn}

A knot is called \emph{topologically planar} if it is equivalent to a knot contained in a real projective plane $\RP^2\subset\RP^3$. There are exactly two classes of such knots:
\begin{itemize}
	\item the class $\langle\mu_0\rangle$ of the trivial knot;
	\item the class $\langle\mu_1\rangle$ of the real locus of a real line.
\end{itemize}
They are distinguished by their classes in $\pi_1(\RP^3)\simeq\Z_2$, and both are topologically hyperboloidal.

Fix a smooth real quadric $Q\simeq\P^1_\R\times\P^1_\R\subset\P^3$ of signature $(2,2)$, and set $\mathcal H=Q(\R)$. Choose real lines $\ell_1,\ell_2$ in the two rulings of $Q$. After orienting their real loci $s_i=\ell_i(\R)$, we obtain a basis
\[ H_1(\mathcal H;\Z)=\Z s_1\oplus\Z s_2. \]
Unlike the meridian--longitude basis used above, this basis comes from the algebraic structure of the quadric. It is particularly convenient for relating knot classes to the bidegrees of algebraic curves on $Q$. The choice of generators is well defined up to signs and permutation. These changes are induced by ambient diffeomorphisms of $\RP^3$, possibly reversing orientation, and therefore do not affect the equivalence classes considered here.

Let $f\colon\S^3\to\RP^3$ be the universal covering and put $\mathbb T=f^{-1}(\mathcal H)$. Each ruling line represents the nonzero element of $\pi_1(\RP^3)$, so the restricted covering $f_{\mathcal H}\colon\mathbb T\to\mathcal H$ corresponds to the subgroup
\[ \left\{(a,b)\in\Z^2\mid a+b\equiv 0\pmod 2\right\} \subset \pi_1(\mathcal H)\simeq H_1(\mathcal H;\Z).\]
In the standard model of a quadric of signature $(2,2)$, $\mathbb T$ is an unknotted torus in $\S^3$.
Moreover, one can choose a meridian--longitude pair $(\mu,\lambda)$ such that
\[(f_{\mathcal H})_*\mu=s_1+s_2,\qquad (f_{\mathcal H})_*\lambda=s_1-s_2.\]
The corresponding meridian discs project to meridian discs of the two solid tori bounded by $\mathcal H$. Consequently, simple closed curves of classes $(1,1)$ and $(1,-1)$, as well as their negatives, belong to $\langle\mu_0\rangle$.

Let $K\subset\mathcal H$ be a knot. After choosing an orientation of $K$, write
\[ [K]=p s_1+q s_2. \]
If this class is nonzero, it is primitive. By changing the signs and order of the ruling generators, we may assume that $0\leq q\leq p$.

Suppose first that $[K]\neq 0$. The inverse image $f^{-1}(K)$ is a $(p+q,p-q)$-torus link. If $p+q$ is odd, it is connected and represents
\[(p+q)\mu+(p-q)\lambda.\]
If $p+q$ is even, it consists of two components, each representing
\[ \frac{p+q}{2}\mu+\frac{p-q}{2}\lambda.\]
Indeed, these assertions follow from the displayed formula for $(f_{\mathcal H})_*$ and from the fact that the class of $K$ in $\pi_1(\RP^3)\simeq\Z_2$ is $p+q\bmod2$. If $[K]=0$, then $K$ bounds a disc on $\mathcal H$, and its inverse image is a two-component unlink.

\begin{prop}[\cite{mikhalkinRealAlgebraicKnots2016}]\label{prop:PlanarHyperboloidalKnots}
	Let $K\subset\mathcal H$ be a knot, and choose the ruling generators so that $[K]=p s_1+q s_2$ with $0\leq q\leq p$. Then:
	\begin{enumerate}
		\item If $(p,q)=(0,0)$ or $(p,q)=(1,1)$, then $K\in\langle\mu_0\rangle$.
		\item If $p=q+1$, then $K\in\langle\mu_1\rangle$.
		\item Otherwise, $p,q$ are coprime, $p>q+1\geq2$, and $K$ is not topologically planar.
	\end{enumerate}
\end{prop}

\begin{proof}
	If $[K]=0$, the knot bounds a disc on $\mathcal H$.	If $(p,q)=(1,1)$, it bounds a meridian disc of one of the solid tori bounded by $\mathcal H$. Both cases therefore give the trivial knot.
	
	If $p=q+1$, the geometric intersection number of $K$ with a	meridian of class $(1,1)$ is
	\[\left|\det \begin{pmatrix}
		p&q\\
		1&1
	\end{pmatrix}\right| = |p-q| = 1.\]
	Pushing $K$ into the corresponding solid torus shows that it is isotopic to its core, hence to a real projective line.
	
	In all remaining cases, primitivity gives $p>q+1\geq2$. The inverse image of $K$ is then a torus link whose two parameters satisfy $p+q>p-q\geq2$. Such a link is neither an unknot nor a two-component unlink. Since these are the inverse images of the two topologically planar knot classes, $K$ cannot be topologically planar.
\end{proof}

The classical classification of torus links extends \Cref{lem:ClassesTorusKnots}: for integers $a\geq b\geq2$, the $(a,b)$-torus links are pairwise inequivalent for distinct pairs $(a,b)$. Applying this to the double covers gives the following result.

\begin{prop}[{\cite[Proposition~8]{mikhalkinRealAlgebraicKnots2016}}]\label{prop:ClassesHyperboloidalKnots}
	For distinct coprime pairs $(p,q)$ satisfying $p>q+1\geq2$, knots on $\mathcal H$ representing the classes $p s_1+q s_2$ are pairwise inequivalent in $\RP^3$.
\end{prop}

\begin{proof}
	Any homeomorphism of $\RP^3$ taking one knot onto another lifts to a homeomorphism of their inverse images in $\S^3$. The torus links $T(p+q,p-q)$ therefore determine distinct equivalence
	classes, since their parameters determine
	\[p=\frac{(p+q)+(p-q)}{2}, \qquad q=\frac{(p+q)-(p-q)}{2}.\]
\end{proof}

\begin{dfn}[Definition~7 in \cite{mikhalkinRealAlgebraicKnots2016}]
	For coprime integers $p,q$ satisfying $p>q+1\geq2$, denote by $h_{p,q}$ the equivalence class of a knot on $\mathcal H$ representing the class $p s_1+q s_2$, for a suitable choice of signs and ordering of the ruling generators. In accordance with our convention, this class includes the mirror image of any of its representatives.
\end{dfn}

\begin{rem}
	In the terminology of low-dimensional topology, topologically hyperboloidal knots are \emph{torus knots in the lens space} $\RP^3$: they can be isotoped onto the splitting torus of a genus-$1$ Heegaard decomposition; see \cite{yeConstrainedKnotsLens2023}. The parameters $(p,q)$ used here describe the slope in the basis of ruling lines, rather than in the usual topological basis. Knots in the class $\langle\mu_1\rangle$ are \emph{core knots}: their exteriors are homeomorphic to a solid torus.
\end{rem}

\begin{thm}\label{Th:TopBlowUpHyperboloidalLink}
	Let $K\subset\RP^3$ be a topologically hyperboloidal knot.
	\begin{enumerate}
		\item If $K\in\langle\mu_0\rangle$, then
		\[\Bl^{\mathrm{top}}_K\RP^3 \approx \RP^3\#\,\S^2\widetilde{\times}\S^1.\]
		\item If $K\in\langle\mu_1\rangle$, then
		\[\Bl^{\mathrm{top}}_K\RP^3 \approx \RP^2\times\S^1.\]		
		\item Otherwise, $K\in h_{p,q}$ for coprime integers $p,q$ satisfying $p>q+1\geq2$, and
		\[ \Bl^{\mathrm{top}}_K\RP^3 \approx M\bigl(\RP^2; 0; \overline{0}; (p-q,\beta_-), (p+q,\beta_+)\bigr),\]
		where $\beta_-$ and $\beta_+$ are uniquely determined by
		\[\begin{aligned}
			&0<\beta_-\leq\frac{p-q}{2},
			&&p\beta_-\equiv\pm1\pmod{p-q},\\
			&0<\beta_+\leq\frac{p+q}{2},
			&&p\beta_+\equiv\pm1\pmod{p+q}.
		\end{aligned}\]
		The signs in these two congruences are chosen independently.
	\end{enumerate}
\end{thm}

\begin{proof}
	The diffeomorphism type of the blow-up depends only on the equivalence class of the knot, so we may choose convenient representatives.
		
	If $K\in\langle\mu_0\rangle$, this is just an application of \Cref{prop:TopBlowUpAlongUnknot}.
		
	If $K\in\langle\mu_1\rangle$, take $K$ to be a real projective line. The Hopf fibration on $\S^3$ descends through the antipodal quotient to a circle bundle $\RP^3\to\S^2$ for which $K$ is a fibre.
	Its Euler number is $\pm 2$, depending on the orientation convention. By \Cref{Th:TopBlowUpSeifertFibrationAlongFibre}, blowing up this fibre gives a circle bundle over $\RP^2$ with $w_1=0$ and obstruction $\overline{0}\in\Z_2$. This is the trivial bundle $\RP^2\times\S^1$.
	
	Now suppose that $K\in h_{p,q}$, and put
	\[ A=p+q,\qquad B=p-q,\qquad d=\gcd(A,B)\in\{1,2\}.\]
	The inverse image of a standard representative of $K$ in $\S^3$	is the torus link $T(A,B)$.	Its components are regular fibres of the Seifert fibration defined by the circle action
	\[(z_1,z_2)\longmapsto \bigl(e^{i(A/d)\theta}z_1,e^{i(B/d)\theta}z_2\bigr)\qquad \text{on }\S^3\subset\C^2.\]
	This action commutes with the antipodal involution, and its orbit decomposition therefore descends to a Seifert fibration of $\RP^3$ having $K$ as a regular fibre. The antipodal involution induces either the identity or an orientation-preserving involution on the upstairs base $\S^2$; the quotient base is therefore again $\S^2$.
		
	To compute the invariants downstairs, it is convenient to write the effective descended action as
	\[t\cdot[z_1,z_2] = [e^{iAt}z_1,e^{iBt}z_2], \qquad t\in\R/\pi\Z,\]
	where $[z_1,z_2]$ denotes the class of $(z_1,z_2)\in\S^3$ modulo the antipodal involution. The action is well defined because $A$ and $B$ have the same parity, and effective because $\gcd(p,q)=1$. Its base is $\S^2$, and its two exceptional fibres are the images of the coordinate circles. Their multiplicities are $A$ and $B$. The fibres are consistently oriented, so $w_1=0$.
	
	At the image of $(1,0)$, the stabiliser is generated by the time $t=\pi/A$. On a transverse disc, this element acts by multiplication by
	\[-e^{i\pi B/A} = e^{i\pi(A+B)/A} =	e^{2\pi ip/A}. \]
	Here the minus sign accounts for the antipodal identification of $(-1,0)$ with $(1,0)$.
	Likewise, at the image of $(0,1)$, the stabiliser generator $t=\pi/B$ acts on a transverse disc by multiplication by $e^{2\pi ip/B}$. Thus the local Seifert invariants may be written as $(B,u)$ and
	$(A,v)$, where
	\[0<u<B,\qquad 0<v<A,\qquad pu\equiv1\pmod B,\qquad pv\equiv1\pmod A.\]
	These coefficients exist and are unique since $\gcd(p,A)=\gcd(p,B)=1$. Write the full Seifert invariants as
	\[\bigl(\S^2;0;e;(B,u),(A,v)\bigr).\]
	We first observe that
	\[p(Au+Bv)\equiv A+B=2p\pmod{AB}.\]
	Since $\gcd(p,AB)=1$, it follows that $Au+Bv\equiv2\pmod{AB}$. Together with $2<Au+Bv<2AB$, this gives
	\[Au+Bv=AB+2.\]
	The standard homology presentation for these Seifert invariants gives
	\[ \bigl|H_1(\RP^3;\Z)\bigr| = \bigl|ABe+Au+Bv\bigr| = \bigl|AB(e+1)+2\bigr|.\]
	The left-hand side equals $2$. Since $AB\geq8$ and $e\in\Z$, we obtain $e=-1$. Consequently, the descended fibration has invariants
	\[\bigl(\S^2;0;-1;(B,u),(A,v)\bigr).\]
	
	Applying \Cref{Th:TopBlowUpSeifertFibrationAlongFibre} to the regular fibre $K$ now yields
	\[\Bl^{\mathrm{top}}_K\RP^3 \approx M\bigl(\RP^2;0;\overline{1};(B,u),(A,v)\bigr).	\]
	It remains to normalise the local invariants. Exactly one of the inequalities $u>B/2$ and $v>A/2$ holds. Indeed, if neither held, then $Au+Bv\leq AB$. If both held, the identity $Au+Bv=AB+2$ would imply
	\[4=A(2u-B)+B(2v-A)\geq A+B=2p\geq6,\]
	a contradiction. Set
	\[\beta_-=\min(u,B-u),\qquad \beta_+=\min(v,A-v).\]
	Exactly one local coefficient is therefore replaced by its complement. For this nonorientable Seifert fibration, this replacement changes the obstruction by $1$ modulo $2$, so its normalised value is $\overline{0}$. The resulting local coefficients satisfy precisely the bounds and congruences in the statement, proving the last assertion.
\end{proof}

\begin{rem}[Orbifold Euler classes of the blow-ups, hyperboloidal case]\label{rem:OrbifoldEulerClassHyperboloidal}
	As in \Cref{rem:OrbifoldEulerClassesTorus}, we may give a reinterpretation of this result in terms of orbifold Euler class. We use the notations introduced there.
	
	For the projective-line class in $\RP^3$, the descended Hopf bundle has	Euler class $2c$ over $\S^2$. Its pullback to $\RP^2$ has Euler class zero, recovering the product $\RP^2\times\S^1$.
	
	Suppose now that $K\in h_{p,q}$, where $\gcd(p,q)=1$ and $p>q+1\geq2$. Put
	\[A = p+q, \qquad B=p-q.\]
	The descended Seifert fibration $\varpi_{p,q} \colon \RP^3 \to \S^2(A,B)$ has Euler class
	$\mathcal E(\varpi_{p,q})$ characterised by	
	\[ p\mathcal E(\varpi_{p,q}) = q_A+q_B \quad \text{in} \quad H^2_{\mathrm{orb}}(\S^2(A,B);\Z) = \left\langle c,q_A,q_B \;\middle|\; Aq_A=c,\; Bq_B=c\right\rangle.\]
	This equation has a unique solution. Indeed, with the notation of the proof, the orbifold Euler class is
	\[\mathcal E(\varpi_{p,q})=-c+u q_B+v q_A.\]
	Using $Bq_B = Aq_A = c$, the local congruences, and the identity $Au + Bv = AB+2$ established there, we obtain
	\[p\mathcal E(\varpi_{p,q}) = q_A+q_B + \left(\frac{pu-1}{B}+\frac{pv-1}{A}-p\right)c = q_A + q_B,\]
	where the last equality uses $A+B=2p$.Moreover,
	\[H^2_{\mathrm{orb}}(\S^2(A,B);\Z) \simeq \Z\oplus\Z_d \qquad d = \gcd(A,B)\in\{1,2\}.\]
	Since $\gcd(p,d)=1$, multiplication by $p$ is injective on this	group. Thus the displayed equation uniquely determines $\mathcal E(\varpi_{p,q})$. 
	
	The induced fibration on the blow-up has base $\RP^2(A,B)$ and Euler class
	\[ \mathcal E(\widetilde\varpi_{p,q}) = f_{A,B}^*\mathcal E(\varpi_{p,q}).	\]		
	More explicitly, if $p+q$ is odd, then $A,B$ are coprime and
	\[ H^2_{\mathrm{orb}}(\S^2(A,B);\Z)=\Z u, \qquad c=ABu, \quad q_A=Bu, \quad q_B=Au.\]
	The antipodal involution acts along the fibres, so the quotient doubles the Euler class:
	\[ \mathcal E(\varpi_{p,q})=2u.\]
	Hence
	\[ H^2_{\mathrm{orb}}(\RP^2(A,B);\Z) = \Z_{2AB}\widehat u, \qquad \mathcal E (\widetilde\varpi_{p,q}) = 2\widehat u.\]
	
	If $p,q$ are both odd, write $A=2a$ and $B=2b$. Then $a,b$ are coprime and have opposite parity. There are classes $v,\tau$ such that
	\[H^2_{\mathrm{orb}}(\S^2(A,B);\Z) = \Z v \oplus \Z_2\tau,\]
	with
	\[q_A = bv + \tau, \qquad q_B = av + \tau, \qquad c = 2abv = \frac{AB}{2}v.\]
	The descended Euler class is $\mathcal E(\varpi_{p,q}) = v$. Consequently,
	\[H^2_{\mathrm{orb}}(\RP^2(A,B);\Z) = \Z_{AB}\, \widehat v \oplus \Z_2\, \widehat\tau, \qquad \mathcal E (\widetilde\varpi_{p,q}) = \widehat v.\]
	In particular $\mathcal E(\widetilde\varpi_{p,q})$ has order $AB=p^2-q^2$ in any case.
\end{rem}

\begin{rem}
	The Seifert fibrations of $\RP^3$ considered above are compatible with its genus-$1$ Heegaard splitting: both solid tori are unions of fibres, and each contains at most one exceptional fibre,
	located at its core. In the standard models, these cores are the images of the coordinate circles
	$\{z_2=0\}$ and $\{z_1=0\}$ under $\S^3\to\RP^3$, and hence are real projective lines.
	
	After a projective change of coordinates, a prescribed real line $c\subset\RP^3$ may therefore be taken as one of these fibres. Applying \Cref{Th:TopBlowUpSeifertFibrationAlongFibre}, including its exceptional-fibre cases, gives Seifert structures on $\Bl^{\mathrm{top}}_c\RP^3\approx\RP^2\times\S^1$ compatible with the chosen fibrations of $\RP^3$ through the blow-down map. If the multiplicity of $c$ is odd, the new base has
	underlying surface $\RP^2$; if it is even, the new base is a disc with reflector boundary. We obtain this way all but one Seifert structures on $\RP^2\times \S^1$ (see \cite[p.~459]{scottGeometries3Manifolds1983}). The only missing one is above a disc with reflector boundary but with no isolated exceptional fibre.
\end{rem}

\section{Some Fano blow-ups of $\P^3$.}

\subsection{About smooth rational curves $c$ lying on quadric surfaces.}

We first recall some general facts about rational curves in $\P^3$ and their relation with quadric surfaces. Assume first that we work over $\C$. Let $c\subset \P^3$ be a smooth rational curve of degree $d$. The quadrics containing $c$ form the linear system $|\pazocal I_c(2)|$. The short exact sequence
\[0\to \pazocal I_c(2)\to \O_{\P^3}(2)\to \O_c(2)\to 0\]
together with $h^0(c,\O_c(2))=h^0(\P^1,\O_{\P^1}(2d))=2d+1$ and $h^0(\P^3,\O_{\P^3}(2))=10$ shows that, if $d\leq 4$, then $c$ is contained in at least one quadric surface. The following fact is well-known.
\begin{prop}\label{prop:RatCurvesOnQuadricCones}
	A quadric cone in $\P^3$ contains no smooth rational curve of degree $d\geq 4$. In addition, any smooth rational curve of degree $3$ in it passes through the vertex and defines a section of the minimal resolution $\mathbb{F}_2\to \P^1$.
\end{prop}

We postpone the proof of the proposition and continue the discussion.

\begin{cor}\label{cor:QuadricSurfaceContainingQuartic}
	Any smooth rational curve $c$ of degree $d\leq 4$ is contained in at least one smooth quadric surface. If $d=4$, this quadric is unique.
\end{cor}

\begin{proof}[Proof of \Cref{cor:QuadricSurfaceContainingQuartic}]
	The existence of a quadric containing $c$ follows from the dimension count above. If $d\leq 3$, that one may choose the quadric to be smooth. For $d=4$, if two distinct quadrics contained $c$, their intersection would have degree $4$. Since $c$ is irreducible of degree $4$, it would be their complete intersection, hence a curve of genus $1$, a contradiction. Thus the quadric containing $c$ is unique. By the proposition, this quadric cannot be a cone; hence it is smooth.
\end{proof}

For $d\geq 5$, lying on a quadric surface is a special property; see \Cref{rem:QuadricSurfacesAndGeneralCurves}. Nevertheless, there are such curves in every degree. For instance, the curve parametrised by
\[ [s:t]\mapsto [s^d:s^{d-1}t:st^{d-1}:t^d] \]
is contained in the smooth quadric $\mathcal Q=\{x_0x_3-x_1x_2=0\}$. Assume now that $c$ lies on a smooth quadric surface $Q\simeq \P^1\times\P^1$. On $Q$, the curve $c$ has bidegree $(a,b)$, with $a+b=d$. Since $c$ is smooth and rational, the genus formula gives
\[ 0 = g(c) = 1+\frac{1}{2}(c^2+c\cdot K_Q) = (a-1)(b-1). \]
Thus, up to exchanging the two rulings, $(a,b)=(d-1,1)$. Equivalently, $c$ is the graph of a morphism $\P^1\to \P^1$ of degree $d-1$.

\begin{rem}\label{rem:QuadricSurfacesAndGeneralCurves}
	If we further assume that $c$ is not a conic, then its normal bundle is isomorphic to $\O_{\P^1}(2d-1)\oplus \O_{\P^1}(2d-1)$. This is the same normal bundle as that of a general rational curve of degree $d$ in $\P^3$. However, a general smooth rational curve of degree $d$ does not lie on any quadric surface for $d\geq 5$. It is mentioned in \cite{eisenbudNormalBundlesSmooth1981} and proved using the number of $d-1$ secants of $c$. We may also argue that for a general rational curve $c$ of degree $d$, the restriction map 
	\[H^0\big(\P^3,\O_{\P^3}(2)\big)\to H^0\big(c,\O_{c}(2)\big)=H^0\big(\P^1,\O_{\P^1}(2d)\big)\]
	has maximum rank. Thus, the dimensions of these spaces imply that, for $d\geq 5$, the former restriction map is injective. But then $H^0\big(\P^3,\pazocal I_c(2)\big)=0$.
\end{rem}

\begin{proof}[Proof of \Cref{prop:RatCurvesOnQuadricCones}]
	Let $c$ be a smooth rational curve of degree $d$. Let $Q\subset \P^3$ be a quadric cone containing $c$, and let $\pi:\mathbb F_2\to Q$ be its minimal resolution. Denote by $E\subset \mathbb F_2$ the exceptional curve, so that $E^2=-2$, and by $f$ a fibre of the ruling. The pullback of a hyperplane section of $Q$ is $H=E+2f$.	
	Let $\widetilde c\subset \mathbb F_2$ be the strict transform of $C$. Write
	\[	\widetilde c\sim aE+bf.	\]
	Since $E^2=-2$, $E\cdot f=1$, and $f^2=0$, we have
	\[d=\deg c=\widetilde c\cdot H =(aE+bf)\cdot(E+2f) = b.\]
	The canonical class of $\mathbb F_2$ is
	\[K_{\mathbb F_2}=-2E-4f.\]
	Therefore, since $\widetilde c\simeq c$ is smooth and rational, adjunction yields  
	\[-2= 2g(\widetilde c)-2=\widetilde c\cdot(\widetilde c+K_{\mathbb F_2})= -2a^2+2ad-2d.\]
	Equivalently,
	\[(a-1)(d-a-1)=0.\]
	Hence either $a=1$ or $a=d-1$.
	The vertex of $Q$ corresponds to the contracted curve $E$. Since $c$ is smooth, its strict transform $\widetilde c$ can meet $E$ in at most one point and transversely. Therefore $0\leq\widetilde c\cdot E=d-2a\leq 1$. In any case, we get $d\leq 3$ and whenever $d=3$ one must have $\widetilde c\cdot E = 1$, meaning $c$ passes through the vertex of the cone and $\widetilde{c}\cdot f = 1$, in other words, $\widetilde{c}$ is a section of $\mathbb{F}_2\to \P^1$.
\end{proof}

Now suppose that $c$ is rational over $\R$ and that it is lying on a smooth \emph{real} quadric surface in $\P^3$. Then we may actually assume 
\[Q=\mathcal{Q}=\{x_0x_3-x_1x_2=0\}\subset \P^3.\]
Indeed, note first that a real smooth non-empty conic is contained in both real quadric surfaces of signature $(2,2)$ and $(3,1)$. And if $d\neq 2$, the class of $c$ in $\Pic Q$, namely $\O_{\P^1\times \P^1}(d-1,1)$, must be fixed by the action of the Galois group $Gal(\C/\R)\simeq \Z_2$, therefore we must have $Q\simeq \P^1_\R\times\P^1_\R$ and the claim follows. 

Choose a line in each ruling of $\mathcal Q$ :
\[\ell_1\in |\O_{\P^1\times \P^1}(1,0)|,\quad \ell_2\in |\O_{\P^1\times \P^1}(0,1)|.\]
In the case where $c$ is a line, choose $\ell_2\neq c$. Denote by $\mathcal H:= \mathcal{Q}(\R)$ the corresponding hyperboloid and by $s_i=\ell_i(\R)$, $i=1,2$ as in \Cref{Sec:HyperboloidalKnots}. On $\mathcal Q$, we have $c\cdot \ell_1=1$ and since both curves are real, they meet in a unique simple real intersection point. Moreover, non-real intersection points come in conjugate pairs, and intersections of even multiplicity may be cancelled by a small smooth perturbation. Consequently, the integer $\left(c(\R)\cdot s_2\right)_{H_1(\pazocal{\mathcal H};\Z)}$ must lie in $\{-d+1,\dots,d-1\}$ and must have the same parity as $d-1$. This parity condition follows also from the characterisation of the class $[c(\R)]\in \pi_1(\RP^3)$ by the degree. We thus obtain the following result. 

\begin{lem}
	Let $c\subset \P^3$ be a smooth real $\R$-rational curve of degree $d$ contained in a real quadric. Then $c(\R)$ is a hyperboloidal knot in $\RP^3$.\\
	More precisely, $c(\R)$ is a knot of type $\langle \mu_0\rangle$, $\langle \mu_1\rangle$ or
	\[h_{a,1} \quad\textnormal{with} \quad a\in \{3,\dots,d-1\} \enspace\textnormal{and} \enspace a\equiv d \pmod 2.\]
\end{lem}
Unfortunately, we have been unable to give a relevant criterion to decide precisely, for a given curve, which class of knots it represents. However, see the discussion below \Cref{Th:BlowUpProjSpaceAlongRatCurves} for the case $d=4$. Also it seems non-trivial to determine all the hyperboloidal knots obtained as real loci of rational curves of a fixed degree.

\begin{rem}
	\begin{enumerate}
		\item The classification of algebraic links in the projective space carried out in \cite{bjorklundRealAlgebraicKnots2011} and \cite{mikhalkinRealAlgebraicKnots2016} shows that the following hyperboloidal knots are realised by rational curves for degree up to $6$ :
		\begin{equation}\label[table]{table:KnotsFromRationalCurves}
			\begin{array}{|c |c |c |c |c |c |c |}
				\hline
				d  & 1 & 2 & 3 & 4 & 5 & 6 \\
				\hline
				\textnormal{Type of } c(\R) & \langle \mu_1\rangle& \langle \mu_0\rangle & \langle \mu_1\rangle & h_{3,1},\; \langle \mu_0\rangle & h_{4,1}, \;\langle \mu_1\rangle & h_{7,5},\; h_{5,3},\; h_{5,1},\;  h_{3,1} \\
				\hline
			\end{array}
		\end{equation}
		The case $\langle \mu_0\rangle$ for degree $4$ is missing in \cite{mikhalkinRealAlgebraicKnots2016}. But the authors told us it was indeed a minor mistake they were already aware of.
		Note that, for degree $6$, only knots in the class $h_{5,1}$ are hyperboloidal, the others being merely topologically hyperboloidal.
		\item In \cite{mikhalkinTopologyMaximallyWrithed2018}, it is proved that real algebraic knots of degree $d$ maximising \textit{encomplexed writhe}, an invariant introduced by Viro in \cite{viroEncomplexingWrithe2000}, are rational and of type $h_{d-1,1}$. Knots and links maximising this invariant are further studied in the series of article \cite{mikhalkinMaximallyWrithedReal2019a}, \cite{mikhalkinRigidIsotopyMaximally2021} and \cite{mikhalkinOsculatingFramingReal2019}.
	\end{enumerate}
\end{rem}

\subsection{Topology of $\Bl_c\P_\R^3$ when $\deg c\leq 4$.}

We apply the topological computations of the preceding sections to real forms of the Fano threefolds in Mori--Mukai Families \textnumero~2.33, 2.30, 2.27 and 2.22. Their complex members are obtained
by blowing up $\P^3_\C$ along a smooth rational curve of degree $1$, $2$, $3$ and $4$, respectively. We describe their real blow-up models and determine the diffeomorphism types of their real loci. We restrict attention to the forms with points, since we're primarily interested in the topologies.

We write
\[\pazocal Q_{r,s}:=\{x_0^2+\dots+x_{r-1}^2-x_r^2-\dots-x_{r+s-1}^2=0\}\subset \P^{r+s-1}.\] 

\begin{lem}\label{lem:RealLocusRealQuadricHypersurfaces}
	For $r,s\geq 1$,
	\[ \pazocal Q_{r,s}(\R)\approx	\faktor{\S^{r-1}\times \S^{s-1}}{\mu_2},\]	
	where $\mu_2=\{-1,1\}$ acts diagonally.
\end{lem}
\begin{proof}
	It suffices to note that 
	\[ \pazocal Q_{r,s}(\R)= \left\{[x:y]\in \P^{r+s-1}(\R)\mid |x|^2=|y|^2\right\}, \]
	where $x\in \R^r$ and $y\in \R^s$. Every real point has a representative satisfying
	$|x|=|y|=1$. And any two representative differs by the diagonal action of an element of $\mu_2$. 
\end{proof}

\begin{thm}\label{Th:BlowUpProjSpaceAlongRatCurves}
	Let $X$ be a smooth Fano threefold over $\R$ with $X(\R)\neq\emptyset$. Suppose that $X_\C$ belongs to one of the following Mori--Mukai families.
	\begin{enumerate}
		\item \textbf{Family \textnumero~2.33.}
		The variety $X$ is the blow-up of $\P^3_\R$ along a real line. All such varieties are $\R$-isomorphic, and
		\[X(\R) \approx \RP^2\times\S^1.\]
		\item \textbf{Family \textnumero~2.30.}
		There are exactly two possibilities up to $\R$-isomorphism:
		\begin{itemize}
			\item the blow-up of $\P^3_\R$ along a smooth real conic with nonempty real locus, equivalently the blow-up of $\pazocal Q_{3,2}$ at a real point, with
			\[X(\R) \approx \RP^3 \# \, (\S^2\widetilde{\times}\S^1);\]
			\item the blow-up of $\P^3_\R$ along a smooth real conic with empty real locus, equivalently the blow-up of $\pazocal Q_{4,1}$ at a real point, with
			\[ X(\R) \approx \RP^3.\]
		\end{itemize}
		\item \textbf{Family \textnumero~2.27.}
		The variety $X$ is the blow-up of $\P^3_\R$ along a real rational normal cubic. All such varieties are $\R$-isomorphic, and
		\[X(\R) \approx \RP^2\times\S^1.\]
		\item \textbf{Family \textnumero~2.22.}
		The variety $X$ is the blow-up of $\P^3_\R$ along a smooth geometrically rational quartic curve $c$. Its real locus is
		\[X(\R) \approx
		\begin{cases}
			\RP^3,
			& c(\R)=\emptyset,\\[2pt]
			\RP^3\#\, \S^2\widetilde{\times}\S^1,
			& c(\R)\in\langle\mu_0\rangle,\\[2pt]
			M\bigl(\RP^2;0;\overline{0};(2,1),(4,1)\bigr),
			& c(\R)\in h_{3,1}.
		\end{cases}\]
		These exhaust the possibilities, and all three occur.
	\end{enumerate}
\end{thm}

\begin{proof}
	By \cite[Lemma~2.5]{abbanKstabilityPointlessPezzo2024}, the extremal birational contraction of $X_\C$ to $\P^3_\C$ descends to a morphism $X\to V$ over $\R$, where $V$ is a real form of $\P^3_\C$. Since $X(\R)\neq\emptyset$, we have $V(\R)\neq\emptyset$. A Severi--Brauer variety with a real point is split, so $V\simeq\P^3_\R$. Thus $X$ is the blow-up of $\P^3_\R$ along a smooth geometrically rational real curve $c$, of degree $1$, $2$, $3$ or $4$, respectively.
	
	The group $\mathrm{PGL}_4(\R)$ acts transitively on real lines. It also acts transitively on each of the two classes of smooth real conics, distinguished by whether their real locus is empty: this follows by taking their spanning planes and applying the classification of nondegenerate real quadratic forms. The alternative descriptions in Family \textnumero~2.30 as blow-ups of quadric threefolds at real points are given by \Cref{prop:BlowUpPointQuadricThreefolds}.
	
	If $c$ is a real rational normal cubic, the real line bundle $\O_c(1)$ has odd degree. Hence $c$ has a real point and is isomorphic to $\P^1_\R$. Its embedding is given
	by the complete linear system $|\O_{\P^1_\R}(3)|$, so all such cubics are projectively equivalent over $\R$. These observations prove the asserted uniqueness statements.
	
	If $c(\R)=\emptyset$, the blow-down induces a diffeomorphism $X(\R) \approx \RP^3$. Otherwise,
	\[X(\R) \approx \Bl^{\mathrm{top}}_{c(\R)}\RP^3.\]
	By \Cref{table:KnotsFromRationalCurves}, the real loci of lines and rational normal cubics belong to $\langle\mu_1\rangle$, while a smooth real conic with real points represents $\langle\mu_0\rangle$. For a smooth geometrically rational quartic with nonempty real locus, the possibilities are $\langle\mu_0\rangle$ and $h_{3,1}$. Applying \Cref{Th:TopBlowUpHyperboloidalLink} gives the stated diffeomorphism types. In particular, for $h_{3,1}$ the normalised Seifert invariants are
	\[ \bigl(\RP^2;0;\overline{0};(2,1),(4,1)\bigr).\]
	The two conic models have different real loci and are therefore not $\R$-isomorphic.
	
	The two nonempty quartic knot types are realised as in \Cref{table:KnotsFromRationalCurves}, see also the discussion below. To realise the empty	quartic case, recall that the product of two anisotropic real conics $\pazocal Q_{3,0}$ identifies with the definite quadric surface $\pazocal Q_{4,0}$ in $\P^3_\R$. On $(\pazocal Q_{3,0})_\C\simeq\P^1_\C$, the morphism
	\[[s:t] \longmapsto [s^3:t^3]\]
	commutes with the real structure $[s:t] \mapsto [-\overline t:\overline s]$. It therefore descends to a degree-$3$ morphism $f\colon \pazocal Q_{3,0}\to \pazocal Q_{3,0}$. Its graph in $\pazocal Q_{4,0}$ is a smooth geometrically rational curve of bidegree $(3,1)$, hence a quartic in $\P^3_\R$, with empty real locus.
\end{proof}

We conclude by distinguishing the two diffeomorphism types obtained by blowing up $\P^3_\R$ along an $\R$-rational quartic. The complex geometry of this family is described in
\cite{cheltsovKstableFanoThreefolds2022}. A smooth rational quartic $c$ lies on a unique smooth quadric surface $\P^1_{u:v}\times \P^1_{x:y}$, by \Cref{cor:QuadricSurfaceContainingQuartic}. When the curve has real points, this quadric is split over $\R$, and the curve is the graph of a degree-$3$ morphism between its two ruling spaces. With the right order on the factors, the first projection identifies $c$ with $\P^1_{u:v}$, while the second projection is the degree-$3$ morphism
\[ \psi\colon\P^1_{u:v}\longrightarrow\P^1_{x:y}.\]

\begin{prop}\label{prop:QuarticKnotDegreeCriterion}
	Choose orientations on the real source and target circles, and let
	\[ \delta = \deg \bigl(\psi_{\R}\colon\RP^1\to\RP^1\bigr).	\]
	Then $\delta\in\{-3,-1,1,3\}$, and
	\[ c(\R)\in\begin{cases}
		h_{3,1},&|\delta|=3,\\
		\langle\mu_0\rangle,&|\delta|=1.
	\end{cases} \]
	If $\psi$ has a real ramification point, then $|\delta|=1$. If it has no real ramification points, its restriction to the real locus is a covering with either one or three sheets; these two	cases give $\langle\mu_0\rangle$ and $h_{3,1}$, respectively.
\end{prop}

\begin{proof}
	The real curve is the graph of $\psi_\R$, so its homology class in the ruling basis is $(1,\delta)$. Since $\psi$ has algebraic degree $3$, its topological degree on
	the real locus is odd and has absolute value at most $3$. The hyperboloidal knot classification therefore gives the asserted alternatives.
		
	A non-real ramification point cannot map to a real point: its conjugate would be another ramification point in the same fibre, whose degree would then be at least $4$. If a ramification point is real, its local degree is $2$ or $3$. In the first case, on one side of its branch value a nearby real fibre has only one real point; the same holds on both sides in
	the second case. This gives $|\delta|=1$ in any case. If there is no real ramification point, there is no real branch value either. The restriction $\psi_\R$ is therefore a covering, and its number of sheets is $|\delta|\in\{1,3\}$.
\end{proof}

Denoting $\operatorname{Rat}_3(\R)$ the space of quotient of real degree-$3$ polynomial functions on $\R^2$, the graph construction identifies real projective equivalence classes of these quartics with
\[\faktor{\operatorname{Rat}_3(\R)}{\bigl(\mathrm{PGL}_2(\R)\times\mathrm{PGL}_2(\R)\bigr)},\]
where the two groups change source and target coordinates independently. We equip this set with the quotient topology. The space $\operatorname{Rat}_3(\R)$ has four connected components,
distinguished by $\delta=-3,-1,1,3$; see \cite[Section~7]{segalTopologySpacesRational1979}. Coordinate changes identify the components with opposite degrees, so the quotient has exactly two connected components, distinguished by $|\delta|=1$ and $|\delta|=3$. 

Representatives of these two components are given by the graphs $\Gamma_{4,1}$ and $\Gamma_{4,3}$ of 
\[\psi_1\colon [u:v]\mapsto [u^3:v^3] \qquad \textnormal{and} \qquad \psi_{3}\colon [u:v]\mapsto [u^3-3uv^2:3u^2v-v^3].\]
Indeed, the map $\psi_1$ is an orientation-preserving homeomorphism of $\P^1(\R)$, hence of topological degree $1$. Its graph represents the class $(1,1)$ in the ruling basis and $\Gamma_{4,1}(\R)$ is of type $\langle \mu_0\rangle$. Writing $[u:v]=[\cos \theta : \sin \theta]$, the triple-angle identities give
\[\psi_3([\cos\theta : \sin\theta]) = [\cos 3\theta : \sin 3\theta].\] 
Its topological degree is $3$ so $\Gamma_{4,3}(\R)$ represents $h_{3,1}$.

\section{Some Fano blow-ups of quadric threefolds.}

\subsection{Quadric threefolds of signature $(3,2)$.}

Let
\[Q=\pazocal Q_{3,2}:=\{x_0^2+x_1^2+x_2^2-x_3^2-x_4^2=0\}\subset\P^4,\qquad Q_0:=Q\cap\{x_2=0\}. \]
We identify $Q_0$ with $\P^1\times\P^1$ through the Segre parametrisation
\[\begin{aligned}
	\Phi:\P^1_{u:v}\times\P^1_{x:y}&\longrightarrow Q_0,\\
	([u:v],[x:y]) &\longmapsto [ux+vy:uy-vx:0:ux-vy:uy+vx].
\end{aligned}\]
Write $\pr_1$ and $\pr_2$ for the two projections. On each real projective line, we use the orientation given by $\theta\mapsto[\cos\theta:\sin\theta]$, with $\theta\in\R/\pi\Z$. This gives a (algebraically) natural basis for $H_1(Q_0(\R);\Z)$ given by the (real loci of the) rulings. The orientation are chosen such that the projections have topological degree $1$ on the corresponding rulings.

\begin{prop}\label{prop:ClassificationCurvesQ32}	
	Let $c\subset Q$ be a smooth geometrically rational curve defined over $\R$, of degree $d\in\{1,2,3\}$. Then an $\R$-automorphism of $Q$ carries $c$ to one of the following
	normal forms.	
	\begin{enumerate}
		\item[\emph{($d=1$)}] There is a single normal form:
		\[R=\{x_2=0,\ x_0=x_3,\ x_1=x_4\}.\]
		Under $\Phi$, this is the fibre $\{[1:0]\}\times\P^1$ of $\pr_1$.
		\item[\emph{($d=2$)}] There are three normal forms, distinguished by the signature of the ambient quadratic form restricted to the vector subspace defining their spanning plane:
		\[\begin{array}{c|c|c}
			\text{Curve} & \text{Spanning plane} & \text{Relative signature}\\
			\hline
			C_{2,+} & x_2=x_4=0 & (2,1)\\
			C_{2,-} & x_1=x_2=0 & (1,2)\\
			\Gamma & x_3=x_4=0 & (3,0).
		\end{array}\]
		Each curve is the intersection of its spanning plane with $Q$. The first two lie in $Q_0$. 
		\item[\emph{($d=3$)}] There are three normal forms contained in $Q_0$ and one contained in a singular hyperplane section. The first three are most conveniently described as graphs. For $\varepsilon\in\{0,2,-2\}$, set	
		\[\nu_\varepsilon([s:t]):=\Phi\big(f_\varepsilon([s:t]),[s:t]\big),	\qquad C_{3,\varepsilon}:=\nu_\varepsilon(\P^1),\]
		where
		\[\begin{array}{c|c|c}
			\varepsilon & f_\varepsilon([s:t]) & \deg\big(f_\varepsilon|_{\P^1(\R)}\big)\\
			\hline
			0  & [s^2:t^2]       & 0\\
			2  & [s^2-t^2:2st]   & 2\\
			-2 & [s^2-t^2:-2st]  & -2.
		\end{array}\]
		Thus the index $\varepsilon$ records the topological degree of the first projection on the real locus. All three curves have bidegree $(1,2)$. In the ambient coordinates, these parametrisations are
		\[\begin{aligned}
			\nu_0([s:t]) &= [s^3+t^3:s^2t-st^2:0:s^3-t^3:s^2t+st^2],\\
			\nu_2([s:t]) &=	[s^3+st^2:-s^2t-t^3:0:s^3-3st^2:3s^2t-t^3],\\
			\nu_{-2}([s:t])	&= [s^3-3st^2:3s^2t-t^3:0:s^3+st^2:-s^2t-t^3].
		\end{aligned}\]
		The remaining normal form lies in the quadric cone $S:=Q\cap\{x_0-x_3=0\}$ and is given by
		\[ C_{3,s} = S\cap \left\{x_0x_2 -(x_1+x_4)^2 = x_0(x_4-x_1) - x_2(x_4+x_1)= 0\right\}\]
		A parametrisation is
		\[\begin{aligned}
			\nu_s:\P^1&\longrightarrow C_{3,s},\\
			[s:t] &\longmapsto [2s^3:t(s^2-t^2):2st^2:2s^3:t(s^2+t^2)].
		\end{aligned}\]
		This curve passes through the vertex $[1:0:0:1:0]$ of $S$. Its strict transform on the
		minimal resolution of $S$ is a section of the ruling, as described in \Cref{prop:RatCurvesOnQuadricCones}.
	\end{enumerate}
\end{prop}

\begin{thm}\label{Th:BlowUpQuadricThreefold32AlongRatCurves}
	Let $Q$ be a smooth real quadric hypersurface of $\P^4$ of signature $(3,2)$ and let $c$ be a smooth real geometrically rational curve in $Q$ with degree $d$. The real locus of $\Bl_c Q$ has the following description:
	\begin{enumerate}
		\item[$(d = 1):$] $\Bl_cQ(\R)\approx \mathbb{RP}^2\times \S^1$.
		\item[$(d = 2):$] \begin{itemize}
			\item if $c\simeq C_{2,+}$ in $Q$, then $\Bl_c Q(\R)\approx \S^2\widetilde{\times}\S^1 \#\, \S^2\widetilde{\times}\S^1$.
			\item if $c\simeq C_{2,-}$ in $Q$, then $\Bl_c Q(\R)\approx M\big(\mathrm{Kl}^2; 0; \bar{1}\big)\approx M_{A_{2,-}}$ where the last manifold is the mapping torus associated to $A_{2,-} = \left( \begin{smallmatrix}
				-1 & 0\\ 1 & 1
			\end{smallmatrix}\right)$.
			\item if $c\simeq \Gamma$ in $Q$, then $\Bl_c Q(\R)\approx \S^2\widetilde{\times}\S^1$.
		\end{itemize}
		\item[$(d = 3):$] \begin{itemize}
			\item if $c\simeq C_{3,0}$ in $Q$, then $\Bl_cQ(\R)\approx \RP^2\times \S^1$.
			\item if $c\simeq C_{3,-2}$ in $Q$, then $\Bl_cQ(\R)\approx \RP^2\times \S^1$.
			\item if $c\simeq C_{3,2}$ in $Q$, then $\Bl_cQ(\R)\approx M\left(\overline{\mathbb D^2};0;0;(2,1),(3,1)\right)\approx M_{\widetilde{g}}$ where $\widetilde{g}$ is an order-$6$ diffeomorphism of $N_3:=\#^3\RP^2$ and $M_{\widetilde{g}}$ denotes its mapping torus. The orbifold $\overline{\mathbb D^2}$ is a closed $2$-disc whose boundary is a reflector circle.
			\item if $c\simeq C_{3,s}$ in $Q$, then $\Bl_cQ(\R)\approx \RP^2\times \S^1$.
		\end{itemize}
	\end{enumerate} 
\end{thm}

Writing
\[v_+ = (x_0,x_1,x_2)\in V_+ = \C^3,\qquad v_-= (x_3,x_4) \in V_- = \C^2,\]
the equation of $Q$, after rescaling, becomes simply $\|v_+\| = \|v_-\| = 1$ in $\P(V_+\oplus V_-)$.
The remaining ambiguity is simultaneous multiplication of $v_+$ and $v_-$ by $-1$, so that  
\[Q(\R)\approx \S^2\widetilde{\times}\S^1:=\faktor{\S^2\times \S^1}{\mu_2}\]
where $\mu_2=\{\pm 1\}$ acts diagonally (this is of course a particular case of \Cref{lem:RealLocusRealQuadricHypersurfaces}).
This non-orientable $3$-manifold has two natural fibration structures corresponding to the projection onto the two factors. 
In the first two paragraphs below, we give an in-depth description of these while relating them to some algebraic constructions. We then split the proofs of both \Cref{prop:ClassificationCurvesQ32} and \Cref{Th:BlowUpQuadricThreefold32AlongRatCurves} according to the degrees in the next three paragraphs.

\subsubsection{The circle bundle $p\colon Q(\R)\to \mathbb{RP}^2$.}

The map 
\[Q(\R)\longrightarrow \mathbb{RP}^2,\qquad [v_+,v_-]\longmapsto [v_+]\] 
is a circle bundle whose fibres may be coherently oriented, \textit{i.e.} its Stiefel--Whitney class is zero and it is principal. There are only two such circle bundle structures on $\mathbb{RP}^2$, classified by the binary Euler class $e\in H^2(\mathbb{RP}^2;\Z)\simeq \Z_2$.
The circle bundle above corresponds to the non-trivial one.

For the algebraic point of view, consider the linear projection
\[p\colon \P^4 = \P(V_+\oplus V_-)\longrightarrow \P^2 = \P(V_+), \qquad [x_0:\dots:x_4] \longmapsto [x_0:x_1:x_2].\]
Its centre is the line $L=\P(V_-)=\{v_+=0\} \subset \P^4$ whose intersection with $Q$ consists of the conjugate pair $p_\pm=[0:0:0:1:\pm \i]$.
In particular, the indeterminacy locus of the induced rational map on $Q$ has no real points and the corresponding map on the real locus is exactly the circle bundle map above. 
We will denote it by $p$ as well.

Note that, blowing up the two points $p_\pm$, the rational map $p$ resolves to a conic bundle 
\[\widehat p\colon \widehat{Q} = \Bl_{p_+,p_-}Q\longrightarrow \P^2.\]
Here is an intrinsic description of $\widehat Q$. 
Decompose the quadratic form defining $Q$ as $q_+-q_-$, with $q_+$ and $q_-$ positive definite on $V_+$ and $V_-$ respectively. 
The blow-up $\Bl_L \P^4$ consists of pairs $([v_+:v_-],[\ell_a])$ such that $v_+\in \ell_a\subset V_+$ where $[\ell_a]\in \P(V_+)$ denotes a one-dimensional subspace.

Fix $[\ell_a]\in \P(V_+)$.
The fibre above is simply $\P(\ell_a\oplus V_-)$. 
Now, as $\ell_a$ varies, the lines $\ell_a$ forms the tautological bundle $\O_{\P(V_+)}(-1)$ and thus
\[\Bl_L \P^4\simeq \P_{\P^2}(\O_{\P^2}(-1)\oplus \O_{\P^2}^{\oplus 2}).\]
More explicitly, if $[a_0:a_1:a_2]\in \P^2$ and $a=(a_1,a_2,a_3)$ is a representative of the corresponding line, then a point in the fibre can be represented by $[\lambda a_0:\lambda a_1:\lambda a_2:y_3:y_4]$ and the fibre coordinates are $[\lambda:y_3:y_4]\in \P^2$.
Hence $\widehat{Q}$ is the relative conic 
\[\widehat{Q} = \{\lambda^2(a_0^2+a_1^2+a_2^2)-y_3^2-y_4^2\}\subset \P_{\P^2}(\O_{\P^2}(-1)\oplus \O_{\P^2}^{\oplus 2}).\] 
Its discriminant is $\{q_+=0\}\subset \P^2$ so every fibre over $\R$ is a smooth real conic.

The restriction of the circle-bundle projection to $Q_0(\R)$ is 
\[p_0\colon Q_0(\R)\to\RP^1, \qquad [x]\longmapsto[x_0:x_1:0]\in \{x_2=0\}\subset \P^2.\]
In the chosen coordinates, $p_0([u:v],[x:y]) = [ux+vy:uy-vx]$. 
Choosing the parametrisation 
\[[u:v]=[\cos\alpha:\sin\alpha],\qquad [x:y]=[\cos\beta:\sin\beta],\]
then $p_0(\alpha,\beta) = [\cos(\alpha-\beta):-\sin(\alpha-\beta)]$.
Hence $p_0$ is again a circle bundle on $Q_0(\R)\approx \RP^1\times\RP^1$ whose fibres have slope $(1,1)$.
Note they are therefore distinct from the fibres of the two rulings, which have slopes $(0,1)$ and $(1,0)$. 
Algebraically, this is also clear.
A fibre of the rational map $Q_0\dashrightarrow\P^1$, $[x]\longmapsto[x_0:x_1]$ is a hyperplane section of $Q_0$, hence a curve of bidegree $(1,1)$. 

\subsubsection{The sphere bundle $r\colon Q(\R)\to \RP^1$.}

the second natural rational projection
\[r\colon Q \dashrightarrow \P(V_-)\simeq\P^1, \qquad [x_0:\cdots:x_4] \longmapsto [x_3:x_4]\]
is defined on the set of real points.
Its complex indeterminacy locus is, however, the smooth isotropic real conic
\[\Gamma = Q\cap\P(V_+) = \{x_3=x_4=0,\; x_0^2+x_1^2+x_2^2=0\} \subset\P^4.\]

We may give an analogue resolution of this rational map.
Consider first the ambient linear projection $\P(V_+\oplus V_-) \dashrightarrow \P(V_-)$.
Its centre is $\P(V_+)$.
Its graph gives an identification
\[\Bl_{\P(V_+)}\P^4 \simeq \P_{\P(V_-)}\left(V_+\otimes\O_{\P(V_-)}\oplus\O_{\P(V_-)}(-1)\right),\]
Indeed, if $[b_3:b_4]\in\P(V_-)$ represents the line $\ell_b\subset V_-$, then the fibre of the graph over $[b_3:b_4]$ is $\P(V_+ \oplus \ell_b)\simeq\P^3$.
Using relative homogeneous coordinates $[y_0:y_1:y_2:\mu]$ on this fibre, the corresponding strict transform $\widetilde Q$ of $Q$ is 
\[\widetilde Q = \left\{y_0^2+y_1^2+y_2^2 - \mu^2(b_3^2+b_4^2) = 0\right\} \subset \P_{\P^1}\left(\O_{\P^1}^{\oplus3}\oplus\O_{\P^1}(-1)\right).\]

The projection of the projective bundle therefore induces $\widetilde r\colon \widetilde Q\longrightarrow \P^1$.
The fibre over \(b=[b_3:b_4]\) is the quadric surface
\[F_b = \left\{y_0^2+y_1^2+y_2^2 - (b_3^2+b_4^2)\mu^2 = 0\right\} \subset \P^3.\]
Thus $\widetilde r$ is a quadric-surface, or degree-$8$ del Pezzo, fibration.

Its discriminant is 
\[\Delta_{\widetilde r} =\{q_-=0\}= \{b_3^2+b_4^2=0\} \subset \P^1.\]
It consists of the two conjugate points $[1:\i]$, $[1:-\i]$.
Since the centre $\Gamma$ has no real points, $\widetilde Q(\R)\approx Q(\R)$. 
For a real point $b\in\mathbb{RP}^1$, one has $b_3^2+b_4^2 > 0$, and the real fibre is a quadric of signature $(3,1)$, and hence $F_b(\mathbf R)\approx \S^2$.
Thus $r$ induces the $\S^2$-bundle structure on $Q(\R)$.

This also makes its monodromy transparent.
After normalising $b_3^2+b_4^2=1$, a choice of lift $b\in \S^1$ of $[b]\in\RP^1$ identifies the fibre with $\S^2$.
Going once around $\RP^1$ changes the lift $b$ to $-b$. 
Since $[y:b]=[-y:-b]$, the corresponding identification of the $\S^2$-fibre is $y\mapsto-y$.
Therefore it corresponds to the mapping torus of the antipodal map on the sphere.

The restriction of $r$ to $Q_0$ is, using the parametrisation,
\[r_0([u:v],[x:y]) = [ux-vy:uy+vx].\]
On real points, writing
\[ [a:b]=[\cos\alpha:\sin\alpha], \qquad [u:v]=[\cos\beta:\sin\beta],\]
we obtain $r_0(\alpha,\beta) = [\cos(\alpha+\beta):\sin(\alpha+\beta)]$.
This is a new circle bundle on $Q_0(\R)\approx \RP^1\times \RP^1$ whose fibres have diagonal slope $(1,-1)$, up to the choice of generators.
Algebraically, a fibre of $r_0$ is again a hyperplane section of $Q_0$ and hence has bidegree $(1,1)$. The corresponding pencil has two base points, namely
\[\Gamma\cap Q_0= \{x_2=x_3=x_4=0,\ x_0^2+x_1^2=0\}=\{r_+,r_-\},\]
where $r_\pm=[1:\pm \i:0:0:0]$.
Of course, the strict transform $\widetilde Q_0\simeq \Bl_{r_+,r_-}Q_0 \subset\widetilde Q$ is simply the hyperplane $\{y_2=0\}$. For each real $b$, $\widetilde Q_0(\R)\cap F_b(\R)$ is an equatorial circle in $F_b(\R)\approx \S^2$.
In this sense the diagonal circle fibration $r_0$ of $Q_0(\R)$ sits naturally inside the $\S^2$-fibration $r$ of $Q(\R)$.

\subsubsection{Case of degree $d=1$.}

\begin{proof}[Proof of the \Cref{prop:ClassificationCurvesQ32}, case $d=1$.]
	There is only one class of lines in $Q$. 
	Indeed, a real line $L\subset Q$ corresponds to a $2$-dimensional totally isotropic subspace $(\R^5,q)$ with $q$ the diagonal quadratic form of signature $(3,2)$. 
	Since $q$ has Witt index $2$, these are maximal isotropic subspaces and Witt's theorem shows that $\mathrm{O}(3,2)$, hence $\operatorname{Aut}_\R(Q)$, acts transitively on them.
\end{proof}

\begin{proof}[Proof of \Cref{Th:BlowUpQuadricThreefold32AlongRatCurves}, case $d=1$.]
	Recall that $R$ is a fibre of the first ruling so $R(\R)$ is oriented such that the second projection $Q_0\simeq \P^1\times \P^1$ is of topological degree $1$. Since the circle fibrations $p_0$ and $r_0$ are given topologically by
	\[p_0(\alpha,\beta)=\alpha-\beta, \qquad r_0(\alpha,\beta)=\alpha+\beta,\]
	their restrictions to $R(\R)$ have degree $\mp1$ and are therefore sections on $Q_0$. The fibres of $r_0$ being the equators of the $\S^2$-fibres of $r\colon Q(\R) \longrightarrow \RP^1$, $R(\R)$ is actually a section of the $\S^2$-bundle.
	
	We may give a quick and totally topological argument yielding a description of the resulting $3$-manifold.   
	We use the $\S^2$-bundle $r\colon Q(\R)\to \RP^1$.
	An analogous line of argument to \ref{prop:TopBlowUpAlongFibreInLocTrivialBundle} shows that the topological blow-up along a section can be performed fibrewise.
	In every $\S^2$-fibre one blows up one real point, and $\Bl^{\mathrm{top}}_{\{pt\}}\S^2 \approx \RP^2$ building an $\RP^2$-bundle structure over $\S^1$ to $\Bl_{R}Q(\R)$.
	Since every diffeomorphism of $\RP^2$ is isotopic to the identity, this bundle is smoothly trivial. 
	Hence
	\[\Bl_{R}Q(\R) \approx \RP^2\times \S^1.\]
	
	We then discuss a more explicit and algebraic proof.
	It has the advantages to give a clearer picture of the exceptional divisor, of the strict transform of $Q_0$ which is used in the proof of \Cref{Th:TopologyOfV5}.
	
	Recall the quadric-surface fibration obtained after blowing up the isotropic conic
	\[\Gamma=\{x_3=x_4=0,\ x_0^2+x_1^2+x_2^2=0\}.\]
	In relative coordinates it is given by
	\[\widetilde{Q}=\{y_0^2+y_1^2+y_2^2-\mu^2(b_3^2+b_4^2)=0\}\subset \P_{\P^1_{b_3:b_4}}\left(\O_{\P^1}^{\oplus3}\oplus\O_{\P^1}(-1)\right),\]
	and the strict transform of $R$ is the section $b=[b_3:b_4] \longmapsto [b_3:b_4:0:1]=:\sigma_b$.
	
	Projection from this section is given by 
	\[ [z_0:z_1:z_2] = [y_0-\mu b_3:y_1-\mu b_4:y_2]\]
	and is resolved by blowing up $R$.
	It yields a morphism
	\[ \Phi\colon \Bl_{R}\widetilde Q\longrightarrow \P^2\times\P^1.\]
	Since the equation becomes 
	\[z_0^2+z_1^2+z_2^2 + 2\mu(b_3z_0+b_4z_1)=0,\]
	the induced map on real loci is a diffeomorphism $\Bl_{R}Q(\R) \approx \RP^2\times\RP^1$.
	
	Alternatively, we blow-up each quadric surface $F_b$ in the real point $\sigma_b$. 
	These are isomorphic to $\P^2$ blown-up in a pair of conjugated points.
	
	Under the above diffeomorphism, the exceptional divisor consists of the tangent directions to the quadric fibres at the section and is therefore mapped onto
	\[ K:=\left\{([z_0:z_1:z_2],[b_3:b_4])\in\RP^2 \times \RP^1 \; \mid \; b_3z_0+b_4z_1=0\right\}.\]
	The fibres of the projection $K\to\RP^1$ are lines passing through $[0:0:1]\in \RP^2$. If $[b_3:b_4]=[\cos t:\sin t]$, the fibre is parametrised by $[u:v] \mapsto [-\sin t\, u:\cos t\, u: v]$.
	After one turn in the base, the parameter $[u:v]$ is changed to $[-u:v]$; hence the fibre sweeps a half turn and $K$ is the mapping torus of a reflection of $\RP^1$, and therefore a Klein bottle.
	Note that the strict transform of $Q_0$ is $\{z_2=0\}$ and describe on each fibre the line at infinity.
	We obtain a diffeomorphism of pairs
	\[\bigl(\Bl_{R}Q(\R),E_R(\R)\bigr) \approx \bigl(\RP^2\times \S^1, K\bigr).\]
\end{proof}

\subsubsection{Case of degree $d=2$.}

\begin{proof}[Proof of the \Cref{prop:ClassificationCurvesQ32}, case $d=2$.]
	Writing $Q = \{q=0\}$, Witt's extension theorem implies that $\operatorname{Aut}_\R Q$ acts transitively on plane sections corresponding to $\P(W)\subset \P^4$ whose induced quadratic form $q_W$ has the same signature.
	Since $q$ has signature $(3,2)$ this gives at most three orbits for smooth conics, with signature $(3,0)$, $(2,1)$ and $(1,2)$.
	The former has no real point and the two latter gives topologically different according to \Cref{Th:BlowUpQuadricThreefold32AlongRatCurves} so these three orbits are disjoint.
\end{proof}

\begin{proof}[Proof of \Cref{Th:BlowUpQuadricThreefold32AlongRatCurves}, case $d=2$.]
	Start with $C_{2,+}$. Its real locus is the equatorial circle in a $\S^2$-fibre of $r$ and is then an unknot in $Q(\R)$, thus \Cref{prop:TopBlowUpAlongUnknot} yields
	\[\Bl_{C_{2,+}} Q (\R) \approx \S^2\widetilde{\times}\S^1 \#\, \S^2\widetilde{\times}\S^1.\]
	
	Then, consider $C_{2,-}$. Its image through $p\colon Q\to \P^2$ is $[1:0:0]$ so $C_{2,-}(\R)$ is just a fibre of this $\S^1$-bundle. Using \Cref{prop:TopBlowUpAlongFibreInLocTrivialBundle} or \Cref{Th:TopBlowUpSeifertFibrationAlongFibre}, the blow-up is the principal and non-trivial circle bundle over the Klein bottle
	\[\Bl_{C_{2,+}} Q (\R)\approx M\big(\mathrm{Kl};0;\bar{1}\big). \]
	Writing the Klein bottle as a mapping torus of the reflection $x\mapsto -x$ of $\S^1=\R/\Z$, it can be viewed as the torus bundle over $\S^1$:
	\[M=\faktor{\S^1\times \S^1\times [0;1]}{(x,z,1)\sim (-x,z+x,0)}\]
	Thus, it is the mapping torus of 
	\[A=\begin{pmatrix}
		-1 & 0\\
		1 & 1
	\end{pmatrix}.\]
	Notice $A^2=I_2$ and $\det A=-1$.
	
	Finally, $\Gamma$ has no real points so $\Bl_{\Gamma} Q (\R)\approx Q(\R)\approx \S^2\widetilde{\times}\S^1$.
\end{proof}

\subsubsection{Case of degree $d=3$.}

\begin{proof}[Proof of the \Cref{prop:ClassificationCurvesQ32}, case $d=3$.]
	If $c\subset Q$ is a twisted cubic, write $H=\langle c\rangle \simeq \P^3$.
	There are already two fundamentally different possibilities.
	
	Suppose first that the quadric surface $Q\cap H$ is smooth. It must then be split over the reals $H\cap Q\simeq \P^1_\R\times \P^1_\R$ for non-split real quadric cannot contain a real curve of odd degree $3$.
	As explained in the above case, $\operatorname{Aut}_\R Q$ acts transitively on such hyperplane sections, so we may assume $c\subset Q_0$. 
	
	Now, the subgroup of $\operatorname{Aut}_\R Q_0$ that extends to $Q$ is slightly smaller than the full real automorphism group $(\operatorname{PGL}_2(\R) \times \operatorname{PGL}_2(\R)) \rtimes \Z_2$ of $\P_\R^1\times \P_\R^1$.
	Indeed, using our Segre parametrisation $\Phi$ we can associate to a point of $Q_0$:
	\[M(X) = \begin{pmatrix}
		x_0+x_3 & x_1+x_4 \\
		x_4-x_1 & x_0-x_3
	\end{pmatrix},\qquad \textnormal{equivalently}\qquad M = \begin{pmatrix}
	u \\ v
	\end{pmatrix}\begin{pmatrix}
	x & y
	\end{pmatrix}\]
	on $Q_0$. 
	The equation of $Q_0$ is given by $4\det M(x) = x_0^2+x_1^2-x_3^2-x_4^2=0$.
	Automorphism of the two factors amount to the coordinate change $M\mapsto AMB^t$ and then multiply the quadratic form $q_0$ of $Q_0$ by $\det A \det B$. 
	The key observation is that such an automorphism extends to real automorphism of 
	\[Q = \{q_0^2+x_2^2 = 0\}\]
	if and only if this multiplier is positive.
	This, together with the observation that the swapping $\tau\colon ([u:v],[x:y]) \mapsto ([x:y],[u:v])$ corresponds to $x_1\mapsto -x_1$, shows that the subgroup of $\operatorname{Aut}_\R Q_0$ extending to elements of $\operatorname{Aut}_\R Q$ is the index $2$ subgroup of automorphisms preserving or reversing both the orientations on the two rulings.

	This gives at once that we may assume $c\in \O_{\P^1\times \P^1}(2,1)$ and that the signed integer $\deg_{\mathrm{top}}\pr_2|_{C(\R)}$ is an invariant, hence $C_{3,0}$, $C_{3,2}$ and $C_{3,-2}$ are not equivalent.
	
	The transitivity proof is elementary. 
	One may write
	\[c = (f(t),t),\quad t\in \P^1 \qquad \textnormal{where}\qquad  f:\P^1\longrightarrow \P^1\]
	has algebraic degree $2$.
	The curves $C_{3,0}$, $C_{3,2}$ and $C_{3,-2}$ correspond respectively to \[f_0\colon [s:t]\longmapsto [s^2:t^2],\quad f_2\colon [s:t]\longmapsto [s^2-t^2:2st], \quad \textnormal{and} \quad f_{-2}\colon [s:t] \longmapsto [s^2-t^2:-2st].\]
	If the ramification points of $f$ are real, send them by a real projective transformation of the source to $0$, $\infty$, and its two branch values by a projective transformation of the target to $0$, $\infty$.
	Then $f(z)=kz^2$ and after rescaling $f=f_0$.
	If the source and target normalising transformations initially have opposite orientation characters, compose the source with $z\mapsto -z$ which fixes $z^2$ but reverses the orientation character. 
	Thus the normalisation can be achieved inside the subgroup extending to $Q$.
	
	If the ramification points are conjugate nonreal, send them to $\{\pm \i\}$, and similarly send the branch points to $\{\pm \i\}$.
	The resulting map is either $f_2$ or $f_{-2}$.	
	Again composing the source with $z\mapsto -z$ if necessary, this may be achieved with elements of $\operatorname{Aut}_\R Q$.\\
	
	Now suppose the hyperplane section $Q\cap H$ is singular. Such a section of $Q$ is a tangent hyperplane section. By Witt's extension theorem, $\operatorname{Aut}_{\R}(Q)$ acts transitively on the real points of $Q$, hence on the set of its real tangent hyperplanes. We may therefore fix $v=[1:0:0:1:0]\in Q$ and consider the tangent hyperplane $H=T_p Q=\{x_0=x_3\}$.
	
	The corresponding hyperplane section is the quadric cone
	\[Q_v = Q \cap H =	\left\{x_0=x_3,\; x_1^2+x_2^2-x_4^2=0\right\},\]
	whose vertex is $v$. Let $c\subset Q_v$ be a smooth twisted cubic. By \Cref{prop:RatCurvesOnQuadricCones}, $c$ passes through the vertex $v$ and its strict transform $\widetilde{c}$ is a section of the minimal resolution $\mathbb{F}_2\to Q_v$. Equivalently, projection from the vertex $v$ identifies $c$ with the base conic
	\[B=\{x_1^2+x_2^2-x_4^2=0\} \subset \P^2_{x_1:x_2:x_4}.\]	
	The group $\operatorname{Aut}_{\R}(B) \simeq \operatorname{PGL}_2(\R)$ is induced by the subgroup of $\operatorname{Aut}_{\R}(Q)$ acting on $(x_1,x_2,x_4)$ and fixing $(x_0,x_3)$. After applying such an automorphism and reparametrising $c$, we may therefore assume that the projection of $c$ onto $B$ is the standard parametrisation
	\[[s:t]\longmapsto	[s^2-t^2:2st:s^2+t^2].\]	
	Choosing $[0:1]$ as the point mapping to the vertex, $c$ consequently has a parametrisation of the form
	\[[s:t] \longmapsto \left[F(s,t):\, s(s^2-t^2):\, 2s^2t:\, F(s,t):\, s(s^2+t^2) \right],\]
	where $F$ is a homogeneous cubic satisfying $F(0,1) \neq 0$.
	Thus $F(s,t)=a t^3+sG(s,t)$ for some $a\neq 0$ and quadratic form $G$.
	
	It remains to eliminate the term $sG$. For this, observe that the stabiliser of $Q_v$ in $\operatorname{Aut}_{\R}(Q)$ contains transformations whose restriction to $Q_v$ is
	\[[x_0:x_1:x_2:x_0:x_4]	\longmapsto [x_0+L(x_1,x_2,x_4):\, x_1:\, x_2:\, x_0+L(x_1,x_2,x_4):\, x_4],\]	
	where $L$ is an arbitrary real linear form. Indeed, if $\delta=x_0-x_3$	and $(a_1,a_2,a_4)\in\R^3$, the linear transformation defined by
	\[\begin{aligned}
		x_1'&=x_1+a_1\delta,\\
		x_2'&=x_2+a_2\delta,\\
		x_4'&=x_4+a_4\delta,\\
		x_0'+x_3' &= x_0+x_3 -2(a_1x_1+a_2x_2-a_4x_4) - \delta(a_1^2+a_2^2-a_4^2),\\
		x_0'-x_3'&=\delta
	\end{aligned}\]	
	preserves $x_0^2+x_1^2+x_2^2-x_3^2-x_4^2$. On $Q_v$, where $\delta=0$, it induces the preceding transformation with $L=-(a_1x_1+a_2x_2-a_4x_4)$.	
	Along $c$, the three functions	
	\[	s(s^2-t^2),\qquad 2s^2t,\qquad s(s^2+t^2)\]	
	span the vector space $s H^0(\P^1,\O_{\P^1}(2))$. Hence we may choose $L$ so that
	\[ L\bigl(s(s^2-t^2),\, 2s^2t,\, s(s^2+t^2)\bigr)=-sG(s,t).\]	
	After applying the corresponding automorphism, the parametrisation becomes
	\[[s:t] \longmapsto \left[at^3:\, s(s^2-t^2):\, 2s^2t:\, at^3:\, s(s^2+t^2)	\right].\]	
	Finally, the subgroup of $\operatorname{Aut}{\R}(Q)$ which, in the coordinates $(x_0,x_3)$, rescales $x_0+x_3$ and $x_0-x_3$ by reciprocal nonzero factors acts on $Q_v$ by a nonzero rescaling of the common coordinate $x_0=x_3$. We may therefore normalise $a=1$.
	Thus every smooth twisted cubic contained in $Q_v$ is equivalent to
	\[[s:t] \longmapsto	\left[t^3:\, s(s^2-t^2):\, 2s^2t:\, t^3:\, s(s^2+t^2)\right]\]	
	and the result follows.
\end{proof}

\begin{lem}
	Consider the real loci of the curves above as unoriented knots in $Q(\R)$. Then $C_{3,0}(\R)$, $C_{3,-2}(\R)$ and $C_{3,s}(\R)$ are ambiently isotopic to $R(\R)$, which is a section of $r$. By contrast, $C_{3,2}(\R)$ is a three-section of $r$: the restriction $r|_{C_{3,2}(\R)}:	C_{3,2}(\R) \to \P^1(\R)$ is a covering of degree $3$.
\end{lem}

\begin{proof}
	Introduce angular coordinates on $Q_0(\R)$ by writing
	\[ [u:v]=[\cos\alpha:\sin\alpha], \qquad [x:y]=[\cos\beta:\sin\beta],\]
	with $\alpha,\beta\in\R/\pi\Z$. The formula for $\Phi$ gives 
	\[r_0(\alpha,\beta)=\alpha+\beta \pmod{\pi}.\]
	In particular, $R(\R)$, which is given by $\alpha=0$, is a section of $r$.

	This also describes the map on first homology. Since the fibres of $r$ are simply connected, $r$ induces an isomorphism on $H_1$, and the inclusion $\iota\colon Q_0(\R)\hookrightarrow Q(\R)$ induces
	\[ \iota_*\colon\Z^2\longrightarrow\Z, \qquad (m,n)\longmapsto m+n.\]
	Its kernel is therefore generated by $(1,-1)$. Geometrically, this class is represented by $C_{2,+}(\R)$, which bounds a disk in the sphere $r^{-1}([1:0])$.
	
	We now establish the isotopy statements using the explicit parametrisations. The curve $C_{3,0}(\R)$ is the graph of the degree-zero map
	\[f_0:\P^1(\R)\longrightarrow\P^1(\R),\qquad [s:t]\longmapsto[s^2:t^2].\]
	This map is homotopic to the constant map with value $[1:0]$. Taking graphs throughout such a homotopy gives an isotopy in $Q_0(\R)$ from $C_{3,0}(\R)$ to $R(\R)$.
	
	For $C_{3,-2}$, setting $[s:t]=[\cos\theta:\sin\theta]$ gives $r\circ\nu_{-2}([\cos\theta:\sin\theta]) =[\cos\theta:-\sin\theta]$. 
	Thus $r$ restricts to a homeomorphism on $C_{3,-2}(\R)$: this curve is itself a section. Any two sections of an $\S^2$-bundle over a circle are isotopic through sections, because $\S^2$ is connected and simply connected. Consequently, $C_{3,-2}(\R)$ is isotopic to $R(\R)$.
	
	The same argument applies to $C_{3,s}(\R)$. Indeed, $r\circ\nu_s([s:t])=[2s^3:t(s^2+t^2)]$.
	In the affine coordinate $z=t/s$, this map is $z\longmapsto\frac{z+z^3}{2}$.
	It is strictly increasing, maps $\R$ onto $\R$, and	extends to a homeomorphism of $\RP^1$ by sending infinity to infinity. Hence $C_{3,s}(\R)$ is a topological section of $r$, and is again isotopic to $R(\R)$.
	
	Finally, the parametrisation of $C_{3,2}$ gives $r\circ\nu_2([\cos\theta:\sin\theta]) = [\cos(3\theta):\sin(3\theta)]$. Thus the restriction of $r$ to $C_{3,2}(\R)$ is precisely the three-fold covering $\theta\mapsto3\theta$ of $\R/\pi\Z$.
\end{proof}

\begin{proof}[Proof of \Cref{Th:BlowUpQuadricThreefold32AlongRatCurves}, case $d=3$.]
	The above lemma settles the cases $C_{3,0}$, $C_{3,-2}$ and $C_{3,s}$ for they are equivalent to the case of the line $R$, yielding
	\[\Bl_{C_{3,0}} Q(\R)\approx \Bl_{C_{3,-2}} Q(\R)\approx \Bl_{C_{3,s}} Q(\R)\approx \RP^2\times \S^1.\]
	For $C_{3,2}$, it tells that each fibre contains three points of $C_{3,2}(\R)$. The same line of argument as for the line $R$ shows that the blow up is a $N_3=\#^3\RP^2$-bundle over $\RP^1$, that is the mapping torus of a diffeomorphism of $N_3$.
	
	To identify this diffeomorphism, choose a trivialisation of the original $\S^2$-bundle $r$ over an interval $[0,\pi]$. Its endpoint gluing is the antipodal map $A:y\mapsto -y$. The three points of $C_{3,2}$ in the fibre can be chosen, on the equator, at angles
	\[\phi_k(t) = \frac{t}{3} + \frac{4k\pi}{3},\qquad k\in \{0,1,2\},\quad t\in [0,\pi].\]
	Thus they move rigidly by a rotation through $t/3$ and at the end of the interval, the rotation is $R_{\frac{\pi}{3}}$. Combining this motion with the antipodal endpoint gluing gives the diffeomorphism $g= A\circ R_{\frac{\pi}{3}}$ of $\S^2$. It preserves the set of the three initial points and cyclically permutes them. Hence it lifts to a diffeomorphism $\widetilde{g}$ of $N_3$ after blowing them up, and
	\[\operatorname{Bl}_{C_{3,2}}Q(\R)\approx M_{\widetilde{g}}\]
	with $\widetilde{g}$ of order $6$. 
	
	To get the claimed Seifert structure, observe that the suspension circles give a Seifert foliation. Having a periodic monodromy, its base orbifold is simply $B=\faktor{N_3}{\langle \widetilde{g}\rangle}$. In our explicit model, the order-$6$ map on $\S^2$ was
	\[g(\theta,z) = (\theta +\frac{4\pi}{3},-z)\]
	in coordinates where $\{\theta=0\}$ describes the equatorial circle.
	
	Consider one of the point $b$ to blow-up. The action of $g^3$ on $\P(T_b\S^2) \approx \mathbb{RP}^1$ is $[u:v]\mapsto [u:-v]$ and it has two fixed directions. The direction $[1:0]$ is tangent to the equator; in a blow-up chart around it the involution still looks like $(u,w)\mapsto (u,-w)$ so its fixed locus continues the strict transform of the equator. This produces the reflector boundary of the base orbifold.
	
	At the other fixed direction $[0:1]$, use the chart $u=vw$ so that the action becomes $(w,v)\mapsto (-w,-v)$ Thus the lifted involution acts locally by $-\operatorname{id}$ yielding an order-$2$ cone point, hence an exceptional fibre of type $(2,1)$.
	
	Finally the north and south poles form one orbit. Their stabilizer is generated by, $g^2$, of order $3$, and locally $g^2$ acts by rotation of order $\frac{2\pi}{3}$, yielding the second cone point, hence a exceptional fibre of type $(3,1)$.	
\end{proof}

\begin{rem}
	For an alternative proof of $\Bl_{C_{3,2}}Q(\R)\approx M\big(\overline{\mathbb D^2};0;0;(2,1),(3,1)\big)$, see \Cref{rem:SarkisovLinkTwistedV5}.
\end{rem}

\subsection{Quadric threefolds of signature $(4,1)$.}

Throughout this subsection, we write $\pazocal Q=\pazocal Q_{4,1}$. Its real locus is diffeomorphic to $\S^3$.

\begin{thm}\label{Th:BlowUpQuadricThreefold41AlongRatCurves}
	\begin{enumerate}
		\item Any two smooth real conics in $\pazocal Q$ with nonempty real loci are equivalent under a real automorphism of $\pazocal Q$. Their real loci are unknots in $\pazocal Q(\R)$. In particular, for every such conic $c_2$,
		\[\bigl(\Bl_{c_2}\pazocal Q\bigr)(\R) \approx \S^2\widetilde{\times}\S^1.\]
		\item The real locus of every real twisted quartic $c_4\subset\pazocal Q$ with $c_4(\R)\ne\emptyset$ is an unknot in $\pazocal Q(\R)$. In particular,
		\[ \bigl(\Bl_{c_4}\pazocal Q\bigr)(\R) \approx \S^2\widetilde{\times}\S^1.\]
	\end{enumerate}
\end{thm}

\begin{rem}
	The equivalence assertion in the first part of the theorem is algebraic. The second part asserts only that the real loci of the twisted quartics are ambiently isotopic: the quartics themselves need not be equivalent under a real automorphism of $\pazocal Q$.
\end{rem}

We begin by recalling the projection from a real point of $\pazocal Q$. The resulting Sarkisov link will be used in the proofs of both parts of the theorem.

\subsubsection{The projection from a point $\varphi_p\colon \pazocal Q(\R)\setminus \{p\}\to \R^3$.}

Fix a point $p\in\pazocal Q(\R)$. After a real projective change of coordinates, we may write
\[\pazocal Q = \{x_0x_4+x_1^2+x_2^2+x_3^2=0\}\subset\P^4, \qquad p=[1:0:0:0:0].\]
The projection from $p$ is the birational map
\[\varphi_p:\pazocal Q \dashrightarrow \P^3, \qquad [x_0:x_1:x_2:x_3:x_4] \longmapsto [x_1:x_2:x_3:x_4].\]
Writing $[u:v:w:t]$ for the coordinates on the target, its inverse is
\[ \varphi_p^{-1}: [u:v:w:t]\longmapsto [-(u^2+v^2+w^2):ut:vt:wt:t^2].\]
The base scheme of this inverse map is the smooth conic
\[ A=\{t=0,\ u^2+v^2+w^2=0\} \subset \Pi:=\{t=0\}\subset\P^3.\]
In particular, $A$ is defined over $\R$ but $A(\R)=\emptyset$.

Blowing up the base schemes of these two maps gives their common graph. Thus $\varphi_p$ fits into the diagram
\par\medskip
\begingroup
\definecolor{mapPi}{RGB}{214,92,92}
\definecolor{mapBeta}{RGB}{92,92,214}
\definecolor{mapPhi}{RGB}{153,92,214}
\noindent
\begin{minipage}[c]{0.72\linewidth}
	\centering
	\resizebox{\linewidth}{!}{%
		\begin{tikzcd}[
			ampersand replacement=\&,
			column sep=0.8em,
			row sep=4em
			]
			\& {E\simeq\P^2}
			\&\&\&\& X
			\&\&\&\& {\widetilde{\Sigma_p}}
			\\
			p
			\&\& {\Sigma_p}
			\&\& {\pazocal Q}
			\&\& {\P^3}
			\&\& \Pi
			\&\& A
			\arrow[hook, from=1-2, to=1-6]
			\arrow[color=mapPi, from=1-2, to=2-1]
			\arrow["{\pi=\Bl_p}"',
			color=mapPi, from=1-6, to=2-5]
			\arrow["{\beta=\Bl_A}",
			color=mapBeta, from=1-6, to=2-7]
			\arrow[hook', from=1-10, to=1-6]
			\arrow["{\P^1\text{-bundle}}",
			color=mapBeta, from=1-10, to=2-11]
			\arrow[hook, from=2-1, to=2-3]
			\arrow[hook, from=2-3, to=2-5]
			\arrow["{\varphi_p}"', dashed,
			color=mapPhi, from=2-5, to=2-7]
			\arrow[hook', from=2-9, to=2-7]
			\arrow[hook', from=2-11, to=2-9]
		\end{tikzcd}%
	}
\end{minipage}\hfill
\begin{minipage}[c]{0.25\linewidth}
	\centering
	\small
	\textbf{Isomorphisms}
	\[
	\begin{gathered}
		\textcolor{mapBeta}{
			E\xrightarrow[\sim]{\beta}\Pi
		}
		\\[1.4em]
		\textcolor{mapPi}{
			\widetilde{\Sigma_p}\setminus E
			\xrightarrow[\sim]{\pi}
			\Sigma_p\setminus\{p\}
		}
	\end{gathered}
	\]
\end{minipage}
\par
\endgroup
\medskip

where $\pi$ is the blow-up of $p$, $\beta$ is the blow-up of $A$ and 
\[\Sigma_p :=\pazocal Q\cap T_p\pazocal Q =\{x_4=0,\ x_1^2+x_2^2+x_3^2=0\}\]
is the tangent quadric cone at $p$. Here, $E$ and $\widetilde{\Sigma_p}$ are the exceptional divisors, with the following descriptions:
\begin{enumerate}
	\item The divisor $E\simeq\P^2_{\R}$ parametrises tangent directions to $\pazocal Q$ at $p$. The restriction $\beta|_E:E\xrightarrow{\sim}\Pi$ is an isomorphism. 	
	\item The exceptional divisor $\widetilde{\Sigma_p}$ of $\beta$ is the strict transform of $\Sigma_p$ under $\pi$. The morphism $\beta|_{\widetilde{\Sigma_p}}:\widetilde{\Sigma_p}\to A$ is a $\P^1$-bundle: its geometric fibres are the strict transforms of the lines on $\pazocal Q$ passing through $p$. More precisely,
	\[ \widetilde{\Sigma_p} \simeq \P_A\bigl(\O_A\oplus\O_A(1)\bigr), \qquad \widetilde{\Sigma_p}_{\C}\simeq\mathbb F_2,\]
	where $\O_A(1)$ is induced by the plane embedding of $A$. The intersection $E\cap \widetilde{\Sigma_p}$ maps isomorphically onto $A$ under $\beta$.
\end{enumerate}

Since $A(\R)=\emptyset$, we also have $\widetilde{\Sigma_p}(\R)=\emptyset$, and $\beta$ induces a diffeomorphism $\beta_{\R}:X(\R)\xrightarrow{\sim}\P^3(\R)$.
Under this identification, $E(\R)$ is the projective plane $\Pi(\R)$. Moreover, the algebraic isomorphism $\varphi_p: \pazocal Q\setminus\Sigma_p \xrightarrow{\sim} \P^3\setminus\Pi$ restricts on real loci to a diffeomorphism $\pazocal Q(\R)\setminus\{p\} \xrightarrow{\sim} \P^3(\R)\setminus\Pi(\R) \approx\R^3$. This is the stereographic description that we shall use below: projection from $p$ identifies the complement of $p$ in the real quadric with an affine three-space, while the exceptional divisor over $p$ records the plane at infinity.

\subsubsection{Case of degree $d=2$.}

Let $c_2\subset\pazocal Q$ be a smooth real conic with nonempty real locus. Its span is a plane $\P(W)\subset\P^4$, and $c_2=\pazocal Q\cap\P(W)$. The restriction to $W$ of the quadratic form defining $\pazocal Q$ is nondegenerate of signature $(2,1)$. Witt's extension theorem therefore implies that any two such planes, and hence any two such conics, are equivalent under a real automorphism of $\pazocal Q$.

Choose $p\in\pazocal Q(\R)\setminus\langle c_2\rangle$. Projection from $p$ restricts to a projective isomorphism from $\langle c_2\rangle$ onto a plane $\Lambda\subset\P^3$. The projected conic has no real point on $\Pi$, since
\[ \varphi_p:\pazocal Q(\R)\setminus\{p\} \xrightarrow{\sim}\P^3(\R)\setminus\Pi(\R).\]
It is therefore a compact conic in the affine plane $\Lambda(\R)\setminus\Pi(\R)\approx\R^2$, where it bounds a disk. The inverse projection carries this disk to an embedded disk in $\pazocal Q(\R)$ with boundary $c_2(\R)$. Then \Cref{prop:TopBlowUpAlongUnknot} completes the proof of the first part of \Cref{Th:BlowUpQuadricThreefold41AlongRatCurves}.

\begin{rem}
	The assertion $\Bl_{c_2} \pazocal Q(\R)$ may also be seen directly from the equation. 
	Indeed, we may choose coordinates in which
	\[\pazocal Q= \{x_0^2+x_1^2+x_2^2+x_3^2=x_4^2\}, \qquad c_2=\pazocal Q\cap\{x_2=x_3=0\}.\]
	The pencil of hyperplanes through $c_2$ induces a morphism
	\[ f:Y:=\Bl_{c_2}\pazocal Q\longrightarrow\P^1,\]
	resolving the rational map $[x_0:\cdots:x_4]\mapsto[x_2:x_3]$. Explicitly,
	\[ Y = \left\{ \bigl([x_0:\cdots:x_4],[s:t]\bigr) \in\pazocal Q\times\P^1 \ \middle|\ tx_2=sx_3 \right\}.\]
	Writing $x_2=sz$ and $x_3=tz$, the fibre over $[s:t]$ is the quadric surface
	\[\{x_0^2+x_1^2+(s^2+t^2)z^2=x_4^2\}\subset \P^3_{x_0:x_1:z:x_4}.\]
	Thus $f$ is a quadric-surface fibration. Its only singular geometric fibres lie over $[1:\i]$ and $[1:-\i]$; both are quadric cones. Every real fibre is smooth, of signature $(3,1)$, and has real locus diffeomorphic to $\S^2$.
	
	Every real point of $\pazocal Q$ has $x_4\ne0$, so we normalise $x_4=1$ and write
	\[ [s:t]=[\cos\theta:\sin\theta], \qquad 0\leq\theta\leq\pi.\]
	The fibre equation becomes $x_0^2+x_1^2+z^2=1$. At the two endpoints of this interval, the coordinates are identified by $(x_0,x_1,z,\pi) \sim (x_0,x_1,-z,0)$. Consequently,
	\[ Y(\R)\approx \frac{\S^2\times[0,\pi]}{(x_0,x_1,z,\pi)\sim(x_0,x_1,-z,0)} \approx \S^2\widetilde{\times}\S^1.\]
\end{rem}

\subsubsection{Case of degree $d=4$.}

Let $c_4\subset\pazocal Q$ be a real rational twisted quartic with $c_4(\R)\ne\emptyset$, and choose $p\in\pazocal Q(\R)\setminus c_4(\R)$. We retain the notation
\[\pazocal Q \xleftarrow{\ \pi\ }X \xrightarrow{\ \beta\ }\P^3\]
introduced for projection from $p$.

We first check that projection embeds $c_4$ in $\P^3$, even over $\C$. Otherwise, $p$ would lie on a secant or tangent line $\ell$ to $c_{4,\C}$ contained in $\pazocal Q_\C$. A rational normal quartic has at most one such line through a point outside the curve: distinct secant or tangent lines meeting away from the curve would contradict the linear independence of subschemes of length at most four on a rational normal quartic. Since $p$ is real, uniqueness would force $\ell$ to be defined over $\R$. Moreover, $\ell$ meets $\pazocal Q$ in $p$ and a subscheme of length two on $c_4$, so $\ell\subset\pazocal Q$. This is impossible, because a quadratic form of signature $(4,1)$ has Witt index one and admits no real projective line on its zero locus. Thus
\[ \Gamma:=\varphi_p(c_4)\subset\P^3 \]
is a smooth, nondegenerate rational quartic.

There is a unique real quadric surface $Q_0\simeq\P^1_{\R}\times\P^1_{\R} \subset\P_\R^3$ containing $\Gamma$, see \Cref{cor:QuadricSurfaceContainingQuartic}. After ordering the rulings, the curve $\Gamma$ has bidegree $(1,3)$. In particular, one ruling projection restricts to an isomorphism $\Gamma\to\P^1_{\R}$.

Consider the plane section
\[ B:=Q_0\cap\Pi. \]
It is a smooth real conic with nonempty real locus. Indeed, a real ruling line on $Q_0$ meets $\Pi$ in a real point, so $B(\R)$ is nonempty. If $\Pi$ were tangent to $Q_0$, its intersection with $Q_0$ would be the union of two real ruling lines. Their intersection degrees with $\Gamma$ are one and three, respectively, so they would meet $\Gamma$ in real points. This contradicts $\Gamma(\R)\cap\Pi(\R)=\emptyset$. The conic $B$ has bidegree $(1,1)$, and its real locus is an essential circle on the torus $Q_0(\R)$.

Since $A(\R)=\emptyset$ whereas $B(\R)\ne\emptyset$, the two plane conics $A$ and $B$ are distinct. Their intersection
\[ Z:=A\cap Q_0=A\cap B \]
is a scheme of length four without real points. In fact, $Z = \Gamma\cap\Pi$ scheme-theoretically. To see this, recall that $\Pi$ is sent to $p$ and $\Gamma$ to $c_4$ by $\varphi_p^{-1}$. But $c_4$ does not pass through $p$, so any intersection point $\Gamma \cap \Pi$ must lie on the indeterminacy locus $A$. Both $\Gamma\cap\Pi$ and $A\cap Q_0$ have length four, so they coincide.

Let $\widehat{Q_0}\subset X$ be the strict transform of $Q_0$ and set $T:=\pi(\widehat{Q_0})\subset\pazocal Q$. Then $\widehat{Q_0}\simeq\Bl_Z Q_0$ and the curve
\[ \widehat{B}:=\widehat{Q_0} \cap E \]
maps isomorphically to $B$ under $\beta$. The surface $\widehat{Q_0}$ is smooth along $\widehat{B}$, and
\[ \widehat{B}^2=B^2-\operatorname{length}(Z)=2-4=-2.\]
The morphism $\pi|_{\widehat{Q_0}}:\widehat{Q_0}\to T$ contracts precisely $\widehat{B}$ to $p$ and is an isomorphism away from $\widehat{B}$.

More explicitly, if $q_0(u,v,w,t)=0$ defines $Q_0$, then
\[ T=\pazocal Q\cap\{q_0(x_1,x_2,x_3,x_4)=0\}.\]
Thus $T$ is a complete intersection of two quadrics in $\P^4$. In the affine chart $x_0=1$ at $p$, the equation of $\pazocal Q$ gives $x_4=-(x_1^2+x_2^2+x_3^2)$ so that $(x_1,x_2,x_3)$ are local coordinates at $p$. With these coordinates, the resulting local equation of $T$ is $q_0(x_1,x_2,x_3,0)$ because terms involving $x_4$ are at least of degree $3$. But the projectivized zero locus of this equation describes also $B$. Since $B$ is smooth and has real points, $p$ is an ordinary double point with two real local sheets.

On real loci, $\beta$ identifies $\widehat{Q_0}(\R)$ with the torus $Q_0(\R)$, and $\pi$ collapses the single essential circle $\widehat{B}(\R)$. Consequently $T(\R)$ is a pinched torus, whose pinched point is $p$.

\par\medskip
\begingroup
\definecolor{mapPi}{RGB}{214,92,92}
\definecolor{mapBeta}{RGB}{92,92,214}
\definecolor{mapPhi}{RGB}{153,92,214}
\definecolor{mapSurface}{RGB}{204,51,204}

\noindent
\begin{minipage}[c]{0.72\linewidth}
	\centering
	\resizebox{\linewidth}{!}{%
		\begin{tikzcd}[ampersand replacement=\&,column sep=0.8em,row sep=1.8em]
			\& {E\simeq\P^2}
			\&\&\&\& X
			\&\&\&\& {\widetilde{\Sigma_p}}
			\\
			{\widehat{B}}
			\&\&\& {\widehat{Q_0}}
			\\
			p
			\&\& {\Sigma_p}
			\&\& {\pazocal Q}
			\&\& {\P^3}
			\&\& \Pi
			\&\& A
			\\
			\&\&\& T
			\&\&\&\& {Q_0}
			\&\& B
			\&\& Z
			\\
			\& {c_4}
			\&\&\&\&\&\&\& \Gamma
			\arrow[hook, from=1-2, to=1-6]
			\arrow[color=mapPi, from=1-2, to=3-1]
			\arrow["{\pi=\Bl_p}"',
			color=mapPi, from=1-6, to=3-5]
			\arrow["{\beta=\Bl_A}",
			color=mapBeta, from=1-6, to=3-7]
			\arrow[hook', from=1-10, to=1-6]
			\arrow["{\P^1\text{-bundle}}",
			color=mapBeta, from=1-10, to=3-11]
			\arrow[hook, from=2-1, to=1-2]
			\arrow[hook, from=2-1, to=2-4]
			\arrow[color=mapPi, from=2-1, to=3-1]
			\arrow[hook, from=2-4, to=1-6]
			\arrow[color=mapPi, from=2-4, to=4-4]
			\arrow[hook, from=3-1, to=3-3]
			\arrow[hook, from=3-1, to=4-4]
			\arrow[hook, from=3-3, to=3-5]
			\arrow["{\varphi_p}"', dashed,
			color=mapPhi, from=3-5, to=3-7]
			\arrow[hook', from=3-9, to=3-7]
			\arrow[hook', from=3-11, to=3-9]
			\arrow[hook, from=4-4, to=3-5]
			\arrow[dashed, color=mapSurface,
			from=4-4, to=4-8]
			\arrow[hook', from=4-8, to=3-7]
			\arrow[hook', from=4-10, to=4-8]
			\arrow[hook', from=4-12, to=3-11]
			\arrow[shift right, hook', from=4-12, to=4-10]
			\arrow[hook', from=4-12, to=5-9]
			\arrow[hook, from=5-2, to=4-4]
			\arrow["\sim"{description},
			color=mapPhi, from=5-2, to=5-9]
			\arrow[hook', from=5-9, to=4-8]
		\end{tikzcd}%
	}
\end{minipage}\hfill
\begin{minipage}[c]{0.25\linewidth}
	\centering
	\small
	\textbf{Isomorphisms}
	\[
	\begin{gathered}
		\textcolor{mapBeta}{
			E \xrightarrow[\beta]{\sim} \Pi
		}
		\\[1.2em]
		\textcolor{mapPi}{
			\widetilde{\Sigma_p}\setminus E
			\xrightarrow[\pi]{\sim}
			\Sigma_p\setminus \{p\}
		}
		\\[1.2em]
		\textcolor{mapBeta}{
			\widehat{B} \xrightarrow[\beta]{\sim} B
		}
	\end{gathered}
	\]
\end{minipage}
\par
\endgroup
\medskip

\begin{rem}
	We obtained that $c_4(\R)$ is a knot lying on a pinched torus $T(\R)$ and avoiding its pinched point. It would be enough to show that $T(\R)$ is unknotted. But we were unable to show this. We go another way. 
\end{rem}

Let $\widetilde c_4\subset X$ be the strict transform of $c_4$. It lies on $\widehat{Q_0}$ and is disjoint from $\widehat{B}$, since $p\notin c_4$. Under $\beta$, its real locus is identified with $\Gamma(\R)$. The degree-one ruling projection shows that this is an essential circle on $\widehat{Q_0}(\R)$. Two disjoint essential circles on a torus are parallel and bound an embedded annulus. Hence there is an annulus in $\widehat{Q_0}(\R)$ with boundary $\widetilde c_4(\R)\sqcup \widehat{B}(\R)$. Contracting its second boundary component gives an embedded topological disk in $T(\R)\subset\pazocal Q(\R)$ with boundary $c_4(\R)$. Thus $c_4(\R)$ is an unknot and
\[ \bigl(\Bl_{c_4}\pazocal Q\bigr)(\R) \approx \S^2\widetilde{\times}\S^1\]
by \Cref{prop:TopBlowUpAlongUnknot}. This completes the proof of \Cref{Th:BlowUpQuadricThreefold41AlongRatCurves}.

\begin{rem}
	The algebraic bidegree $(1,3)$ of $\Gamma$ should be distinguished from the winding numbers of $\Gamma(\R)$ on the real torus.	Here $\Gamma(\R)$ is disjoint from $B(\R)$, so the two circles are parallel. After choosing orientations of the two real ruling factors, their common slope is $(1,1)$ or $(1,-1)$.
\end{rem}

\begin{rem}\label{RemarkWeakDelPezzo}
	Suppose that $Z$ is reduced. Then it consists of four distinct geometric points, forming two conjugate pairs. They are in a special configuration: they all lie on the smooth $(1,1)$-curve $B$. The resulting blow-up $\widehat{Q_0}$ is therefore a smooth weak del Pezzo surface of degree four. Moreover $\pi|_{\widehat{Q_0}}$ is its anticanonical morphism. It contracts the unique $(-2)$-curve $\widehat{B}$ to the ordinary double point $p$ of the degree-four del Pezzo surface $T$, see \cite{demazureSurfacesPezzoII1980}.
\end{rem}

\section{Topology of real forms of the quintic del Pezzo threefold $V_5$.}

\subsection{Classical properties of $V_5$ and its two real forms.}

We recall the basic and classical geometry of $V_5$, with an emphasis on the descriptions that will be useful for studying its real topology. In particular, its realisation as an orbit closure will allow us to interpret its two real structures in terms of configurations of points on the projective line. This will be enough to obtain the topology for the non-trivial one. The Sarkisov link recalled below allows to compute the topology of the trivial one.

Up to isomorphism, there is a unique smooth complex Fano threefold of Picard rank one, index two and anticanonical degree $40$. It is called the \emph{quintic del Pezzo threefold} and is denoted by $V_5$; it belongs to family \textnumero~1.15 in the classification of smooth Fano threefolds.
Its Picard group is generated by a very ample divisor $H$ satisfying
\[-K_{V_5}=2H, \qquad H^3=5.\]
The complete linear system $|H|$ embeds $V_5$ into $\P^6$, where it can be realised as a smooth linear section
\[V_5\simeq \operatorname{Gr}(2,5)\cap\P^6\subset\P^9\]
of the Grassmannian in its Pl\"ucker embedding. In particular, its ideal in $\P^6$ is generated by the five quadrics obtained by restricting the Pl\"ucker relations. The threefold $V_5$ is rational, as can for instance be seen from the Sarkisov link recalled below.

\subsubsection{The Mukai--Umemura description as an orbit closure.}

A complementary description, due to Mukai--Umemura \cite{MukaiUmemura}, makes the symmetries of $V_5$ explicit. Let $R_6=\C[x,y]_6$ be the seven-dimensional space of binary sextics. The action
\[\begin{pmatrix}a&b\\ c&d\end{pmatrix}\cdot P(x,y)=P(ax+cy,bx+dy)\] 
of $\operatorname{GL}_2(\C)$ induces an action of $\operatorname{PGL}_2(\C)$ on $\P(R_6)\simeq\P^6$. It is enlightening to think of a projective binary sextic as its divisor of roots. More precisely, a factorisation
\[ P(x,y)=\lambda\prod_{i=1}^6(\beta_i x-\alpha_i y),\qquad \lambda\in\C^*,\]
associates to $[P]$ the effective divisor
\[ D_P=\sum_{i=1}^6[\alpha_i:\beta_i] \in\operatorname{Sym}^6(\P^1).\]
This gives an identification $\P(R_6)\simeq\operatorname{Sym}^6(\P^1)$.

Consider the sextic
\[ \varphi=xy(x^4+y^4).\]
Under stereographic projection, its six roots are the vertices of a regular octahedron on the Riemann sphere. The following description realises $V_5$ as a compactification of the space of configurations projectively equivalent to this octahedron.

\begin{prop}[Mukai--Umemura, {\cite{MukaiUmemura}}; see also {\cite[Proposition~1.1]{duboulozCylindersMoriFiber2019}}]
	There is an isomorphism
	\[V_5\simeq \overline{\operatorname{PGL}_2(\C)\cdot[\varphi]} \subset\P(R_6).\]
	The action has exactly three orbits:
	\begin{enumerate}
		\item the open orbit $\Omega_3=\operatorname{PGL}_2(\C)\cdot[\varphi]$, whose stabiliser is the octahedral group $O\simeq\mathfrak{S}_4$, of order $24$;		
		\item the two-dimensional orbit $\Omega_2=\operatorname{PGL}_2(\C)\cdot[xy^5]$, neither open nor closed, and whose stabiliser is the diagonal torus;		
		\item the closed one-dimensional orbit $\Omega_1=\operatorname{PGL}_2(\C)\cdot[x^6]$, whose stabiliser is the Borel subgroup of upper triangular matrices.
	\end{enumerate}
	Moreover, this action identifies $\operatorname{Aut}_{\C}(V_5)$ with $\operatorname{PGL}_2(\C)$.
\end{prop}

Here $O$ is the group of rotations of the octahedron. Its inverse image under $\operatorname{SL}_2(\C)\to\operatorname{PGL}_2(\C)$ is the
\emph{binary octahedral group} $2O$, of order $48$; this is the stabiliser of $[\varphi]$ for the $\operatorname{SL}_2(\C)$-action.

The two boundary orbits have a simple interpretation in terms of roots: $\Omega_2$ consists of divisors $5p+q$ with $p\ne q$, while $\Omega_1$
consists of divisors $6p$. In particular, $\Omega_1$ is a rational normal sextic curve, and
\[V_5\setminus\Omega_3=\overline{\Omega_2}=\Omega_2\sqcup\Omega_1.\]
Geometrically, this boundary is the surface swept out by the tangent lines to $\Omega_1$. It is an anticanonical divisor, singular along $\Omega_1$; see the discussion following \cite[Proposition~1.1]{duboulozCylindersMoriFiber2019}.

\subsubsection{The two real forms.}

The orbit description also provides concrete models of the real forms of $V_5$. Recall that, up to equivalence, the projective line has two
real structures:
\[c([\alpha:\beta])=[\overline{\alpha}:\overline{\beta}],\qquad a([x:y])=[-\overline{\beta}:\overline{\alpha}].\]
The fixed locus of $c$ is $\P^1(\R)$, whereas $a$ is the antipodal involution of the Riemann sphere and has no fixed points. On projective binary sextics, they induce the involutions 
\[\sigma_0([P])=[\overline P(x,y)],\qquad \sigma_1([P])=[\overline P(-y,x)],\]
where $\overline P$ denotes coefficientwise conjugation. Both fix $[\varphi]$ and normalise the $\operatorname{PGL}_2(\C)$-action, so they preserve its orbit closure $V_5$.

These give the two real forms of $V_5$, up to real isomorphism; see \cite[Proposition~10.4]{terpereauRealFormsMori2024}. The connection with
the two real structures on $\P^1$ is also reflected in the closed orbit: restriction to $\Omega_1\simeq\P^1$ identifies $\operatorname{Aut}_{\C}(V_5)$ with $\operatorname{Aut}_{\C}(\P^1)$.

\begin{dfn}
	We denote by $X_5$ the real form defined by $\sigma_0$ and by $Y_5$ the real form defined by $\sigma_1$. We call them the \emph{trivial} and \emph{non-trivial} real forms, respectively. Thus $X_5$ is obtained from the restriction of the real structure of $\P^6$ through its half-anticanonical embedding.
\end{dfn}

The projection from a line described below is defined over $\R$ so both $X_5$ and $Y_5$ are rational over $\R$.

\begin{lem}
	The groups of real automorphisms are
	\[\operatorname{Aut}_{\R}(X_5)\simeq\operatorname{PGL}_2(\R),	\qquad \operatorname{Aut}_{\R}(Y_5) \simeq\operatorname{PSU}_2 \simeq\operatorname{SO}_3(\R).\]
	For either real form, each of the three complex orbit strata is	defined over $\R$.
\end{lem}

\begin{proof}
	Real automorphisms are precisely the complex automorphisms commuting with the chosen real structure. In terms of the projective line, the centralisers of $c$ and $a$ in $\operatorname{PGL}_2(\C)$ are, respectively, $\operatorname{PGL}_2(\R)$ and $\operatorname{PSU}_2\simeq\operatorname{SO}_3(\R)$. Both real structures preserve the orbit strata, as is immediate from their descriptions by root divisors.
\end{proof}

The descent of this stratification should be distinguished from an orbit decomposition on real points: the real points of a complex
orbit need not form a single orbit under the real automorphism group. For instance, $[xy(x^4+y^4)]$ and $[xy(x^4-y^4)]$ belong to the same
complex open orbit, but to different $\operatorname{PGL}_2(\R)$-orbits, since they have respectively two and four real roots. For the non-trivial form, the antipodal description gives a direct identification of the entire real locus.

\begin{cor}\label{cor:TopologyTwistedV5}
	The two boundary strata of $Y_5$ have no real points, and $\operatorname{SO}_3(\R)$ acts transitively on $Y(\R)$. Consequently, there are diffeomorphisms
	\[Y_5(\R)\approx\faktor{\operatorname{SO}_3(\R)}{O} \approx\faktor{\operatorname{SU}_2}{2O} \approx \faktor{\S^3}{2O}.\]
	In particular, $Y_5(\R)$ is a closed orientable spherical three-manifold with fundamental group $2O$.
\end{cor}

\begin{proof}
	A point of $Y(\R)$ corresponds to a root divisor invariant under the antipodal involution $a$. Such a divisor cannot have the form $6p$ or $5p+q$ with $p\ne q$: in either case, invariance would force $p$ to be fixed by $a$. Thus every real point lies	in the open orbit. To prove transitivity, let $D_0$ be the root divisor of $\varphi$ and write an $a$-invariant divisor in the open orbit as $D=g(D_0)$, for a M\"obius transformation $g$. Then $b=g^{-1}ag$ is a fixed-point-free antiholomorphic involution preserving $D_0$. Since $a$ also preserves $D_0$, we can write $b=ha$ with $h\in O$. Rotations commute with the antipodal map, so $b^2=1$ implies $h^2=1$. If $h\ne1$, it is a rotation through angle $\pi$, and $ha$ is reflection in a great circle, which has fixed points. Hence $h=1$, and $g$ commutes with $a$.
	It therefore belongs to	$\operatorname{PSU}_2\simeq\operatorname{SO}_3(\R)$.
	
	This proves transitivity, with stabiliser $O$. Taking the inverse image of $O$ under the double covering $\operatorname{SU}_2\to\operatorname{SO}_3(\R)$ gives $2O$, and the remaining assertions follow from $\operatorname{SU}_2\approx S^3$.
\end{proof}

\subsubsection{The Hilbert scheme of lines.}

The geometry of lines provides another way to distinguish the two real forms of $V_5$. Throughout, a \emph{line} means a curve of degree one with respect to the half-anticanonical class $H$, and we denote its
Hilbert scheme by $\mathcal H(V_5)$. There are two possible normal bundles:
\[\NN_{\ell/V_5}\simeq\O_{\P^1}\oplus \O_{\P^1}\quad\text{or}\quad \NN_{\ell/V_5}\simeq\O_{\P^1}(-1)\oplus\O_{\P^1}(1).\]
We call the corresponding lines \emph{ordinary} and \emph{special}, respectively. In both cases $h^0(\NN_{\ell/V_5})=2$ and $h^1(\NN_{\ell/V_5})=0$, so the special lines also represent smooth points of the Hilbert scheme.

The classical description of this scheme is particularly simple:
\[\mathcal H(V_5)\simeq\P^2,\]
and the special lines are parametrised by a smooth conic $\mathcal C\subset\mathcal H(V_5)$. Equivalently, one may identify $\mathcal H(V_5)$ with the projective space of binary quadratics, with
$\mathcal C$ corresponding to the squares of linear forms. Thus the conic records the degeneration of a pair of distinct roots in $\P^1$ to a double point.

Let $p:\mathcal U\to\mathcal H(V_5)$ be the universal family and $e:\mathcal U\to V_5$ its evaluation morphism. The latter is finite of degree three. Over $\Omega_3$, $\Omega_2$ and $\Omega_1$, its fibres parametrise, respectively, three ordinary lines, one ordinary and one special line, and a single special line, with total length three in each case. The special lines are precisely the tangent lines to the
closed orbit $\Omega_1$; their union is the boundary $\overline{\Omega_2}$, and distinct special lines are disjoint. We refer to \cite[Section~1.2]{duboulozCylindersMoriFiber2019} for these facts.

The conic also describes incidence between lines. For a line $\ell$, let
\[M_\ell= \overline{\{[m]\in\mathcal H(V_5)\mid 	m\ne\ell,\ m\cap\ell\ne\varnothing\}}.\]
Then $M_\ell$ is the polar line of $[\ell]$ with respect to $\mathcal C$; see \cite[Proposition~2.4]{duboulozCylindersMoriFiber2019}. In particular, it is tangent to $\mathcal C$ at $[\ell]$ when $\ell$ is special. When $\ell$ is ordinary, it meets $\mathcal C$ at two distinct points, corresponding to the two special lines meeting $\ell$. The intersection map
\[ \rho_\ell:M_\ell\longrightarrow\ell, \qquad [m]\longmapsto m\cap\ell, \]
is then a double cover, ramified precisely at these two points. Indeed, the two lines other than $\ell$ through a general point of $\ell$ coalesce exactly when one of them becomes special.

Over $\R$, both Hilbert schemes remain projective planes:
\[\mathcal H(X_5)\simeq\P^2_{\R}, \qquad \mathcal H(Y_5)\simeq\P^2_{\R}.\]
This is a special case of \cite[Lemma~2.1]{duboulozCylindersMoriFiber2019}. Their distinguished conics, however, have different real structures. In suitable real coordinates they are
\[\mathcal C_{X_5}=\{b^2-4ac=0\},\qquad \mathcal C_{Y_5}=\{u^2+v^2+w^2=0\}. \]
The first is the discriminant conic of real binary quadratics. For the second, a real point of the Hilbert scheme corresponds to an antipodal pair on the Riemann sphere, which can never be a double point. Consequently, $X_5$ contains real special lines, while $Y_5$ contains none; both varieties nevertheless contain real ordinary lines.

The action of $\operatorname{PGL}_2(\R)$ on $\mathcal H(X_5)(\R)$ has exactly three orbits:
\begin{enumerate}
	\item $b^2-4ac=0$: the special lines, let $\ell_{\mathrm{sp}}$ be any of them;
	\item $b^2-4ac>0$: the ordinary lines meeting two real special lines, let $\ell_{\mathrm{ext}}$ be any of them;
	\item $b^2-4ac<0$: the ordinary lines whose two incident special lines are non-real and complex conjugate, let $\ell_{\mathrm{int}}$ be any of them.
\end{enumerate}
We shall call the latter two types \emph{exterior} and \emph{interior}, respectively. These names refer to the position of $[\ell]$ relative to the oval $\mathcal C_{X_5}(\R)$: its interior is an open disk, whereas its exterior is an open M\"obius band. 

On $Y_5$, the group $\operatorname{SO}_3(\R)$ acts transitively on $\mathcal H(Y_5)(\R)\approx\RP^2$, and the two special lines meeting any real line are complex conjugate.

Finally, every real point of either threefold lies on at least one real line: conjugation acts on the length-three fibre of $e$, and its non-real points occur in pairs of equal multiplicity.

\subsubsection{Projection from a line and the associated Sarkisov link.}\label{sec:ProjectionFromLineV5}

We now turn to the birational construction associated with a real line $\ell$ on $Z$, where $Z$ denotes either $X_5$ or $Y_5$. Projection from $\ell\subset Z\subset\P^6_{\R}$ is resolved by its blow-up $\alpha:\widetilde Z=\operatorname{Bl}_\ell(Z)\longrightarrow Z$.
Writing $S'$ for the exceptional divisor and $H^*=\alpha^*H$, the linear system $|H^*-S'|$ defines the second morphism from $\widetilde Z$ in the diagram 
\[ \xymatrix{ & S' \ar[dl] \ar@{^{(}->}[rr]& & \widetilde Z \ar[dl]_{\Bl_\ell\, = \, \alpha} \ar[dr]^{\beta\,=\, \Bl_\Gamma} & & T_\ell' \ar[dr] \ar@{_{(}->}[ll]& \\
	\ell \ar@{^{(}->}[r] & T_\ell \ar@{^{(}->}[r] & Z \ar@{-->}[rr]^{|H-\ell|} & & Q & S \ar@{_{(}->}[l] & \ar@{_{(}->}[l]\Gamma} \]
Here $Q$ is a smooth quadric threefold, and $\beta$ is the blow-up of a smooth twisted cubic $\Gamma\subset Q$. This is the classical Sarkisov link of type~II obtained by projection from a line; see \cite[Section~1.1.1]{duboulozCylindersMoriFiber2019} and \cite[Proposition~2.1.6(2)]{takagiScorzaQuarticsTrigonal2008}.

To describe its exceptional loci, let $T_\ell\subset Z$ be the surface swept out by the lines meeting $\ell$. It is a hyperplane section singular along $\ell$, with multiplicity two there. Its strict transform $T_\ell'$ is the exceptional divisor of $\beta$, and 
\[ \beta^*\O_Q(1)=H^*-S',\qquad T_\ell'\sim H^*-2S',\qquad H^* \sim 2\beta^*\O_Q(1)-T_\ell'.\]
The last equality describes the inverse rational map: it is given by quadrics in $\P^4$ containing $\Gamma$. Moreover, 
\[S:=\beta(S')=Q\cap\langle\Gamma\rangle\]
is a hyperplane section containing $\Gamma$, and the link restricts to an isomorphism $Z\setminus T_\ell\simeq Q\setminus S$.

In all cases the target is the real quadric of signature $(3,2)$. Indeed, the fibres of $S'\to\ell$ map to real lines on $Q$, so its defining quadratic form has Witt index two. In particular, projection from a real ordinary line recovers the rationality of both $X_5$ and $Y_5$ over $\R$.
For an ordinary line, $S'\simeq\P^1\times\P^1$, and $\beta|_{S'}$ identifies it with the smooth quadric surface $S=Q_0$, with the notations of \Cref{prop:ClassificationCurvesQ32} and using an automorphism of $Q$ if necessary. 

In addition, the three choices of line on $X_5$ therefore give the following real pictures.
\begin{enumerate}
	\item \emph{A special line.}
	Here $S'\simeq\mathbb F_2$, and $\beta|_{S'}$ contracts its negative section onto the vertex of the quadric cone $S$. Thus $S$ is a tangent hyperplane section of the smooth threefold $Q$. The twisted cubic $\Gamma$ passes through the vertex. Up to isomorphism of $Q$, $\Gamma = C_{3,s}$.
	\item \emph{An exterior ordinary line.}
	The double cover $\rho_\ell$ has two real ramification points, corresponding to the two real special lines meeting $\ell$. In suitable coordinates on its source and target it is $[s:t] \longmapsto [s^2:t^2]$. On real points its image is a closed interval in $\ell(\R)\approx \S^1$. Up to isomorphism in $Q$, $\Gamma = C_{3,0}$.
	\item \emph{An interior ordinary line.}
	The two ramification points of $\rho_\ell$ are complex conjugate. A real normal form is $[s:t] \longmapsto [s^2-t^2:2st]$. Its restriction to real points is an unramified double covering of circles. Up to isomorphism in $Q$, $\Gamma = C_{3,-2}$.
\end{enumerate}

For $Y_5$, every real line is ordinary, so $S'\simeq\mathbb F_0$, $T_\ell'\simeq\mathbb F_1$, and $\rho_\ell$ again induces an unramified double covering of real circles. This agrees with the interior case on $X_5$ at the level of the double cover, but the full embedded pairs $(Q,\Gamma)$ are not isomorphic over $\R$: otherwise the inverse construction would identify $X_5$ and $Y_5$. Up to isomorphism of $Q$, $\Gamma = C_{3,2}$.

\subsection{Computation of the topology.}

\begin{thm}\label{Th:TopologyOfV5}
	Let $X_5$ and $Y_5$ be the real forms of $V_5$ corresponding to the trivial and nontrivial real structures, respectively. Then
	\[X_5(\R)\approx\mathbb{RP}^3,\]
	whereas
	\[Y_5(\R) \approx \faktor{\operatorname{SO}_3(\R)}{O} \approx \faktor{\S^3}{2O} \approx M\big(\S^2;0;-1;(2,1);(3,1);(4,1)\big).\]
	Here $O\simeq\mathfrak S_4$ is the rotational octahedral group,	of order $24$, and $2O$ is the binary octahedral group, of order $48$, obtained as its inverse image under the double covering
	$\operatorname{SU}_2\to\operatorname{SO}_3(\R)$.
\end{thm}

\begin{proof}
	The identification of $Y_5(\R)$ with $\operatorname{SO}_3(\R)/O\approx\S^3/2O$ was established in
	\Cref{cor:TopologyTwistedV5}. Its displayed Seifert description is the classical one induced by the Hopf fibration, with base orbifold $\S^2(2,3,4)$. This Seifert fibration is unique up to fibre-preserving homeomorphism; the remark below gives	a geometric description of it.
	
	To determine $X_5(\R)$, we use the Sarkisov link of \Cref{sec:ProjectionFromLineV5}. Start with
	\[ C_{3,0}\subset Q_0\subset\pazocal Q_{3,2} \]
	in the notation of \Cref{prop:ClassificationCurvesQ32}, and let $\beta:\widetilde X_5\to\pazocal Q_{3,2}$ be the blow-up of $C_{3,0}$. The other morphism in the link, $\alpha:\widetilde{X_5}\to X_5$, contracts the strict transform $S'=Q_0'$ onto an ordinary line $\ell\subset X_5$.
	
	First, we identify the ruling contracted by $\alpha$. Since $C_{3,0}$ is a divisor on $Q_0$, the restriction $\beta|_{Q_0'}$ identifies $Q_0'$ with $Q_0\simeq\P^1\times\P^1$. Under this identification,
	\[\begin{aligned}
		\NN_{Q_0'/\widetilde{X_5}}	&\simeq \NN_{Q_0/\pazocal Q_{3,2}} \otimes\O_{Q_0}(-C_{3,0})\\
		&\simeq \O_{\P^1\times\P^1}(1,1) \otimes\O_{\P^1\times\P^1}(-1,-2)\\
		&\simeq\O_{\P^1\times\P^1}(0,-1).
	\end{aligned}\]
	An exceptional fibre has normal degree $-1$ in this line bundle. Thus $\alpha|_{Q_0'}$ is the first projection: its fibres are the curves $\{p\}\times\P^1$.
	
	Set $T=Q_0'(\R)$. The topological description in the proof of \Cref{Th:BlowUpQuadricThreefold32AlongRatCurves} gives an identification of pairs
	\[ \big(\widetilde{X_5}(\R),T\big) \approx \big(\RP^2\times\S^1, \RP^1\times\S^1\big). \]
	In particular, $T$ is a one-sided torus. A closed tubular neighbourhood $N$ of $T$ is the product of a M\"obius band with $\S^1$, and its complement
	\[
	V:=\widetilde{X_5}(\R)\setminus\operatorname{Int}(N)
	\approx\mathbb D^2\times\S^1
	\]
	is a solid torus. Blowing down $T$ replaces $N$ by a tubular neighbourhood $V'$ of $\ell(\R)$. Since $\ell$ is ordinary, this neighbourhood is also a solid torus. Hence $X_5(\R)=V\cup_{\partial V}V'$ is a lens space and it remains to determine the meridian of $V'$ on $\partial V$.
	
	Let $a,b\in H_1(T;\mathbb Z)$ be the classes of the two product circles $\RP^1\times\{*\}$ and $\{*\}\times\S^1$. In these product coordinates, the ruling	contracted by $\alpha$ has class
	\[ c=a+b,\]
	after choosing orientations. This is the slope identification between the algebraic ruling and the topological product structure used in the construction above.
	
	The boundary projection
	\[q:\partial N=\partial V\longrightarrow T\]
	is the connected double covering associated with the normal	bundle of $T$. Choose the product meridian $\mu$ and longitude $\lambda$ of $V$. Since the boundary of a M\"obius band covers
	its core twice, these classes satisfy
	\[q_*(\mu)=2a, \qquad q_*(\lambda)=b.\]
	A meridian $m$ of the solid torus $V'$ projects twice onto a contracted circle: locally, this is the usual map from the boundary of a normal disk to its projective line of directions.
	Consequently,
	\[ q_*(m)=2c=2a+2b, \qquad\text{and hence}\qquad m=\mu+2\lambda.\]
	The two meridians therefore have geometric intersection number $2$. The resulting genus-one Heegaard splitting is that of $L(2,1)$, and we conclude that
	\[ X_5(\R)\approx L(2,1)\approx\RP^3.\]
\end{proof}

\begin{rem}
	In this computation, the choice of $C_{3,0}$ is a matter of convenience. We could just as well have chosen $C_{3,-2}$ or $C_{3,s}$.
\end{rem}

\begin{rem}\label{rem:SarkisovLinkTwistedV5}
	The Sarkisov link also allows us to study the real locus of $\Bl_{C_{3,2}}\pazocal Q_{3,2}$. For this purpose, we describe a Seifert fibration on $Y_5(\R)$ for which the real locus of
	a line is an exceptional fibre of multiplicity $4$.
	
	Using the description as a quasi-homogeneous space, $Y_5(\R)\approx \faktor{\operatorname{SU}(2)}{2O}$. Consider the maximal torus
	\[D = \left\{ d_\theta = \begin{pmatrix} 
		e^{\i\theta} & 0\\0 & e^{-\i\theta}
	\end{pmatrix} \;\middle|\; \theta\in\R \right\}	\subset \operatorname{SU}(2).\]
	Left multiplication by $D$ defines a circle action on $\faktor{\operatorname{SU}(2)}{2O}$.
	Its kernel is the central subgroup $\{\pm I\}$, so the effective acting circle is
	\[ T:=\faktor{D}{\{\pm I\}} \approx \S^1.\]
	
	To describe the orbit space, use the Hopf fibration
	\[\operatorname{SU}(2) \longrightarrow \leftfaktor{D}{\operatorname{SU}(2)} \approx \S^2.\]
	Right multiplication by $2O$ preserves its fibres, and the induced action on $\S^2$ factors through $\faktor{2O}{\{\pm I\}}=O$. Thus the Hopf fibration descends to
	\[ \pi\colon Y_5(\R)\approx \faktor{\operatorname{SU}(2)}{2O} \longrightarrow \faktor{D\backslash\operatorname{SU}(2)}{2O} \approx \faktor{\S^2}{O}. \]
	As an orbifold, $\faktor{\S^2}{O}$ is a sphere with three cone points, of orders $2$, $3$ and $4$.
		
	More explicitly, the stabiliser of $g2O\in \faktor{\operatorname{SU}(2)}{2O}$ under the effective circle action is
	\[ T_{g2O} = \faktor{\bigl(D\cap g2Og^{-1}\bigr)}{\{\pm I\}}. \]
	These stabilisers are finite cyclic groups. Suppose $T_{g2O}$ has order $m$. Choose a small disk transverse to the $T$-orbit. By the slice theorem, a neighbourhood of the orbit is 	\[N_g \approx \faktor{\S^1\times \mathbb{D}^2}{\Z_m}\]
	where the standard generator of $\Z_m$ acts by $(u,z)\mapsto (e^{2\pi \i/m}u,e^{2\pi \i/m}z)$ up to signs.
	
	We now identify the fibre of multiplicity $4$ that is relevant to the Sarkisov link. Recall that
	\[\varphi=xy(x^4+y^4).\]
	If $D_{\C}$ denotes the complexification of $D$, then
	\[ \overline{D_{\C}\cdot[\varphi]} = \P\bigl(\langle x^5y,xy^5\rangle\bigr) =:\ell_{\C}.\]
	This line is preserved by the real structure defining $Y_5$, and therefore descends to a real line $\ell\subset Y_5$. Its real locus is exactly	the orbit $D\cdot[\varphi]$.	
	Indeed, the action is given explicitly by
	\[ d_\theta\cdot\varphi = e^{4\i\theta}x^5y+e^{-4\i\theta}xy^5.\]
	In addition,
	\[d_\theta\cdot[\varphi]=[\varphi] \quad\Longleftrightarrow\quad e^{8\i\theta}=1.\]
	The stabiliser in $D$ therefore has order $8$, and its image in $T$ has order $4$. Thus $\ell(\R)=D\cdot[\varphi]$ is an exceptional Seifert fibre of multiplicity $4$.
	
	The desired topology is therefore obtained by performing the topological blow-up of $Y_5(\R)$ along this exceptional fibre. The description in \Cref{Th:TopBlowUpSeifertFibrationAlongFibre} applies with multiplicity $4$.
\end{rem}

\begin{rem}
	The multiplicities $2$, $3$ and $4$ also admit a geometric interpretation in terms of octahedral symmetry. A fundamental domain for
	\[Y_5(\R)\approx \faktor{\operatorname{SO}_3(\R)}{O} \approx \faktor{\S^3}{2O}\]
	can be represented by a truncated cube, with six octagonal faces and eight triangular faces. Opposite octagonal faces are identified with a twist of $\pi/4$, and opposite triangular faces with a twist of $\pi/3$, using the standard face-pairings induced by $2O$.
	
	The three types of rotational axes of the cube correspond to the three exceptional multiplicities:
	\[\begin{array}{c|c}
		\text{Axis through opposite}
		& \text{Order of the cyclic rotation subgroup}\\
		\hline
		\text{face centres}  & 4\\
		\text{vertices}      & 3\\
		\text{edge midpoints}& 2.
	\end{array}\]
	In the truncated-cube picture, these directions pass through the centres of the octagonal faces, the centres of the triangular faces, and the midpoints of the edges shared by two octagonal faces, respectively.
	
	More precisely, let $A$ be one of these axes, and let $T_A\subset\operatorname{SO}_3(\R)$ be the circle subgroup of rotations about $A$. If the rotations in $O$ about $A$ form a cyclic group of order $m$, then the orbit of the identity coset is
	\[ T_A\cdot[I] \simeq T_A/(T_A\cap O), \qquad T_A\cap O \simeq \Z_m. \]
	Choosing $T_A$ as the acting circle makes this orbit an	exceptional Seifert fibre of multiplicity $m$. This realises the multiplicities $4$, $3$ and $2$ directly from the three types of symmetry axes. Changing the axis changes the circle action, but yields an equivalent
	Seifert fibration.
\end{rem}

\printbibliography[title={References}]

@article{abbanKstabilityPointlessPezzo2024,
  title = {K-Stability of Pointless {{del Pezzo}} Surfaces and {{Fano}} 3-Folds},
  author = {Abban, Hamid and Cheltsov, Ivan and Kishimoto, Takashi and Mangolte, Frédéric},
  date = {2026},
  journaltitle = {Forum of Mathematics, Sigma},
  shortjournal = {Forum Math. Sigma},
  volume = {14},
  eprint = {2411.00767},
  eprinttype = {arXiv},
  pages = {e142},
  doi = {10.1017/fms.2026.10266},
  url = {https://doi.org/10.1017/fms.2026.10266}
}

@article{bjorklundRealAlgebraicKnots2011,
  title = {Real Algebraic Knots of Low Degree},
  author = {Björklund, Johan},
  date = {2011},
  journaltitle = {Journal of Knot Theory and Its Ramifications},
  shortjournal = {J. Knot Theory Its Ramif.},
  volume = {20},
  number = {9},
  pages = {1285--1309},
  issn = {0218-2165},
  doi = {10.1142/S0218216511009248},
  url = {https://doi.org/10.1142/S0218216511009248}
}

@article{cheltsovKstableFanoThreefolds2022,
  title = {K-Stable {{Fano}} Threefolds of Rank 2 and Degree 30},
  author = {Cheltsov, Ivan and Park, Jihun},
  date = {2022},
  journaltitle = {European Journal of Mathematics},
  shortjournal = {Eur. J. Math.},
  volume = {8},
  number = {3},
  pages = {834--852},
  issn = {2199-6768},
  doi = {10.1007/s40879-022-00569-x},
  url = {https://doi.org/10.1007/s40879-022-00569-x}
}

@article{comessattiSullaConnessioneSuperficie1914,
  title = {Sulla Connessione Delle Superficie Razionali Reali},
  author = {Comessatti, Annibale},
  date = {1914},
  journaltitle = {Annali di Matematica Pura ed Applicata},
  shortjournal = {Ann. Mat. Pura Ed Appl.},
  volume = {23},
  number = {1},
  pages = {215--284},
  issn = {0373-3114},
  doi = {10.1007/BF02419577},
  url = {https://doi.org/10.1007/BF02419577}
}

@incollection{demazureSurfacesPezzoII1980,
  title = {Surfaces de Del Pezzo. II. Éclater n points dans P\textsuperscript{2}},
  booktitle = {Séminaire sur les singularités des surfaces},
  author = {Demazure, Michel},
  editor = {Demazure, Michel and Pinkham, Henry Charles and Teissier, Bernard},
  date = {1980},
  series = {Lecture notes in mathematics},
  volume = {777},
  pages = {23--35},
  publisher = {Springer},
  location = {Berlin},
  doi = {10.1007/BFb0085876},
  url = {https://doi.org/10.1007/BFb0085876},
  isbn = {978-3-540-38586-8},
  langid = {french}
}

@article{duboulozCylindersMoriFiber2019,
  title = {Cylinders in {{Mori}} Fiber Spaces: Forms of the Quintic {{del Pezzo}} Threefold},
  author = {Dubouloz, Adrien and Kishimoto, Takashi},
  date = {2019},
  journaltitle = {Annales de l'Institut Fourier},
  shortjournal = {Ann. Inst. Fourier},
  volume = {69},
  number = {6},
  pages = {2377--2393},
  doi = {10.5802/aif.3297},
  url = {https://doi.org/10.5802/aif.3297}
}

@article{eisenbudNormalBundlesSmooth1981,
  title = {On the Normal Bundles of Smooth Rational Space Curves},
  author = {Eisenbud, David and family=Ven, given=Antonius, prefix=van de, useprefix=true},
  date = {1981},
  journaltitle = {Mathematische Annalen},
  shortjournal = {Math. Ann.},
  volume = {256},
  number = {4},
  pages = {453--463},
  issn = {1432-1807},
  doi = {10.1007/BF01450541},
  url = {https://doi.org/10.1007/BF01450541}
}

@article{epsteinPeriodicFlowsThreeManifolds1972,
  title = {Periodic Flows on Three-Manifolds},
  author = {Epstein, D. B. A.},
  date = {1972},
  journaltitle = {Annals of Mathematics},
  shortjournal = {Ann. Math.},
  volume = {95},
  number = {1},
  eprint = {1970854},
  eprinttype = {JSTOR},
  pages = {66--82},
  issn = {0003486X, 19398980},
  doi = {10.2307/1970854},
  url = {https://doi.org/10.2307/1970854}
}

@article{geigesSeifertFibrationsLens2017,
  title = {Seifert Fibrations of Lens Spaces},
  author = {Geiges, Hansjörg and Lange, Christian},
  date = {2018},
  journaltitle = {Abhandlungen aus dem Mathematischen Seminar der Universität Hamburg},
  shortjournal = {Abh. Aus Dem Math. Semin. Univ. Hambg.},
  volume = {88},
  number = {1},
  eprint = {1608.06844},
  eprinttype = {arXiv},
  eprintclass = {math.GT},
  pages = {1--22},
  doi = {10.1007/s12188-017-0188-z},
  url = {https://doi.org/10.1007/s12188-017-0188-z}
}

@incollection{KhBourbaki,
  title = {Variétés de {{Fano}} Réelles [d'après {{C}}. {{Viterbo}}]},
  booktitle = {Séminaire {{Bourbaki}}, Vol. 1999/2000, Exposés 865–879},
  author = {Kharlamov, Viatcheslav},
  date = {2002},
  series = {Astérisque},
  pages = {189--206},
  publisher = {Société Mathématique de France},
  location = {Paris},
  url = {https://www.numdam.org/item/SB_1999-2000__42__189_0/},
  eid = {276},
  note = {Exposé no. 872}
}

@online{kineiderClassificationSeifertFibrations2026,
  type = {Survey},
  title = {Classification of {{Seifert}} Fibrations through Orbifold {{Euler}} Classes.},
  author = {Kineider, Guillaume},
  date = {2026},
  url = {https://www.i2m.univ-amu.fr/perso/guillaume.kineider/seifert-circle-bundles.pdf},
  pubstate = {prepublished}
}

@article{kollarRealAlgebraicThreefolds1998,
  title = {Real Algebraic Threefolds {{I}}. {{Terminal}} Singularities},
  author = {Kollár, János},
  date = {1998},
  journaltitle = {Collectanea Mathematica},
  shortjournal = {Collect. Math.},
  volume = {49},
  number = {2--3},
  pages = {335--360},
  issn = {0010-0757},
  url = {https://eudml.org/doc/41339}
}

@article{kollarRealAlgebraicThreefolds1999,
  title = {Real Algebraic Threefolds {{II}}. {{Minimal}} Model Program},
  author = {Kollár, János},
  date = {1999},
  journaltitle = {Journal of the American Mathematical Society},
  shortjournal = {J. Am. Math. Soc.},
  volume = {12},
  number = {1},
  pages = {33--83},
  issn = {0894-0347},
  doi = {10.1090/s0894-0347-99-00286-6},
  url = {https://doi.org/10.1090/s0894-0347-99-00286-6}
}

@article{kollarRealAlgebraicThreefolds1999a,
  title = {Real Algebraic Threefolds {{III}}. {{Conic}} Bundles},
  author = {Kollár, János},
  date = {1999},
  journaltitle = {Journal of Mathematical Sciences},
  shortjournal = {J. Math. Sci.},
  volume = {94},
  number = {1},
  pages = {996--1020},
  issn = {1573-8795},
  doi = {10.1007/BF02367244},
  url = {https://doi.org/10.1007/BF02367244}
}

@incollection{kollarRealAlgebraicThreefolds2016,
  title = {Real Algebraic Threefolds {{IV}}. {{Del Pezzo}} Fibrations},
  booktitle = {Complex Analysis and Algebraic Geometry: A Volume in Memory of {{Michael Schneider}}},
  author = {Kollár, János},
  editor = {Peternell, Thomas and Schreyer, Frank-Olaf},
  date = {2000},
  eprint = {math/9901022},
  eprinttype = {arXiv},
  pages = {317--346},
  publisher = {Walter de Gruyter},
  location = {Berlin},
  doi = {10.1515/9783110806090-014},
  url = {https://doi.org/10.1515/9783110806090-014},
  isbn = {978-3-11-080609-0}
}

@incollection{kollarTopologyRealComplex,
  title = {The Topology of Real and Complex Algebraic Varieties},
  booktitle = {Taniguchi Conference on Mathematics Nara '98},
  author = {Kollár, János},
  date = {2001},
  series = {Advanced Studies in Pure Mathematics},
  volume = {31},
  pages = {127--145},
  publisher = {Mathematical Society of Japan},
  location = {Tokyo},
  doi = {10.2969/aspm/03110127},
  url = {https://doi.org/10.2969/aspm/03110127}
}

@article{krasnovIntersectionsTwoReal2018,
  title = {On Intersections of Two Real Quadrics},
  author = {Krasnov, V. A.},
  date = {2018},
  journaltitle = {Izvestiya: Mathematics},
  shortjournal = {Izv. Math.},
  volume = {82},
  number = {1},
  pages = {91--139},
  issn = {1064-5632, 1468-4810},
  doi = {10.1070/IM8535},
  url = {https://doi.org/10.1070/IM8535}
}

@article{krasnovRigidIsotopyClassification2006,
  title = {Rigid Isotopy Classification of Real Three-Dimensional Cubics},
  author = {Krasnov, Vyacheslav A.},
  date = {2006},
  journaltitle = {Izvestiya: Mathematics},
  shortjournal = {Izv. Math.},
  volume = {70},
  number = {4},
  pages = {731--768},
  issn = {1064-5632},
  doi = {10.1070/IM2006v070n04ABEH002326},
  url = {https://doi.org/10.1070/IM2006v070n04ABEH002326}
}

@incollection{KS04,
  title = {Real Fano 3-Folds of Type V$_{22}$},
  booktitle = {The {{Fano}} Conference},
  author = {Kollár, János and Schreyer, Frank-Olaf},
  date = {2004},
  pages = {515--531},
  publisher = {Università di Torino},
  location = {Turin}
}

@article{mangolteUniruledSixmanifoldsContain2012,
  title = {Do Uniruled Six-Manifolds Contain {{Sol Lagrangian}} Submanifolds?},
  author = {Mangolte, Frédéric and Welschinger, Jean-Yves},
  date = {2012},
  journaltitle = {International Mathematics Research Notices},
  shortjournal = {Int. Math. Res. Not.},
  volume = {2012},
  number = {7},
  eprint = {1001.2927},
  eprinttype = {arXiv},
  eprintclass = {math.SG},
  pages = {1569--1602},
  issn = {1687-0247, 1073-7928},
  doi = {10.1093/imrn/rnr063},
  url = {http://arxiv.org/abs/1001.2927},
  urldate = {2026-11-13},
  note = {Comment: 23 pages, 2 figures}
}

@article{mikhalkinBlowupEquivalenceSmooth1997,
  title = {Blowup Equivalence of Smooth Closed Manifolds},
  author = {Mikhalkin, Grigory},
  date = {1997},
  journaltitle = {Topology},
  volume = {36},
  number = {1},
  pages = {287--299},
  issn = {00409383},
  doi = {10.1016/0040-9383(95)00062-3},
  url = {https://doi.org/10.1016/0040-9383(95)00062-3}
}

@article{mikhalkinMaximallyWrithedReal2019a,
  title = {Maximally Writhed Real Algebraic Links},
  author = {Mikhalkin, Grigory and Orevkov, Stepan},
  date = {2019},
  journaltitle = {Inventiones mathematicae},
  shortjournal = {Invent. Math.},
  volume = {216},
  number = {1},
  pages = {125--152},
  issn = {1432-1297},
  doi = {10.1007/s00222-018-0844-7},
  url = {https://doi.org/10.1007/s00222-018-0844-7}
}

@article{mikhalkinOsculatingFramingReal2019,
  title = {On Osculating Framing of Real Algebraic Links},
  author = {Mikhalkin, Grigory and Orevkov, Stepan},
  date = {2019},
  journaltitle = {Arnold Mathematical Journal},
  shortjournal = {Arnold Math. J.},
  volume = {5},
  number = {4},
  pages = {393--399},
  doi = {10.1007/s40598-019-00120-z},
  url = {https://doi.org/10.1007/s40598-019-00120-z}
}

@article{mikhalkinRealAlgebraicKnots2016,
  title = {Real Algebraic Knots and Links of Small Degree},
  author = {Mikhalkin, Grigory and Orevkov, Stepan},
  date = {2016},
  journaltitle = {Journal of Knot Theory and Its Ramifications},
  shortjournal = {J. Knot Theory Its Ramif.},
  volume = {25},
  number = {12},
  eprint = {1608.03648},
  eprinttype = {arXiv},
  pages = {1642010},
  doi = {10.1142/S0218216516420104},
  url = {https://doi.org/10.1142/S0218216516420104}
}

@article{mikhalkinRigidIsotopyMaximally2021,
  title = {Rigid Isotopy of Maximally Writhed Links},
  author = {Mikhalkin, Grigory and Orevkov, Stepan},
  date = {2021},
  journaltitle = {Algebraic Geometry},
  shortjournal = {Algebr. Geom.},
  volume = {8},
  number = {3},
  pages = {268--285},
  issn = {22142584},
  doi = {10.14231/AG-2021-006},
  url = {https://doi.org/10.14231/AG-2021-006}
}

@article{mikhalkinTopologyMaximallyWrithed2018,
  title = {Topology of Maximally Writhed Real Algebraic Knots},
  author = {Mikhalkin, G. B. and family=Orevkov, given=S. Yu., given-i=S{{Yu}}},
  date = {2018},
  journaltitle = {Doklady Mathematics},
  shortjournal = {Dokl. Math.},
  volume = {97},
  number = {1},
  pages = {28--31},
  issn = {1064-5624, 1531-8362},
  doi = {10.1134/S106456241801009X},
  url = {https://doi.org/10.1134/S106456241801009X}
}

@incollection{MukaiUmemura,
  title = {Minimal Rational Threefolds},
  booktitle = {Algebraic Geometry (Tokyo/Kyoto, 1982)},
  author = {Mukai, Shigeru and Umemura, Hiroshi},
  editor = {Raynaud, Michel and Shioda, Tetsuji},
  date = {1983},
  series = {Lecture Notes in Mathematics},
  volume = {1016},
  pages = {490--518},
  publisher = {Springer},
  location = {Berlin},
  doi = {10.1007/BFb0099976},
  url = {https://doi.org/10.1007/BFb0099976},
  isbn = {978-3-540-38676-6}
}

@article{scottGeometries3Manifolds1983,
  title = {The Geometries of 3-Manifolds},
  author = {Scott, Peter},
  date = {1983},
  journaltitle = {Bulletin of the London Mathematical Society},
  shortjournal = {Bull. Lond. Math. Soc.},
  volume = {15},
  number = {5},
  pages = {401--487},
  issn = {1469-2120},
  doi = {10.1112/blms/15.5.401},
  url = {https://doi.org/10.1112/blms/15.5.401}
}

@article{segalTopologySpacesRational1979,
  title = {The Topology of Spaces of Rational Functions},
  author = {Segal, Graeme},
  date = {1979},
  journaltitle = {Acta Mathematica},
  shortjournal = {Acta Math.},
  volume = {143},
  pages = {39--72},
  issn = {0001-5962, 1871-2509},
  doi = {10.1007/BF02392088},
  url = {https://doi.org/10.1007/BF02392088}
}

@misc{takagiScorzaQuarticsTrigonal2008,
  title = {Scorza Quartics of Trigonal Spin Curves and Their Varieties of Power Sums},
  author = {Takagi, Hiromichi and Zucconi, Francesco},
  date = {2008},
  eprint = {0801.1760},
  eprinttype = {arXiv},
  eprintclass = {math.AG},
  doi = {10.48550/arXiv.0801.1760},
  url = {https://arxiv.org/abs/0801.1760},
  howpublished = {arXiv:0801.1760}
}

@article{terpereauRealFormsMori2024,
  title = {Real Forms of {{Mori}} Fiber Spaces with Many Symmetries},
  author = {Terpereau, Ronan and Zimmermann, Susanna},
  date = {2026},
  journaltitle = {Annali della Scuola Normale Superiore di Pisa, Classe di Scienze},
  shortjournal = {Ann. Della Scuola Norm. Super. Pisa Cl. Sci.},
  eprint = {2403.14493},
  eprinttype = {arXiv},
  eprintclass = {math.AG},
  doi = {10.2422/2036-2145.202405_021},
  url = {https://doi.org/10.2422/2036-2145.202405_021},
  note = {Published online 22 January 2026}
}

@incollection{viroEncomplexingWrithe2000,
  title = {Encomplexing the Writhe},
  booktitle = {Topology, Ergodic Theory, Real Algebraic Geometry: {{Rokhlin}}'s Memorial},
  author = {Viro, Oleg},
  editor = {Turaev, Vladimir and Vershik, Anatoly},
  date = {2001},
  series = {American Mathematical Society Translations, Series 2},
  volume = {202},
  eprint = {math/0005162},
  eprinttype = {arXiv},
  pages = {241--256},
  publisher = {American Mathematical Society / ResearchGate},
  location = {Providence, RI},
  doi = {10.1090/trans2/202/17},
  url = {https://doi.org/10.1090/trans2/202/17}
}

@misc{weiWuClassesTopological2011,
  title = {On the {{Wu}} Classes of the Topological Blow Up},
  author = {Wang, Wei},
  date = {2011},
  eprint = {1107.1065},
  eprinttype = {arXiv},
  eprintclass = {math.AT},
  doi = {10.48550/arXiv.1107.1065},
  url = {https://arxiv.org/abs/1107.1065},
  howpublished = {arXiv:1107.1065}
}

@article{yeConstrainedKnotsLens2023,
  title = {Constrained Knots in Lens Spaces},
  author = {Ye, Fan},
  date = {2023},
  journaltitle = {Algebraic \& Geometric Topology},
  shortjournal = {Algebr. Geom. Topol.},
  volume = {23},
  number = {3},
  eprint = {2007.04237},
  eprinttype = {arXiv},
  eprintclass = {math.GT},
  pages = {1097--1166},
  issn = {1472-2739, 1472-2747},
  doi = {10.2140/agt.2023.23.1097},
  url = {https://doi.org/10.2140/agt.2023.23.1097}
}
	
\end{document}